\documentclass[12pt]{article}
\usepackage[letterpaper]{geometry}

\usepackage{graphicx,color}
\usepackage{graphics} 
\usepackage{float}
\usepackage{multirow}
\usepackage{epstopdf}
\usepackage{subfigure}
\usepackage{rotating}
\usepackage{ctable}
\usepackage[square,sort,comma,numbers]{natbib}
\usepackage{amsthm}

\newtheorem{theorem}{Theorem}

\newcolumntype{C}{>{\centering\arraybackslash} m{0.3\textwidth} }
\newcolumntype{R}{>{\centering\arraybackslash} m{0.005\textwidth} }

\newcolumntype{N}{>{\centering\arraybackslash} m{0.2\textwidth} }
\newcolumntype{V}{>{\centering\arraybackslash} m{0.005\textwidth} }

\usepackage{amsmath,amssymb,amsfonts,amstext}
\usepackage{bbm}
\usepackage{mathrsfs}
\usepackage{dsfont}

\newcommand{\ind}{{\rm i_{\rm k}}}

\newcommand{\inv}[1]{\left[{\bf W^{(\rm #1)}}\right]^{\rm -1}}

\newcommand{\Nobs}{\textsc{n}_{{\rm obs}}}

\newcommand{\Nens}{\textsc{n}_{{\rm ens}}}

\newcommand{\Nproc}[1]{{\rm N_{proc}^{\rm #1}}}

\newcommand{\nc}[1]{{\rm C_{p}^{\rm #1}}}

\newcommand{\Nstate}{\textsc{n}_{{\rm state}}}

\newcommand{\BO}{{\mathcal O}}

\newcommand{\XA}{{\bf X}^{\rm A}}

\newcommand{\XB}{{\bf X}^{\rm B}}

\newcommand{\DXB}{{\bf \Delta X}^{\rm B}}

\newcommand{\xb}[1]{{\bf x}^{\rm B}_{{\rm #1}}}

\newcommand{\Xmean}{{\bf \overline{X}}^{\rm B}}

\newcommand{\xmean}{{\bf \overline{x}}^{\rm B}}

\newcommand{\xmeana}{{\bf \overline{x}}^{\rm A}}

\newcommand{\ONES}{\mathbf{1}}

\newcommand{\PB}{{\bf P}^{\rm B}}

\newcommand{\Q}{{\bf Q}}

\newcommand{\R}{{\bf R}}

\renewcommand{\S}{{\bf S}}

\newcommand{\Z}{{\bf Z}}

\newcommand{\sZ}[1]{{\bf Z}^{\left(\rm #1 \right)}}

\newcommand{\h}[1]{{\bf h}^{\left( \rm #1\right)}}

\newcommand{\z}[1]{{\bf z}_{{\rm #1}}}

\newcommand{\U}[1]{{\bf U}^{\left ( \rm #1 \right)}}

\renewcommand{\u}[2]{{\bf u}_{\rm #1}^{(\rm #2)}}

\newcommand{\um}[1]{{\bf u}_{\rm #1}}

\newcommand{\sz}[2]{{\bf z}_{\rm #1}^{(\rm #2)}}

\newcommand{\D}{{\bf D}}

\renewcommand{\d}[1]{{\bf d}_{\rm #1}}

\newcommand{\W}[1]{{{\bf W}^{\left ( {\rm #1} \right )}}}

\newcommand{\V}{{\bf V}}

\renewcommand{\v}[1]{{\bf v}_{{\rm #1}}}

\newcommand{\KAL}{{\bf K}}

\newcommand{\Y}{{\bf Y}}

\newcommand{\y}[1]{{\bf y}_{{\rm #1}}}

\newcommand{\xt}{{\bf x}^{\rm true}}

\newcommand{\Mo}{{\cal M}}

\newcommand{\Ho}{{\cal H}}

\newcommand{\Normal}{\cal N}

\newcommand{\Lo}{{\bf H}}

\renewcommand{\Re}{\mathbbm{R}}

\newcommand{\E}{{\bf \Upsilon}}

\newcommand{\G}[1]{\bf G^{(\rm #1)}}
\newcommand{\gj}[2]{\bf {g_{\rm #1}^{(\rm #2)}}}

\newcommand{\e}[1]{{\bf \upsilon}_{{\rm #1}}}

\newcommand{\A}{{\bf A}}
\renewcommand{\L}{{\bf L}}
\newcommand{\M}{{\bf M}}
\newcommand{\N}{{\bf N}}

\newcommand{\invM}[1]{\left( #1 \right )^{-1}}

\newcommand{\invS}[1]{ {#1^{-1}}}

\newcommand{\SMF}[4]{\invS{#1}-\invS{#1} \cdot #2  \cdot \invM{\invS{#3} + #4 \cdot #1 \cdot #2} \cdot #4 \cdot \invS{#1}}

\newcommand{\F}[1]{{\bf F}^{\left({\rm #1}\right)}}
\newcommand{\g}[1]{{\bf g}^{\left({\rm #1}\right)}}
\newcommand{\f}[2]{{{\bf f}_{{\rm #2}}^{\left({ \rm #1}\right)}}}

\newcommand{\x}{{\bf x}}

\newcommand{\fun}[1]{{\cal S}\left(#1 \right )}

\newcommand{\ifun}[3]{{\cal S}_{\star}\left(#1 ,#2, #3\right )}

\renewcommand{\k}{{\rm k}}

\newcommand{\level}[3]{
\text{ Level #1}:~~ \left\{ 
\begin{array}{l}
\h{#1} = \left( 1+\v{#1}^{ \bf T } \cdot \u{#1}{#2}\right)^{-1}\, \u{#1}{#2} \\
\sZ{#1} = \sZ{#2}-\h{#1} \cdot \left ( \v{#1}^{\bf T} \cdot \sZ{#2} \right ) \\
\u{i}{#1} = \u{i}{#2}-\h{#1} \cdot \left ( \v{#1}^{\bf T} \cdot \u{i}{#2} \right )\,, ~~ {\rm i}={\rm #3},\dots, \Nens \\
\end{array}
\right.
}

\newcommand{\T}[3]{{\rm T_{\rm #3} \left( #1,#2 \right )}}

\newcommand{\rv}{{\bf r}}

\newcommand{\diag}[2]{{\bf diag} \left ( #1_{1}, #1_{2}, \ldots ,  #1_{#2}\right )}

\newcommand{\Lorenz}{{\bf \frac{dx_{\rm i}}{dt}} = \begin{cases}
 \left ({{\bf x_{2}}-{\bf x_{\Nstate-1}}} \right ) \cdot {\bf x_{\Nstate}}-{\bf x_{1}}+{\rm F} & \text{ for $\rm i = 1$} \\
 \left ({{\bf x_{\rm i+1}}-{\bf x_{\rm i-2}}} \right ) \cdot {\bf x_{\rm i-1}}-{\bf x_{\rm i}}+{\rm F} & \text{ for $2 \le \rm  i 
\le \Nstate-1$} \\
\left ({{\bf x_{\rm 1}}-{\bf x_{\Nstate-2}}} \right ) \cdot {\bf x_{\Nstate-1}}-{\bf x_{\Nstate}}+{\rm F} & \text{ for $\rm i =  \Nstate$}
\end{cases}} 

\newcommand{\Bi}[2]{{\bf B_{\rm #1}^{\rm ( #2 )}}}

\begin{document}

\thispagestyle{empty}
\setcounter{page}{0}

\begin{Huge}
\begin{center}
Computational Science Laboratory Technical Report CSL-TR-00-2013\\
\today
\end{center}
\end{Huge}
\vfil
\begin{huge}
\begin{center}
Elias D. Nino, Adrian Sandu and Jeffrey Anderson
\end{center}
\end{huge}

\vfil
\begin{huge}
\begin{it}
\begin{center}
``An Efficient Implementation of the Ensemble Kalman Filter Based on an Iterative Sherman-Morrison Formula''
\end{center}
\end{it}
\end{huge}
\vfil
{\tiny
\textbf{Cite as:} Elias D. Nino-Ruiz, Adrian Sandu, Jeffrey Anderson, ``An efficient implementation of the ensemble Kalman filter based on an iterative Sherman–Morrison formula'',\textit{Statistics and Computing}, ISSN:0960-3174, PP: 1--17, Feb 2014.}

\begin{large}
\begin{center}
Computational Science Laboratory \\
Computer Science Department \\
Virginia Polytechnic Institute and State University \\
Blacksburg, VA 24060 \\
Phone: (540)-231-2193 \\
Fax: (540)-231-6075 \\ 
Email: {\rm enino@vt.edu},{\rm sandu@cs.vt.edu} \\
Web: {\rm http://csl.cs.vt.edu}
\end{center}
\end{large}

\vspace*{1cm}

\begin{center}
\begin{tabular}{c}
\includegraphics[width=1.2cm]{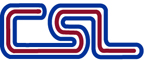}\\ 
.  \\
\includegraphics[width=1.2cm]{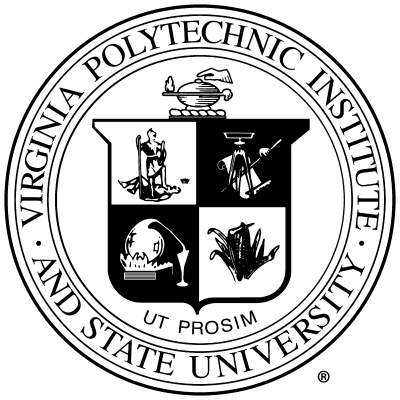} \\
\end{tabular}
\end{center}

\newpage

\title{An Efficient Implementation of the Ensemble Kalman Filter Based on an Iterative Sherman-Morrison Formula}

\author{Elias D. Nino$^\dagger$, Adrian Sandu$^\dagger$, and Jeffrey Anderson$^\ddagger$ \\
$^\dagger$Computational Science Laboratory, Department of Computer Science \\
Virginia Polytechnic Institute and State University \\
Blacksburg, VA 24060, USA \\
enino@vt.edu, sandu@cs.vt.edu \\
$^\ddagger$Data Assimilation Research Section \\
Institute for Mathematics Applied to Geosciences \\
National Center for Atmospheric Research \\
Boulder, CO 80307-3000 \\
jla@ucar.edu}

\date{}
\maketitle
\tableofcontents

\begin{abstract}
We present a practical implementation of the ensemble Kalman (EnKF) filter based on an iterative Sherman-Morrison formula. 
The new direct method exploits the special structure of the ensemble-estimated error covariance matrices  in order to efficiently solve the linear systems  involved in the analysis step of the EnKF. The computational complexity of the proposed implementation is equivalent to that of the best EnKF implementations available in the literature when the number of observations is much larger than the number of ensemble members. Even when this condition is not fulfilled, the proposed method is expected to perform well since it does not employ matrix decompositions. Moreover, the proposed method provides the best theoretical complexity when compared to generic formulations of matrix inversion based on the Sherman Morrison formula. The stability analysis of the proposed method is carried out and a pivoting strategy is discussed in order to reduce the accumulation of round-off errors without increasing the computational effort. A parallel implementation is discussed as well. Computational experiments carried out using the Lorenz 96 and then oceanic quasi-geostrophic models reveal that the proposed algorithm yields the same accuracy as other EnKF implementations, but is considerably faster.
\end{abstract}


\section{Introduction}
\label{intro}
The ensemble Kalman filter (EnKF) is a well-established, sequential Monte Carlo method to estimate the state and parameters of non-linear, large dynamical models \cite{Even09A} such as those found in atmospheric \cite{Ott04}, oil reservoir \cite{Even09B}, and oceanic \cite{Even08} simulations. The popularity of EnKF owes to its simple conceptual formulation and the relative ease implementation \cite{Evensen2009}. EnKF represents the error statistics by an ensemble of model states, and the evolution of error statistics is obtained implicitly via the time evolution of the ensemble during the forecast step. In the analysis step, information from the model and the measurements is combined in order to obtain an improved estimate of the true vector state. This process is repeated over the observed time period. In typical data assimilation applications, the dimension of state space (number of variables) ranges between $\BO(10^{7})$ and $\BO(10^{9})$, and 
the dimension of the observation space between $\BO(10^{5})$ and $\BO(10^{7})$. Consequently, the dimension of the linear systems solved
during the analysis step is very large, and the computational cost considerable. In order to address this challenge we propose an efficient implementation of the EnKF analysis step based on an iterative application of the Sherman-Morrison formula.

The paper is structured as follows. Section \ref{sec:enkf} discusses the conceptual formulation of the EnKF and several efficient implementations available in the literature. Section \ref{sec:enkfSherman} presents the novel implementation of the EnKF based on iterative Sherman-Morrison formula, in which the special structure of the measurements error covariance matrix is exploited. Computational cost and stability analyses are carried out for this approach, and pivoting and parallelization ideas are discussed. Section \ref{sec:results} reports numerical results of the proposed algorithm applied to the Lorenz 96 and quasi-geostrophic models. Conclusions are presented in Section \ref{sec:conclusions}. 

\section{Formulation of the EnKF}\label{sec:enkf}

EnKF consists of two steps: the forecast and the analysis.
An EnKF cycle starts with the matrix $\XB \in \Re^{\Nstate \times \Nens}$ whose columns ${\xb{i}} \in \Re^{\Nstate \times 1}$ form an ensemble of model states,
all corresponding to the same model time $t_{\rm current}$:
\begin{eqnarray} \nonumber
\displaystyle {\XB} &=& \left ( {\xb{1}},{\xb{2}}, \ldots , {\xb{\Nens}}\right ) \in \Re^{\Nstate \times \Nens}   \,. 
\end{eqnarray}
Typically $\xb{i}$ is an ensemble of model forecasts.
Here $\Nstate$ is the size of the model state vector, and $\Nens$ is the number of ensemble members. Each ensemble member $\xb{i}$ differs from the true state 
of the system $\xt \in \Re^{\Nstate \times 1}$, and we denote by ${ \rm \xi_i} \in \Re^{\Nstate \times 1}$ the corresponding error. The statistics of the ensemble of states is consistent with the background probability distribution.

The ensemble mean $\xmean \in \Re^{\Nstate \times 1}$ and the ensemble covariance matrix $\PB \in \Re^{\Nstate \times \Nstate}$ can be written as follows:
\begin{subequations}
\label{ensemble-background}
\begin{eqnarray} 
\label{ensemble-background-mean}
\displaystyle \xmean &=& \frac{1}{\Nens}  \sum_{i=1}^{\Nens} {\xb{i}} = \frac{1}{\Nens} \left ( \XB \cdot \ONES_{\Nstate \times 1} \right ) \in \Re^{\Nstate \times 1}, \\ \label{ensemble-background-mean-matrix}
\displaystyle \Xmean &=& \xmean \otimes {\ONES}_{\Nens \times 1}^{{\bf T}} \in \Re^{\Nstate \times \Nens}, \\ 
\label{ensemble-background-covariance}
{\PB} &=& \frac{1}{{\Nens}-1} \cdot \left ( {\XB} - {\Xmean} \right ) \cdot \left ( {\XB} - {\Xmean} \right )^{{\bf T}} + \Q \in \Re^{\Nstate \times \Nstate}\,.
\end{eqnarray}
\end{subequations}
Here $\ONES_{\Nens \times 1} \in \Re^{\Nens \times 1}$ is a vector whose entries are all equal one.  $\Q$ is the covariance matrix of model errors. In the typical case where $\XB$ is an ensemble of model forecasts, the explicit addition of the matrix $\Q$ to the covariance formula is not necessary.
Instead, the effect of model errors can be accounted for by adding random vectors $\xi_{{\rm i}} \sim {\Normal} \left ( 0, \Q \right )$ to model states:
$\xb{i} \leftarrow \xb{i}  + \xi_{{\rm i}}$. Prior to any measurement, the forecast step provides the best estimation to the true vector state $\xt$ \cite{Suarez12}.

The vector of observations ${\bf y} \in \Re^{\Nobs \times 1}$ is available at $t_{\rm current}$, where  $\Nobs$ is the number of data points. 
The observations are related to the model state by the relation
\[
{\bf y} = \Lo\, \x + v
\]
where $\Lo \in \Re^{\Nobs \times \Nstate}$ is the observation operator which maps the model space state into the observed space, and
$v \sim {\Normal} \left ( 0 , \R \right )$ is a vector of observation errors, accounting for both instrument and representativeness errors.

In order to account for observation errors one forms the matrix $\Y \in \Re^{\Nobs \times \Nens}$ whose columns $\y{i} \in \Re^{\Nobs \times 1}$ are perturbed measurements \cite{Kova11}:
\begin{eqnarray}
 \nonumber
\displaystyle \Y &=& \left ( {\bf y} + \e{1}, {\bf y} + \e{2}, \ldots , {\bf y} + \e{\Nens} \right )   \\ 
\nonumber
&=&\left ( \y{1} , \y{2}, \ldots , \y{\Nens} \right )  \in \Re^{\Nobs \times \Nens},
\end{eqnarray}
The vectors $\e{i} \in \Re^{\Nobs \times 1}$ represent errors in the data, and are drawn from a normal distribution $\e{i} \sim {\Normal} \left ( 0 , \R \right )$. 
We denote
\begin{eqnarray} \nonumber
\displaystyle {\E} = \left ( \e{1}, \e{2}, \ldots, \e{\Nens}\right ) \in \Re^{\Nobs \times \Nens},
\end{eqnarray}
and the ensemble representation of the measurements error covariance matrix is
\begin{eqnarray} \nonumber
\displaystyle \R = \frac{1}{\Nobs-1} \cdot \left ( \E \cdot \E^{{\bf T}} \right ) \in \Re^{\Nobs \times \Nobs}\,.
\end{eqnarray}

The EnKF {\em analysis step} produces an ensemble of improved estimates (analyses) $\XA \in \Re^{\Nstate \times \Nens}$ by applying the Kalman filter
to each of the background ensemble members:
\begin{subequations}
\begin{eqnarray} 
\label{EnKF:ensemble-analysis}
\displaystyle \XA &=& \XB + \KAL \cdot \left ( {{\Y}-{\Lo} \cdot  {\XB} } \right ) \in \Re^{\Nstate \times \Nens}, \\
\label{EnKF:kalman-gain}
\KAL &=& {\PB}  \cdot {\Lo}^{{\bf T}} \cdot \left ( {\Lo}  \cdot \PB \cdot {\Lo}^{{\bf T}} + \R \right )^{-1} \in \Re^{\Nstate \times \Nobs},
\end{eqnarray}
\end{subequations}
where the matrix $\KAL \in \Re^{\Nstate \times \Nobs}$ is the Kalman gain and quantifies the contribution of the background--observations difference to the analysis. 

The EnKF {\em forecast step} uses the dynamical model operator $\Mo$ to evolve each member of the ensemble 
$\XA$ from the current time $t_{\rm current}$ to the next time $t_{\rm next}$ where observations are available:
\begin{eqnarray} 
\label{forecast-step}
\displaystyle {\XB}(t_{\rm next}) = {\Mo}_{t_{\rm current} \rightarrow t_{\rm next}} \left ( {\XA}(t_{\rm current})\right ) \in \Re^{\Nstate \times \Nens}\,.
\end{eqnarray}
The forecast ensemble ${\XB}$ is the background for the new EnKF cycle at $t_{\rm next}$. The analysis and forecast steps are repeated.

\subsection{Efficient implementations of the analysis step}
\label{enkf-analysis}

From equations \eqref{EnKF:ensemble-analysis}--\eqref{EnKF:kalman-gain} the analysis step can be written as 
\begin{eqnarray}
\label{EnKF:analysis-rewritten}
\XA= \XB + {\PB}  \cdot {\Lo}^{{\bf T}} \cdot \Z\,,
\end{eqnarray}
where $\Z \in \Re^{\Nobs \times \Nens}$ is the solution of the following linear system:
\begin{eqnarray}
\label{EnKF:system-to-solve}
\displaystyle \underbrace{\left ( {\Lo}  \cdot \PB \cdot {\Lo}^{{\bf T}} + \R \right )}_{\bf W \in \Re^{\Nobs \times \Nobs}} \cdot \Z = \left ( {{\Y}-{\Lo}  \cdot {\XB}} \right ) \in \Re^{\Nobs \times \Nens}\,.
\end{eqnarray}
A direct solution of this linear system can be obtained using the Cholesky decomposition for matrix inversion \cite{Lin99,Sch90,Xia10}. 
While this is a numerically stable and accurate approach  \cite{Gill96,Mei83,Ste97}, its application to \eqref{EnKF:system-to-solve} leads to the following complexity \cite{Jan06} of the analysis step:
\begin{eqnarray}
\label{EnKF:overall-complexity}
\BO \left( \Nobs^3 + \Nobs^2 \cdot \Nens + \Nobs \cdot \Nens^2 + \Nstate \cdot \Nens^2 \right )\,.
\end{eqnarray}
This is an acceptable complexity for a large number of degrees of freedom ($\Nstate$), but not for a large number of observations ($\Nobs$). An alternative is to solve \eqref{EnKF:system-to-solve}, and the overall analysis step, using Singular Value Decomposition (SVD) based methods. Those methods exploit the special structure of the data error covariance matrix $\R$, which is often (block) diagonal
and can be easily factorized:
\begin{eqnarray} \nonumber
\displaystyle \R = \diag{\R}{\rm N_{block}}, \text{ with $\R_{\k} \in \Re^{\rm N_{\k} \times N_{\k}}$}\, , 1 \le \k \le {\rm N_{\rm block}}\,,
\end{eqnarray}
where $\rm N_{block}$ is the number of blocks in the matrix $\R$ and:
\begin{eqnarray*}
\displaystyle \Nobs = \sum_{\k = 1}^{\N_{\rm block}} {\rm N_{\k}}\,.
\end{eqnarray*}
The matrix $\R$ is a covariance matrix, and in practice it is always positive definite.

The observation operator $\Lo \in \Re^{\Nobs \times \Nstate}$ is sparse or can be applied efficiently to a state vector. Then, when $\R$ is diagonal, we can express the system matrix \eqref{EnKF:system-to-solve} as follows:
\begin{eqnarray}
\nonumber
\widehat{\S} &=& \left ( \XB-\Xmean \right ) \in \Re^{\Nstate \times \Nens}\,,\\
\label{EnKF:SVD-rewritten-W}
\displaystyle {\bf W} &=& \sqrt{\R} \cdot \left [ \frac{1}{\Nens-1} \cdot \sqrt{\R^{-1}} \cdot \Lo \cdot \widehat{\S} \cdot \left( \Lo \cdot \widehat{\S} \right )^{\rm T} \cdot \sqrt{\R^{-1}} + {\bf I} \right ] \cdot \sqrt{\R}\,.
\end{eqnarray}
Employ the singular value decomposition
\begin{eqnarray}
\label{EnKF:SVD-single-value-decomposition}
\displaystyle \sqrt{\R^{-1}} \cdot \Lo \cdot \widehat{\S}= {\bf U} \cdot {\bf \Sigma} \cdot {\bf V^{\rm T}} \in \Re^{\Nobs \times \Nens},\,
\end{eqnarray}
where ${\bf U \in \Re^{\Nobs \times \Nobs}}$ and ${\bf V \in \Re^{\Nens \times \Nens}}$ are orthogonal square matrices, and ${\bf \Sigma} = {\bf diag}(\sigma_1,\sigma_2,\ldots,\sigma_{\Nens}) \in \Re^{\Nobs \times \Nens}$ is the diagonal matrix holding the singular values of ${\bf W} \in \Re^{\Nobs \times \Nobs}$. The linear system \eqref{EnKF:system-to-solve} can be written as follows \cite{Jan06}:
\begin{eqnarray}
\label{EnKF:SVD-linear-system-rewritten}
\displaystyle \left[{\sqrt{\R}} \cdot {\bf U} \cdot \left( \frac{\Sigma^2}{\Nens-1} + {\bf I} \right ) \cdot {\bf U^{\rm T}} \cdot {\sqrt{\R}} \right ] \cdot \Z = \left ( {{\Y}-{\Lo}  \cdot {\XB}} \right ) \in \Re^{\Nobs \times \Nens},
\end{eqnarray}
which yields the solution:
\begin{eqnarray}
\label{SVD:solution-linear-system}
\displaystyle \Z = \sqrt{\bf R^{-1}} \cdot {\bf U} \cdot {\bf diag} \left\{ \left( \frac{\sigma^2_i}{\Nens-1}+1 \right)^{-1} \right\} \cdot {\bf U^{\rm T}} \cdot \sqrt{\bf R^{-1}} \cdot  \left ( {{\Y}-{\Lo}  \cdot {\XB}} \right )\,.
\end{eqnarray}
The overall complexity of the analysis step
\begin{eqnarray}
\label{EnKF:SVD-overall-complexity}
\displaystyle \BO \left( \Nens^2 \cdot \Nobs + \Nens^3 + \Nstate \cdot \Nens^2 \right )
\end{eqnarray}
is suitable for large $\Nstate$ and $\Nobs$, assuming $\Nens$ remains small.
Many algorithms in the literature employ SVD ~\eqref{EnKF:SVD-linear-system-rewritten} for  the solution of the linear system \eqref{EnKF:system-to-solve} \cite{Evensen2009}. The analysis step is written in terms of the solution ~\eqref{SVD:solution-linear-system} in order to minimize the number of matrix computations. Due to this, the solution of the linear system and the improvement of the forecast ensemble are performed as a single step. The ensemble adjustment Kalman filter (EAKF) and the ensemble transform Kalman filter (ETKF) are based on this idea \cite{Tippett2003}. Other efficient implementations of the ensemble Kalman filter make use of SVD decompositions in order to derive pseudo-inverses, furthermore; these algorithms compute the inverse in the $\Nens$-dimensional ensemble space rather than $\Nobs$-dimensional measurement space. Thus, in practice, when $\Nobs \gg \Nens$, those algorithms exhibit a good performance.
All these methods have the overall complexity (number of long operations) given in  \eqref{EnKF:SVD-overall-complexity}.

A different approach is to employ iterative methods for solving the linear system \eqref{EnKF:system-to-solve}, for instance the conjugate gradient method \cite{Reid72,Golub89,Cohen72,Eis81} for $\Nens$ right-hand sides. However, each iteration costs $\BO \left( \Nens^2 \cdot \Nobs \right )$, therefore iterative methods do not seem to be competitive for the solution of \eqref{EnKF:system-to-solve}.

The well-established EnKF implementations presented above employ a Cholesky or SVD decomposition, which require considerable computational effort.
The next section discusses an efficient implementation of the ensemble Kalman filter which does not require any decomposition prior to the solution of the linear system \eqref{EnKF:system-to-solve}. 

\section{Iterative Implementation of the EnKF Analysis Step}
\label{sec:enkfSherman}

We make the assumptions \cite{Tippett2003,Jan06} that, in practice:
\begin{itemize}
 \item The data error covariance matrix $\R$ has a simple structure (e.g., is block diagonal). 
 \item The observation operator $\Lo$ is sparse or can be applied efficiently.
 \item The variables $\Nobs$ and $\Nstate$ are very large.
\end{itemize}
 Moreover, we consider the following situations:
\begin{itemize}
 \item In many real applications of the EnKF $\Nobs \gg \Nens$, and the number of variables ranges between $\BO(10^7)$ and $\BO(10^9)$.
 \item When many computational resources are available or when the number of components in the model state is relatively small  ($\Nstate \sim \BO(10^5)$, the number of ensemble members can be increased considerably in order to provide more accurate statistics. In this case $\Nobs \sim \Nens$. 
 \end{itemize}

Taking into in account the previous assumptions, we now derive the implementation of the EnKF. We define the matrix of member deviations $\S \in \Re^{\Nstate \times \Nens}$ as follows:
\begin{equation} 
\label{EnSM:Matrix-member-deviations} 
\displaystyle \S = \frac{1}{\sqrt{\Nens-1}} \cdot \left( \xb{1}-\xmean,\xb{2}-\xmean, \ldots, \xb{\Nens}-\xmean \right) \in \Re^{\Nstate \times \Nens}\,,
\end{equation}
which allows to write the ensemble covariance matrix as
\begin{eqnarray}
\label{EnKSM:covariance-matrix}
\displaystyle {\PB}=\S \cdot {\S}^{\bf T} \in \Re^{\Nstate \times \Nstate}\,.
\end{eqnarray}
By replacing equation \eqref{EnKSM:covariance-matrix} in \eqref{EnKF:system-to-solve}, the linear system solved during the analysis step is written as follows:
\begin{subequations}
\label{EnKSM:system-rewritten}
\begin{eqnarray}
\label{EnKSM:system-rhs}
&& \D = \Y-\Lo \cdot \XB \in \Re^{\Nobs \times \Nens} \,, \\
\label{EnKSM:member-observations}
&& \V = {\Lo} \cdot \S  \in \Re^{\Nobs \times \Nens}\,,\\
\label{EnKSM:system-to-solve-V}
&& \left( \R + \V \cdot \V^{\bf T} \right) \cdot \Z = \D \in \Re^{\Nobs \times \Nens}\,.
\end{eqnarray}
\end{subequations}
Note that
\begin{eqnarray} 
\nonumber
{\bf W} = \R + \V \cdot \V^{\bf T} 
&=& \R + \sum_{{\bf i}=1}^{\Nens}\v{i} \cdot \v{i}^{\bf T}
\end{eqnarray}
can be computed recursively via the sequence of matrices $\W{\k} \in \Re^{\Nobs \times \Nobs}$:
\begin{eqnarray}\nonumber
\displaystyle \W{0} &=& \R, \\ \nonumber
			  \W{k} &=& \R + \sum_{\rm i=1}^{{\bf k}} \v{ i} \cdot \v{i}^{{\bf T}} = \W{k-1}+\v{k} \cdot \v{k}^{{\bf T}}\,, \quad 1 \le {\bf k} \le \Nens\,,
\end{eqnarray}
therefore
\begin{eqnarray}
\label{EnKSM:sequence-to-replace}
{\bf W} = \W{\Nens} &=& \R + \sum_{{\bf i=1}}^{\Nens} \v{i} \cdot \v{i}^{{\bf T}} = \W{\Nens-1}+\v{\Nens} \cdot \v{\Nens}^{{\bf T}} \,.
\end{eqnarray}
By replacing equation \eqref{EnKSM:sequence-to-replace} in \eqref{EnKSM:system-to-solve-V} we obtain:
\begin{eqnarray}
\label{EnKSM:system-to-solve-W}
\displaystyle \left( \W{Nens-1} + \v{\Nens} \cdot \v{\Nens}^{{\bf T}}  \right) \cdot \Z = \D\,.
\end{eqnarray}
Theorem \ref{Theo:sequence-non-singular} shows that the matrix \eqref{EnKSM:sequence-to-replace} is non-singular. The linear system \eqref{EnKSM:system-to-solve-W} can be solved by making use again of the Sherman-Morrison formula \cite{Sher89}:
\begin{equation} 
\label{SWM}
\displaystyle \invM{\A+\L \cdot \M \cdot \N} = \SMF{\A}{\L}{\M}{\N}
\end{equation}
with $\A=\inv{\Nens-1}$, $\L = 1$, $\M = \v{Nens}$ and $\N=\v{Nens}^{\bf T}$. The solution of \eqref{EnKSM:system-to-solve-W} is computed as follows:
\begin{eqnarray}
\label{EnKSM:solution-sherman-morrison}
\displaystyle \Z = \F{\Nens}-\g{\Nens} \cdot \invM{1+\v{\Nens}^{\bf T} \cdot \g{\Nens}} \cdot \v{\Nens}^{\bf T} \cdot \F{\Nens} \in \Re^{\Nobs \times \Nens},
\end{eqnarray}
where $\F{\Nens} \in \Re^{\Nobs \times \Nens}$ and $\g{\Nens} \in \Re^{\Nobs \times \Nens}$ are given by the solution of the following linear systems:
\begin{subequations}
\label{EnKSM:subsystems}
\begin{equation}
\label{EnKSM:subsystem-F}
\displaystyle \W{\Nens-1} \cdot \F{\Nens} = \D \in \Re^{\Nobs \times \Nens}\,,
\end{equation}
\begin{equation}
\label{EnKSM:subsystem-g}
\displaystyle \W{\Nens-1} \cdot \g{\Nens} = \v{\Nens} \in \Re^{\Nobs \times 1}\,.
\end{equation}
\end{subequations}
Note that \eqref{EnKSM:subsystem-F} can be written as follows:
\begin{eqnarray}
\label{EnKSM:subsystem-f}
\displaystyle \W{\Nens-1} \cdot \f{\Nens}{i} = \d{i} \in \Re^{\Nobs \times 1}\,, \quad 1 \le {\bf i} \le \Nens\,,
\end{eqnarray}
where $\f{\Nens}{i} \in \Re^{\Nobs \times 1}$ and $\d{i} \in \Re^{\Nobs \times 1}$ are the i-th columns of the matrices $\F{\Nens}$ and $\D$, respectively. Following \eqref{EnKSM:solution-sherman-morrison}, the i-th column of the matrix $\Z$ is given by:
\begin{eqnarray}
\label{EnKMS:z-columns}
\displaystyle \z{i} = \f{\Nens}{i}-\g{\Nens} \cdot \invM{1+\v{\Nens}^{\bf T} \cdot \g{\Nens}} \cdot \v{\Nens}^{\bf T} \cdot \f{\Nens}{i} \in \Re^{\Nobs \times 1}\,.
\end{eqnarray}

By equation \eqref{EnKMS:z-columns}, the computation of $\Z$ involves the solution of the linear systems \eqref{EnKSM:subsystem-g} and \eqref{EnKSM:subsystem-f}. We apply the Sherman-Morrison formula  \eqref{SWM} again. The solution of the linear system \eqref{EnKSM:subsystem-f} can be obtained as follows:
\begin{eqnarray} \nonumber
\displaystyle \f{\Nens}{i} &=& \f{\Nens-1}{i}-\g{\Nens-1} \cdot   \invM{1+\v{\Nens-1}^{\bf T} \cdot \g{\Nens-1}} \cdot \\
 & &\qquad \cdot \v{\Nens-1}^{\bf T} \cdot \f{\Nens-1}{i}\,,
\end{eqnarray}
where $\f{\Nens-1}{i} \in \Re^{\Nobs \times 1}$ and $\g{\Nens-1} \in \Re^{\Nobs \times 1}$ are the solutions of the following linear systems, respectively:
\begin{eqnarray}
\label{EnKSM:subsystems-for-f} \nonumber
\displaystyle \W{\Nens-2} \cdot \f{\Nens-1}{i} &=& \d{i} \in \Re^{\Nobs \times 1}, \\ \nonumber
\W{\Nens-2} \cdot \g{\Nens-1} &=& \v{\Nens-1} \in \Re^{\Nobs \times 1}.
\end{eqnarray}

The linear system \eqref{EnKSM:subsystem-g} can be solved similarly. Note that the solution of each linear system involves the computation of two new linear systems, derived from the matrix sequence  \eqref{EnKSM:sequence-to-replace}. each of the new linear systems can be solved by applying recursively the Sherman-Morrison formula. For simplicity we denote by ${\bf f}$ and ${\bf g}$ the solutions of the new linear systems in each recursively application of the Sherman-Morrison formula. We have that:
{\allowdisplaybreaks
\begin{eqnarray}
\label{EnKSM:recursive-solution} \nonumber
\displaystyle \W{\Nens} \cdot \Z &=& \D \,, \\ \nonumber
\inv{\Nens} \cdot \d{i} &=& \underbrace{\inv{\Nens-1} \cdot \d{i}}_{\bf f}-\underbrace{\inv{\Nens-1} \cdot \v{\Nens}}_{\bf g} \\ \nonumber
&\cdot & \invM{1+\v{\Nens}^{\bf T} \cdot \underbrace{\inv{\Nens-1} \cdot \v{\Nens}}_{\bf g}} \\ \nonumber
&\cdot &\v{\Nens}^{\bf T} \cdot \underbrace{\inv{\Nens-1} \cdot \d{i}}_{\bf f} \in \Re^{\Nobs \times 1}, 1 \le {\rm i} \le \Nens \\ \nonumber
\inv{\Nens-1} \cdot \d{i} &=& \underbrace{\inv{\Nens-2} \cdot \d{i}}_{\bf f}-\underbrace{\inv{\Nens-2}\ \cdot \v{\Nens-1}}_{\bf g} \\ \nonumber
&\cdot & \invM{1+\v{\Nens-1}^{\bf T} \cdot \underbrace{\inv{\Nens-2} \cdot \v{\Nens-1}}_{\bf g}} \\ \nonumber
&\cdot & \v{\Nens}^{\bf T} \cdot \underbrace{\inv{\Nens-2} \cdot \d{i}}_{\bf f} \in \Re^{\Nobs \times 1}, 1 \le {\rm i} \le \Nens \\ \nonumber
\inv{\Nens-1} \cdot \v{\Nens} &=& \underbrace{\inv{\Nens-2} \cdot \v{\Nens}}_{\bf f}-\underbrace{\inv{\Nens-2}\ \cdot \v{\Nens-1}}_{\bf g} \\ \nonumber
&\cdot & \invM{1+\v{\Nens-1}^{\bf T} \cdot \underbrace{\inv{\Nens-2} \cdot \v{\Nens-1}}_{\bf g}} \\ \nonumber
&\cdot & \v{\Nens-1}^{\bf T} \cdot \underbrace{\inv{\Nens-2} \cdot \v{\Nens-1}}_{\bf f} \in \Re^{\Nobs \times 1}, \\ \nonumber
\vdots & =& \vdots \\ \nonumber
\inv{k} \cdot \x &=& \underbrace{\inv{k-1} \cdot \x}_{\bf f}-\underbrace{\inv{k-1}\ \cdot \v{k}}_{\bf g} \\ \nonumber
&\cdot & \invM{1+\v{k}^{\bf T} \cdot \underbrace{\inv{\k-1} \cdot \v{k}}_{\bf g}} \\ \nonumber
&\cdot & \v{k}^{\bf T} \cdot \underbrace{\inv{k-1} \cdot \x}_{\bf f} \in \Re^{\Nobs \times 1}, \\ \nonumber
\inv{k-1} \cdot \x &=& \underbrace{\inv{k-2} \cdot \x}_{\bf f}-\underbrace{\inv{k-2}\ \cdot \v{k-1}}_{\bf g} \\ \nonumber
&\cdot & \invM{1+\v{k-1}^{\bf T} \cdot \underbrace{\inv{k-2} \cdot \v{k-1}}_{\bf g}} \\ \nonumber
&\cdot & \v{k-1}^{\bf T} \cdot \underbrace{\inv{k-1} \cdot \x}_{\bf f} \in \Re^{\Nobs \times 1}, \\ \nonumber
\inv{k-1} \cdot \v{k} &=& \underbrace{\inv{k-2} \cdot \v{k}}_{\bf f}-\underbrace{\inv{k-2}\ \cdot \v{k-1}}_{\bf g} \\ \nonumber
&\cdot & \invM{1+\v{k-1}^{\bf T} \cdot \underbrace{\inv{k-2} \cdot \v{k-1}}_{\bf g}} \\ \nonumber
&\cdot & \v{k-1}^{\bf T} \cdot \underbrace{\inv{k-2} \cdot \v{k}}_{\bf f} \in \Re^{\Nobs \times 1}, \\ \nonumber
\vdots & =& \vdots \\ \nonumber
\inv{1} \cdot \x &=& \underbrace{\inv{0} \cdot \x}_{\bf f}-\underbrace{\inv{0} \cdot \v{1}}_{\bf g} \cdot  \invM{1+\v{1}^{\bf T} \cdot \underbrace{\inv{0} \cdot \v{1}}_{\bf g}} \\ \nonumber
&\cdot & \v{1}^{\bf T} \cdot \underbrace{\inv{0} \cdot \x}_{\bf f}, \\ \nonumber
\inv{0} \cdot \x &=& \invS{\R} \cdot \x \in \Re^{\Nobs \times 1}, \\ \nonumber
\inv{0} \cdot \v{1} &=& \invS{\R} \cdot \v{1} \in \Re^{\Nobs \times 1},
\end{eqnarray}
}
where $\x \in \Re^{\Nobs \times 1}$ can be either, a column of matrix $\D \in \Re^{\Nobs \times \Nens}$ or $\V \in \Re^{\Nobs \times \Nens}$. 

We note that:
\begin{itemize}
\item The computation of $\inv{k} \cdot \x$ involves the solution of the linear systems $\W{k-1} \cdot {\bf f} = \x$ and $\W{k-1} \cdot {\bf g} = \v{k}$.
\item Since the recursion is based on the sequence of matrices defined in \eqref{EnKSM:sequence-to-replace}, the base case is the linear system $\invS{\R} \cdot \x$ in which the matrix $\R$ is (block) diagonal.
\end{itemize}

From the previous analysis we derive a recursive Sherman-Morrison formula as follows. Define
\begin{eqnarray}
\label{EnKSM:recursive-SMF}
\displaystyle \fun{\x,\k} = 
\begin{cases}
{\bf z} = \invS{\R} \cdot \x\,, & \textnormal{for }\k = 0\,, \\
{\bf f} = \fun{\x,\k-1}\,; & \\
{\bf g} = \fun{\v{k},\k-1}\,, & \textnormal{for } 1\le \k \le \Nens\,, \\
{\bf z} = {\bf f}-{\bf g}\cdot \invM{1+\v{k}^{\bf T} \cdot {\bf g}} \cdot \v{k}^{\bf T} \cdot {\bf f}\,;
\end{cases}
\end{eqnarray}
where $\x \in \Re^{\Nobs \times 1}$. the columns of $\Z \in \Re^{\Nobs \times \Nens}$ are computed as follows:
\begin{eqnarray} \nonumber
\displaystyle \z{i} = \fun{\d{i},\Nens} \in \Re^{\Nobs \times 1}\,,\quad 1 \le {\rm i} \le \Nens.
\end{eqnarray}

\begin{figure}[H]
  \centering
	      \includegraphics[width=1\textwidth]{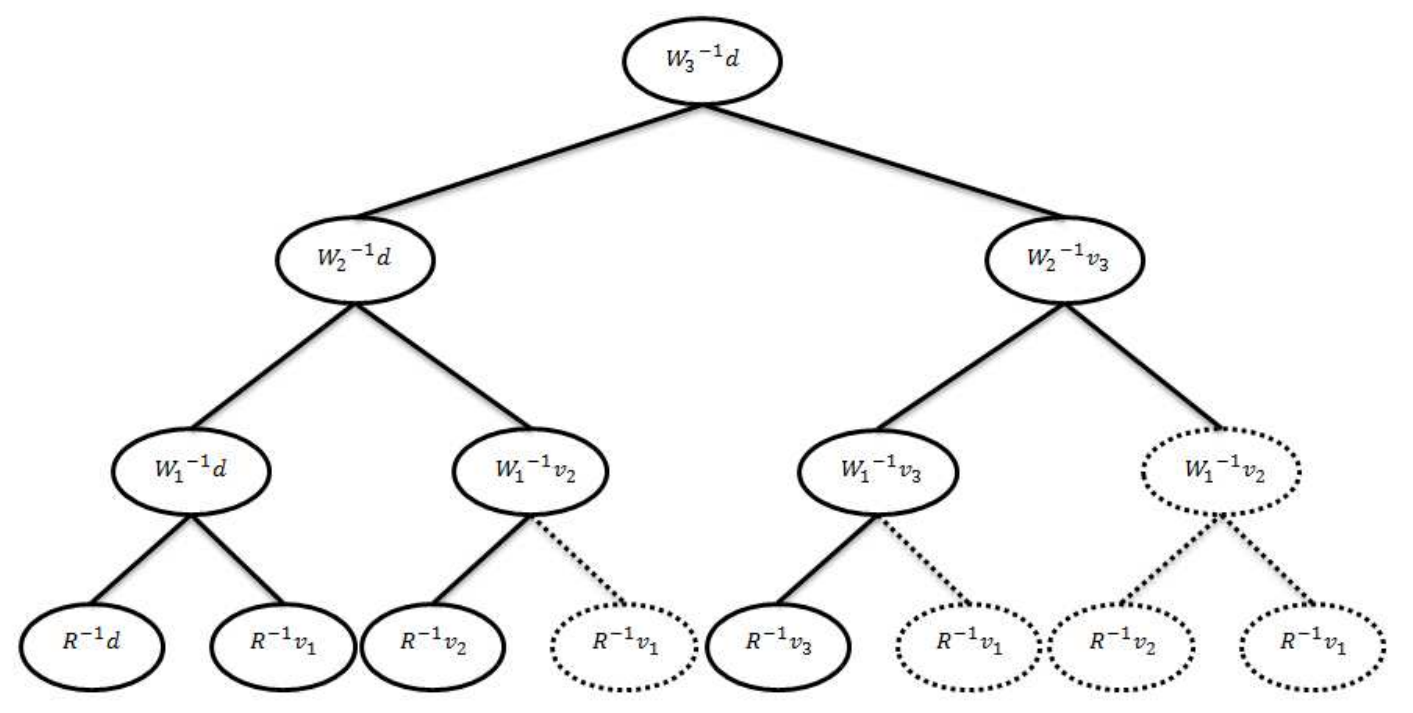}
  \caption{\label{Fig:sherman-full}The recursive Sherman-Morrison formula ($\fun{\bullet}$) applied to solve the linear system $\W{3} \cdot \Z=\D \in \Re^{\Nobs \times 3}$.   Here ${\bf d}$ is any column of matrix $\D \in \Re^{\Nobs \times 3}$. Dashed nodes represent repeated computations.}
\end{figure}

The recursive the computations performed by $\fun{\bullet}$ can be represented as a tree in which the solution ${\bf z} \in \Re^{\Nobs \times 1}$
 of each node depends on the computations of its left ($\bf f \in \Re^{\Nobs \times 1}$) and right ($\bf g \in \Re^{\Nobs \times 1}$) children (i.e., on the solutions of two linear systems). Figure \ref{Fig:sherman-full} illustrates the derivation of linear systems in order to solve $\W{\Nens} \cdot {\bf z} = {\bf d}$ for $\Nens=3$ and ${\bf d}\in \{ \d{1},\d{2},\d{3} \}$, ${\bf d} \in \Re^{\Nobs \times 1}$.
%

We see that $\fun{ \bullet}$ solves multiple times identical linear systems. For instance, the repeated computations performed in order to solve $\W{3} \cdot \Z = \D \in \Re^{\Nobs \times 3}$ are represented in Figure \ref{Fig:sherman-full} as dashed nodes. There, for instance, the linear system $\R \cdot {\bf g} = \v{1}$ is solved four times in the last level. The total number of linear systems to solve is $\BO(\Nens \cdot 2^{\Nens})$, i.e., it increases exponentially with regard to the number of ensemble members if identical computations are not avoided. Next subsection discusses how to achieve this and obtain an efficient implementation of the recursive Sherman-Morrison formula.
%
%

\subsection{An iterative Sherman-Morrison formula for matrix inversion}
\label{Sec:IterativeShermanMorrison}

In order to avoid identical computations  in Figure \ref{Fig:sherman-full} we can solve the linear systems from the last level of the tree up to the root level. We denote by ${\bf U} \in \Re^{\Nobs \times \Nens}$ and $\Z \in \Re^{\Nobs \times \Nens}$ the matrices holding partial results of the computations with regard to $\V$ and $\D$, respectively. 

\begin{figure}[H]
  \centering
	      \includegraphics[width=1\textwidth]{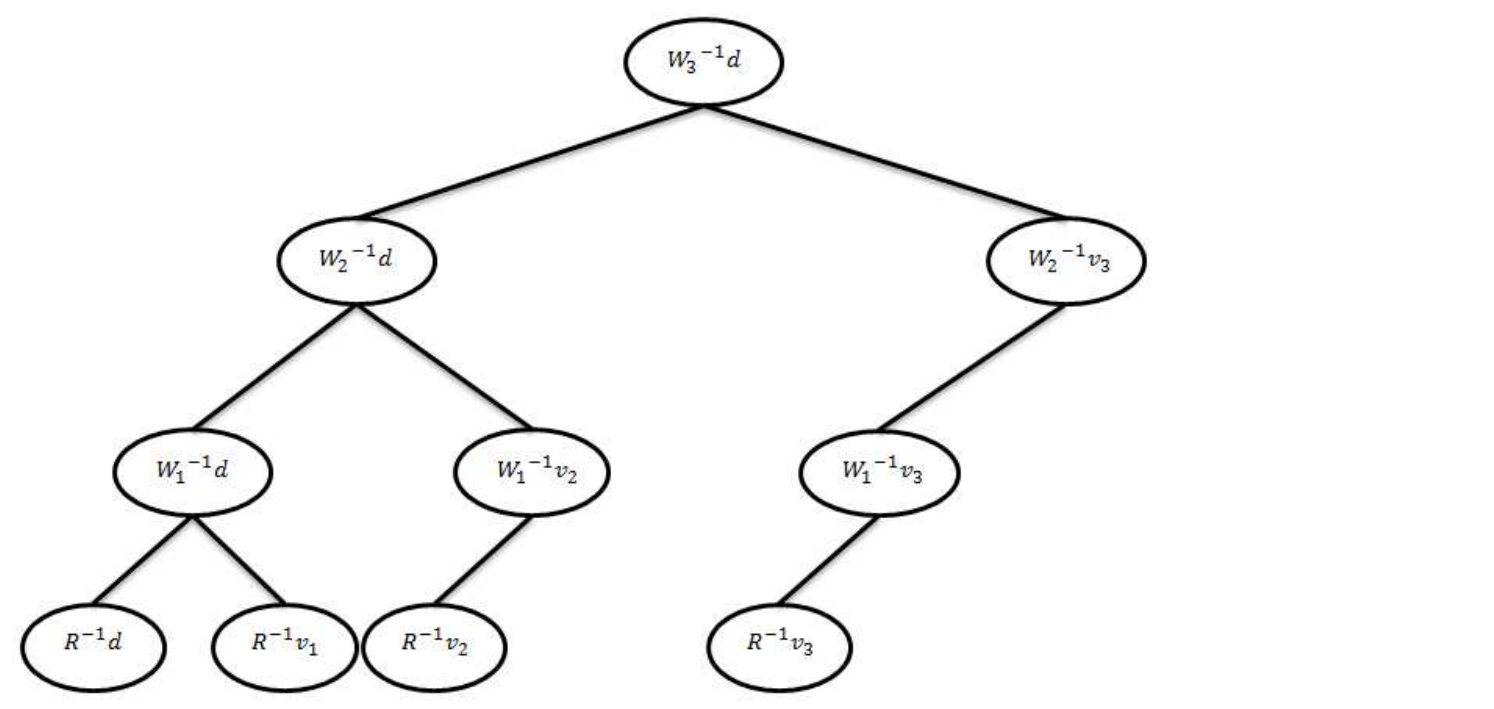}
  \caption{\label{Fig:sherman-optimal}Necessary computations for the solution of $\W{3} \cdot \Z=\D \in \Re^{\Nobs \times 3}$ using the recursive 
  Sherman-Morrison formula. This iterative version avoids all redundant computations.}
\end{figure}

Level $0$ can be computed as follows without any repeated effort: 
\begin{eqnarray} \nonumber
\displaystyle \U{0} &=& \left( \invS{\R} \cdot \v{1},\invS{\R} \cdot \v{2},\invS{\R} \cdot \v{3}\right ) \\ \nonumber
&=& \left( \inv{0} \cdot \v{1},\inv{0} \cdot \v{2},\inv{0} \cdot \v{3}\right ) \\ \nonumber
&=& \left( \u{1}{0},\u{2}{0},\u{3}{0}\right ), \\ \nonumber
\sZ{0} &=& \left( \invS{\R} \cdot \d{1},\invS{\R} \cdot \d{2},\invS{\R} \cdot \d{3}\right ) \\ \nonumber
&=& \left( \inv{0} \cdot \d{1},\inv{0} \cdot \d{2},\inv{0} \cdot \d{3}\right ) \\ \nonumber
&=& \left( \sz{1}{0},\sz{2}{0},\sz{3}{0}\right )\,. 
\end{eqnarray}
We make use of the Sherman-Morrison formula \eqref{SWM} and compute level 1 as follows:
\begin{eqnarray} \nonumber
\displaystyle \h{1} &=& \u{1}{0} \cdot \frac{1}{1+\v{1}^{\bf T} \cdot \inv{0} \cdot \v{1} } \\ \nonumber
&=& \u{1}{0} \cdot \frac{1}{1+\v{1}^{\bf T} \cdot \u{1}{0} }, \\ \nonumber
 \U{1} &=& \left( \u{1}{0}, \u{2}{0}-\h{1} \cdot \left( \v{1}^{\bf T} \cdot \u{2}{0} \right ), \u{3}{0}-\h{1} \cdot \left( \v{1}^{\bf T} \cdot \u{3}{0} \right )\right) \\ \nonumber
 &=& \left( \inv{0} \cdot \v{1},\inv{1} \cdot \v{2},\inv{1} \cdot \v{3}\right ) \\ \nonumber 
 &=& \left( \u{1}{1}, \u{2}{1}, \u{3}{1}\right), \\ \nonumber
 \sZ{1} &=& \left( \sz{1}{0}-\h{1} \cdot \left( \v{1}^{\bf T} \cdot \sz{1}{0}\right),\sz{2}{0}-\h{1} \cdot \left( \v{1}^{\bf T} \cdot \sz{2}{0}\right),\sz{3}{0}-\h{1} \cdot \left( \v{1}^{\bf T} \cdot \sz{3}{0}\right) \right ) \\ \nonumber
 &=& \left( \inv{1} \cdot \d{1},\inv{1} \cdot \d{2},\inv{1} \cdot \d{3}\right ) \\ \nonumber 
 &=& \left( \sz{1}{1},\sz{2}{1},\sz{3}{1} \right)\,.
\end{eqnarray}
Note that $\u{1}{0}$ has not been updated since it is not needed in the computations of the next levels. Similarly, the computations at level 2 make use of the Sherman-Morrison formula \eqref{SWM}:
\begin{eqnarray} \nonumber
\displaystyle \h{2} &=& \u{2}{1} \cdot \frac{1}{1+\v{2}^{\bf T} \cdot \inv{1} \cdot \v{2} } \\ \nonumber
&=& \u{2}{1} \cdot \frac{1}{1+\v{2}^{\bf T} \cdot \u{2}{1} }, \\ \nonumber
 \U{2} &=& \left( \u{1}{1}, \u{2}{1}, \u{3}{1}-\h{1} \cdot \left( \v{1}^{\bf T} \cdot \u{3}{1} \right )\right) \\ \nonumber
 &=& \left( \inv{0} \cdot \v{1},\inv{1} \cdot \v{2},\inv{2} \cdot \v{3}\right ) \\ \nonumber 
 &=& \left( \u{1}{2}, \u{2}{2}, \u{3}{2}\right), \\ \nonumber
 \sZ{2} &=& \left( \sz{1}{1}-\h{2} \cdot \left( \v{2}^{\bf T} \cdot \sz{1}{1}\right),\sz{2}{1}-\h{2} \cdot \left( \v{2}^{\bf T} \cdot \sz{2}{1}\right),\sz{3}{1}-\h{2} \cdot \left( \v{2}^{\bf T} \cdot \sz{3}{1}\right) \right ) \\ \nonumber
 &=& \left( \inv{2} \cdot \d{1},\inv{2} \cdot \d{2},\inv{2} \cdot \d{3}\right ) \\ \nonumber 
 &=& \left( \sz{1}{2},\sz{2}{2},\sz{3}{2} \right)\,.
\end{eqnarray}
The vectors $\u{1}{1}$ and $\u{2}{1}$ are not required for the computations of level 3, and they are not updated. Making use of the Sherman-Morrison formula
\eqref{SWM} once again, the root level is computed as follows:
\begin{eqnarray} \nonumber
\displaystyle \h{3} &=& \u{3}{1} \cdot \frac{1}{1+\v{3}^{\bf T} \cdot \invS{\W{1}} \cdot \v{3} } \\ \nonumber
 \sZ{3} &=& \left( \sz{1}{2}-\h{2} \cdot \left( \v{3}^{\bf T} \cdot \sz{1}{2}\right),\sz{2}{2}-\h{2} \cdot \left( \v{3}^{\bf T} \cdot \sz{2}{2}\right),\sz{3}{2}-\h{2} \cdot \left( \v{3}^{\bf T} \cdot \sz{3}{2}\right) \right ) \\ \nonumber
 &=& \left( \inv{3} \cdot \d{1},\inv{3} \cdot \d{2},\inv{3} \cdot \d{3}\right ) \\ \nonumber 
 &=& \inv{3} \cdot \D\,.
\end{eqnarray}
The computations performed by this iteration are shown in Figure \ref{Fig:sherman-optimal}. 

Some key features of the iteration are highlighted next.
\begin{itemize}
\item The number of iterations is $\Nens$.
\item At level 0 matrices $\sZ{0}$ and $\U{0}$ are computed as follows:
\begin{eqnarray*} \nonumber
\displaystyle \sZ{0} &=& \invS{\R} \cdot \D \in \Re^{\Nobs \times \Nens}, \\ \nonumber
\U{0} &=& \invS{\R} \cdot \V \in \Re^{\Nobs \times \Nens}.
\end{eqnarray*}
\item The matrix $\W{k} \in \Re^{\Nobs \times \Nobs}$ is never stored in memory. It can be represented implicitly by matrix $\V \in \Re^{\Nobs \times \Nens}$. 
This implicit representation realizes considerable memory savings, especially when  $\Nobs \gg \Nens$.
\item At iteration $\k$, only the columns $\u{i}{k}$ with $\k<{\rm i} \le \Nens$ are updated.
\end{itemize}

In summary, the solution of the linear system \eqref{EnKSM:system-to-solve-W} is obtained by the following iteration:
\begin{eqnarray}
\label{EnKSM:iterative-sherman-derivation} \nonumber
\displaystyle 
&& \text{ Level 0}: ~~ \left\{
\begin{array}{l}
\sZ{0} = \invS{\R} \cdot \D \\ \nonumber
\U{0} = \invS{\R} \cdot \V 
\end{array}
\right.\,, \\ \nonumber
& & \level{1}{0}{2}, \\ \nonumber
& & \quad \vdots \\ \nonumber
& & \level{k}{k-1}{k+1}, \\ \nonumber
& & \quad \vdots \\ \nonumber
& & \text{ Level $\Nens$}:~~\left\{ 
\begin{array}{l}
\h{\Nens} = \left(1+\v{\Nens}^{ \bf T } \cdot \u{\Nens}{\Nens-1} \right)^{-1}\,\u{\Nens}{\Nens-1}  \\
\Z = \sZ{\Nens} = \sZ{\Nens-1}-\h{\Nens} \cdot \left ( \v{\Nens}^{\bf T} \cdot \sZ{\Nens-1} \right ) \\
\end{array}
\right. \, .
\end{eqnarray}
Since the matrix $\R$ has a simple structure its inverse is easy to obtain. In the case of $\R$ (block) diagonal:
\begin{eqnarray} \nonumber
\R^{-1} = \diag{\invS{\R}}{{\rm N_{block}}} \in \Re^{\Nobs \times \Nobs},
\end{eqnarray}
and in general, under the assumptions done ($\R$ is easy to decompose), the computations $\sZ{0} = \invS{\R} \cdot \D \in \Re^{\Nobs \times \Nens}$ and $\U{0} = \invS{\R} \cdot \V \in \Re^{\Nobs \times \Nens}$ can be performed with no more than $\BO(\Nens^2 \cdot \Nobs)$ long operations.

Putting it all together, we define the iterative Sherman-Morrison formula $\ifun{\R}{\V}{\D}$ as follows:
\begin{itemize}
\item \textbf{Step 1.} Compute the matrices $\sZ{0} \in \Re^{\Nobs \times \Nens}$ and $\U{0} \in \Re^{\Nobs \times \Nens}$ as follows:
\begin{subequations}
\label{SM-step-1}
\begin{eqnarray} 
\displaystyle 
\sZ{0} &=& \invS{\R} \cdot \D,\\ 
\U{0} &=& \invS{\R} \cdot \V,
\end{eqnarray}
\end{subequations}
where $\invS{\R}$ is computed according to its special structure. 

\item \textbf{Step 2.} For $\k = 1$ to $\Nens$ compute:
\begin{subequations}
\label{SM-step-2}
\begin{eqnarray} 
\label{SM-step-2-h}
\h{\k} &=& \left(1+\v{\k}^{ \bf T} \cdot \u{\k}{\k-1}\right)^{-1} \, \u{\k}{\k-1}  \in \Re^{\Nobs \times 1}, \\  
\label{SM-step-2-Z}
\sZ{\k} &=& \sZ{\k-1} - \h{\k} \cdot \left( \v{\k}^{\bf T} \cdot \sZ{\k-1} \right) \in \Re^{\Nobs \times \Nens}, \\ 
\label{SM-step-2-u}
\u{\rm i}{\k} &=& \u{\rm i}{\k-1} - \h{\k} \cdot \left( \v{k}^{\bf T} \cdot \u{\rm i}{\k-1} \right) \in \Re^{\Nobs \times 1}, \text{ $\k+1 \le {\rm i} \le \Nens$}.
\end{eqnarray} 
\end{subequations}
\end{itemize}

We now use the iterative Sherman-Morrison formula in the analysis step to obtain
an efficient implementation of the Ensemble Kalman filter (SMEnKF). This filter is as follows.
%
%
The background ensemble states $\XB$ are obtained from the forecast step \eqref{forecast-step},
the ensemble mean $\xmean$ is given by \eqref{ensemble-background-mean}, and the ensemble deviations
form the mean $\S$ are given by \eqref{EnSM:Matrix-member-deviations}. 
The analysis is obtained as follows:
\begin{eqnarray} \nonumber
\displaystyle \D &=& \Y-\Ho \left (  \XB \right ) \in \Re^{\Nobs \times \Nens}, \\ \nonumber
\V &=& \Ho \left ( \S \right ) \in \Re^{\Nobs \times \Nens}, \\ \nonumber
\Z &=& \ifun{\R}{\V}{\D} \in \Re^{\Nobs \times \Nens}, \\ \nonumber
\XA &=& \XB+\S \cdot \V^{\bf T} \cdot \Z \in \Re^{\Nstate \times \Nens},
\end{eqnarray}
where the function $\Ho \left ( {\bf G} \right ) \in \Re^{\Nobs \times \Nens}$ is an efficient implementation of the observation operator
applied to several state vectors, represented by ${\bf G} \in \Re^{\Nstate \times \Nens}$. 

\subsubsection{Inflation aspects}
Inflation increases periodically the ensemble spread, such as to compensate for the small ensemble size, to
simulate the existence of model errors, and to avoid filter divergence \cite{LiHong2009}.	
All the inflation techniques applied in traditional EnKF can be used, virtually without modification, 	
in the context of SMEnKF. For example, after the forecast step, one can increase the spread of the ensemble
\[
\x_{\rm i} \leftarrow \xb{i} + \alpha \, (\x_{\rm i}-\xb{i})\,, \quad 1 \le \rm i \le \Nens\,,
\]
such as the ensemble covariance $\PB$ is increased by a factor $\alpha^2$  \cite{WuZheng2011}. 

\subsubsection{Localization aspects}
Using \eqref{EnKF:kalman-gain}, the analysis step can be written as follows:
\begin{eqnarray}  \nonumber
\displaystyle \XA &=& \XB+ \DXB \\ \nonumber
\DXB &=& \S \cdot \V^{\bf T} \cdot \Z \\   \nonumber
&=& {\PB}  \cdot {\Lo}^{{\bf T}} \cdot \left ( {\Lo}  \cdot \PB \cdot {\Lo}^{{\bf T}} + \R \right )^{-1} \D \in \Re^{\Nstate \times \Nens}\\\
\label{EnKSM:Localization}
&\approx& \PB_L  \cdot {\Lo}^{{\bf T}} \cdot \left ( {\Lo}  \cdot \PB_L \cdot {\Lo}^{{\bf T}} + \R \right )^{-1} \D \in \Re^{\Nstate \times \Nens}\,.
\end{eqnarray}
Localization techniques are explained in detail in \cite{Anderson200799}.
Localization replaces the ensemble based $\PB$ by a matrix $\PB_L = \rho \circ \PB$ in \eqref{EnKSM:Localization}, where $\rho$ is a localization matrix and 
 $\circ$ represents the Schur product. 
 
Clearly  localization in the form \eqref{EnKSM:Localization} requires the full covariance matrix, and cannot be applied in the context of the iterative Sherman-Morrison implementation. 
Applying SMEnKF with a single data point $\y{i}$ leads to a correction $\DXB_{\{\rm i\}}$, which can be localized
by multiplication with a {\em diagonal} matrix $\hat{\Delta}_{\{\rm i\}}$ that scales down components with the negative exponential of their distance to the observation $i$ location, and sets them to zero if outside the radius of influence:
\[
\DXB_{\{\rm i\}} = \hat{\Delta}_{\{\rm i\}} \cdot \S \cdot \V^{\bf T} \cdot \Z_{\{\rm i\}}
\]
This can be applied {\em in succession} for all data points to obtain a fully localized solution.

We discuss next a general approach to perform {\em partial localization}.
Let $\x_{\rm i}$ be an individual component of the state vector and $\y{j}$ an individual observation. Define the impact factor $\delta_{i,j} \in [0,1]$
of the information in $\y{j}$ on the state point $\x_{\rm i}$. For example, one can use a
correlation length, and a radius about the measurement location outside which  the impact factor is zero.
Define the influence matrix $\Delta = (\delta_{i,j}) \in \Re^{\Nstate \times \Nobs}$, and
replace \eqref{EnKSM:Localization} with the following partial localization formula
\begin{eqnarray*}  
\DXB  &\approx& \PB_L  \cdot {\Lo}^{{\bf T}} \cdot \left ( {\Lo}  \cdot \PB \cdot {\Lo}^{{\bf T}} + \R \right )^{-1} \d{\ell} \in \Re^{\Nstate \times \Nens} \\
&=& \Delta \circ \left( \S \cdot \V^{\bf T} \right) \cdot \Z \,.
\end{eqnarray*}
The (i, $\ell$)-th entry contains the i-th component the correction vector for the $\ell$-th ensemble member and reads
\begin{eqnarray}  
\label{localized-correction-component}
\DXB_{\rm i, \ell} &=& \sum_{k=1}^{\Nens}  \S_{\rm i,\rm k} \, \sum_{j=1}^{\Nobs} \delta_{\rm i, \rm j}\, \V_{\rm k,\rm j}  \, \Z_{\rm j, \ell} 
\,,~~ 1 \le \rm i \le \Nstate\,,~ 1 \le \ell \le \Nens\,.
\end{eqnarray}
The components of the correction matrix \eqref{localized-correction-component} are independent of one another, and can be evaluated in parallel
after the system solution $\Z$ has been computed.

\subsection{Computational complexity}
\label{Sec:computational-complexity}

In the complexity analysis of the iterative Sherman-Morrison formula we count only the long operations (multiplications and divisions). Moreover, as discussed before, we make the assumptions presented in \cite{Tippett2003,Jan06}, namely, the data error covariance matrix $\R \in ^{\Nobs \times \Nobs}$ is inexpensive to decompose, and the observation operator $\Lo$ can be applied efficiently to any vector. We now analyze each each step of the iterative Sherman-Morrison formula when $\R$ is diagonal, the extension to nondiagonal data error covariance matrices is inmediate.

In the first step \eqref{SM-step-1} each row ${\rm i}$ of matrices $\D \in \Re^{\Nobs \times \Nens}$ and $\V \in \Re^{\Nobs \times \Nens}$ is divided by the corresponding component $\rv_{\rm i} \in \Re^{\Nobs}$ in order to obtain $\sZ{0} \in \Re^{\Nobs \times \Nens}$ and $\U{0} \in \Re^{\Nobs \times \Nens}$ respectively. This yields to $\Nobs \cdot \Nens$ number of long operation for each matrix, therefore:
\begin{eqnarray}
\label{EnKSM:operations-count-step-one}
\displaystyle \T{\Nens}{\Nobs}{step1} = 2 \cdot \Nobs \cdot \Nens.
\end{eqnarray}

In the second step  \eqref{SM-step-2} we compute the vector $\h{\k} \in \Re^{\Nobs}$ \eqref{SM-step-2-h}, and the matrices $\sZ{\k} \in \Re^{\Nobs \times \Nens}$ \eqref{SM-step-2-Z} and $\U{\k} \in \Re^{\Nobs \times \Nens}$ \eqref{SM-step-2-u}. The number of long operations for each of one are as follows:
\begin{eqnarray} \nonumber
\displaystyle
\h{\k} &=& \overbrace{\u{\k}{\k-1} \cdot \frac{1}{1+\underbrace{\v{\k}^{ \bf T} \cdot \u{\k}{\k-1}}_{\Nobs}}}^{\Nobs}\,,  \\ \nonumber
\sZ{\k} &=& \sZ{\k-1} - \overbrace{\h{\k} \cdot \left( \underbrace{\v{\k}^{\bf T} \cdot \sZ{\k-1}}_{\Nens \cdot \Nobs} \right)}^{\Nens \cdot \Nobs}\,, \\ \nonumber
\u{\rm i}{\k} &=& \u{\rm i}{\k-1} - \overbrace{\h{\k} \cdot \left( \underbrace{\v{k}^{\bf T} \cdot \u{\rm i}{\k-1}}_{\Nobs} \right)}^{ \Nobs}\,.
\end{eqnarray}
Since the second step \eqref{SM-step-2} is performed $\Nens$ times, the number of long operations can be expressed as:
\begin{eqnarray} 
 \nonumber
\displaystyle \T{\Nens}{\Nobs}{step2} &=& \sum_{\k = 1}^{\Nens} \left( \underbrace{2 \cdot \Nobs}_{\h{k}}+\underbrace{2 \cdot \Nobs \cdot \Nens}_{\sZ{k}}+\overbrace{\sum_{\rm j = 1}^{\k-1}\left (2 \cdot \Nobs \right)}^{\u{\rm j}{\k}}\right ) \\ \nonumber
&=& 2 \cdot \Nens \cdot \Nobs + 2 \cdot \Nens^2 \cdot \Nobs + \sum_{\rm k = 1}^{\Nens}(\k-1) \cdot 2\cdot \Nobs \\ \label{EnKSM:step2-long-operations}
&=& 3 \cdot \Nens^2 \cdot \Nobs + \Nens \cdot \Nobs\,.
\end{eqnarray}
Consequently, from \eqref{EnKSM:operations-count-step-one}--\eqref{EnKSM:step2-long-operations}, we have
\begin{eqnarray} \nonumber
\T{\Nens}{\Nobs}{SMF} &=& \underbrace{2 \cdot \Nobs \cdot \Nens}_{\T{\Nens}{\Nobs}{step1}} + \underbrace{3 \cdot \Nens^2 \cdot \Nobs + \Nens \cdot \Nobs}_{\T{\Nens}{\Nobs}{step2}} \\ \nonumber
&=&3 \cdot ( \Nens^2 \cdot \Nobs + \Nens \cdot \Nobs)\,,
\end{eqnarray}
which yields a complexity of
\begin{eqnarray}
\label{EnKSM:upper-bound} 
\BO \left( \Nens^2 \cdot \Nobs \right )\,.
\end{eqnarray}
Note that when $\R$ is not diagonal, under the assumptions done, the computations  \eqref{SM-step-1} of $\sZ{0}$ and $\U{0}$ can be efficiently performed in $\BO(\Nobs \cdot \Nens^2)$ long operations;  the overall effort becomes $3 \cdot ( \Nens^2 \cdot \Nobs + \Nens^2 \cdot \Nobs)$. This leads to the same complexity  \eqref{EnKSM:upper-bound} for $\R$ diagonal, block diagonal, or in general easy to decompose. 

The overall complexity of the analysis step for the iterative formula
\begin{eqnarray} \nonumber
\label{EnKSM:long-computations-complexity-analysis}
\displaystyle \XA = \XB + \underbrace{\S \cdot \overbrace{\V^{\rm T} \cdot \underbrace{\Z}_{\BO \left( \Nens^2 \cdot \Nobs \right )}}^{ \BO \left (\Nens^2 \cdot \Nobs \right )}}_{\BO \left( \Nens^2 \cdot \Nstate \right)}
\end{eqnarray}
is:
\begin{eqnarray}
\label{EnKSM:overall-complexity-analysis}
\displaystyle \BO \left ( \Nens^2 \cdot \Nobs + \Nens^2 \cdot \Nstate \right ),
\end{eqnarray}
The complexity of the proposed implementation of the EnKF is equivalent to the upper bounds os the methods described in \cite{Tippett2003}, as detailed in the Table  \ref{Tab:Analysis-Complexity-Contrast}. The term $\Nens^3$ does not appear in the upper-bound of the proposed method even when $\R$ is not diagonal. This term can affect the performance of the EnKF when $\Nens \sim \Nobs$.

\begin{table}[H]
\centering
\begin{tabular}{|l|l|} \hline
\bf Analysis method & \bf Computational cost \\ \hline
Direct \cite{Tippett2003} & $\BO \left( \Nens^2 \cdot \Nobs + \Nens^3 + \Nens^2 \cdot \Nstate\right )$ \\
Serial \cite{Anderson07} (for each observation) & $\BO \left( \Nens \cdot \Nobs + \Nens \cdot \Nobs \cdot \Nstate \right )$ \\
ETKF \cite{Anderson01} & $\BO \left( \Nens^2 \cdot \Nobs + \Nens^3 + \Nens^2 \cdot \Nstate\right )$ \\
EAKF \cite{Anderson01} & $\BO \left( \Nens^2 \cdot \Nobs + \Nens^3 + \Nens^2 \cdot \Nstate\right )$ \\
{\bf Proposed EnKF Implementation} & $\BO \left( \Nens^2 \cdot \Nobs + \Nens^2 \cdot \Nstate\right )$ \\ \hline
\end{tabular}
\caption{Summary of computational costs of the analysis steps for several ensemble filters. The costs are functions of the ensemble size $\Nens$, number of observations $\Nobs$ and state dimension $\Nstate$.}
\label{Tab:Analysis-Complexity-Contrast}
\end{table}

Maponi \cite{Maponi2007276} proposed a general approach based on the Sherman Morrison formula  to solve linear systems. The application of this generic algorithm 
to \eqref{EnKSM:system-to-solve-W} leads to an increased computational cost as the special structure of the system (and special structure of $\R$) are not exploited. The generic algorithm 
applied to EnKF analysis 
\begin{eqnarray}
\label{W-single-representation}
\displaystyle \W{\Nens} \z{i} = \d{i} \in \Re^{\Nobs}\, , \text{ for $1 \le \rm i \le \Nens$}\, ,
\end{eqnarray}
uses the decomposition \cite[Remark 1]{Maponi2007276}:
\begin{eqnarray*}
 \displaystyle \W{\Nens} = \W{0}+\sum_{\rm i=1}^{\Nobs} \um{\rm i} \cdot \v{\rm i}^{\rm T} \, ,
\end{eqnarray*}
where $\W{0}=\diag{w}{\Nobs} \in \Re^{\Nobs \times \Nobs}$ is a diagonal matrix holding the diagonal entries of $\W{\Nens}$, $\um{i}$ is the i-th column of ${\bf U} = \W{\Nens}-\W{0} \in \Re^{\Nobs \times \Nobs}$ and $\v{i}={\bf e}_{\rm i}$ is the i-th element of the canonical basis in $\Re^{\Nobs}$. Thus, according to \cite[Corollary 4]{Maponi2007276}, each linear system \eqref{W-single-representation} can be solved with $\BO(\Nobs^3)$ long operations, leading to a total of
\begin{eqnarray*}
 \displaystyle \BO\left( \Nens \cdot \Nobs^3 \right )\, .
\end{eqnarray*}
Therefore the computational cost of the analysis step is:
\begin{eqnarray}
\label{Maponi:computational-cost}
 \displaystyle \BO \left( \Nens \cdot \Nobs^3 + \Nens \cdot \Nstate \right )\, ,
\end{eqnarray}
which is larger than the computational cost of our proposed EnKF implementation when either $\Nobs \gg \Nens$ or $\Nens \sim \Nobs$. Moreover, according to \citep[Theorem 3]{Maponi2007276}, when $\Nobs \gg \Nens$, the solution of linear system ~\eqref{W-single-representation} can be computed with no more than $\BO( \Nens^2 \cdot \Nobs + \Nens^2)$ long operations. The resulting computational cost of the analysis step is:
\begin{eqnarray*}
 \displaystyle \BO \left( \Nens^3 \cdot \Nobs + \Nens^3 + \Nens \cdot \Nstate \right )\, ,
\end{eqnarray*}
which is similar to the computational costs of the ETKF and EAKF methods when $\Nobs \gg \Nens$. In this case, it is unclear how to construct the matrix $\bf U \in \Re^{\Nobs \times \Nens}$ according to Maponi's method; $\bf U$ can not be chosen as we propose since $\W{0}$ must be diagonal and $\V$ differs from our definition in ~\eqref{EnKSM:member-observations}. In addition, Maponi's algorithm requires the explicit representation in memory of the matrix $\W{\Nens}$, which, in practice, is $\BO(10^7 \times 10^{7})$ dimensional. In contradistinction, $\W{\Nens}$ is not required explicitly in memory by our iterative Sherman Morrison formula. 

Lastly, the stability conditions of Maponi's method are not discussed in \cite{Maponi2007276}. Furthermore, the sequence of matrices $\W{\k}$ are not proved to be non-singular, which is crucial for the well-performance of that method. On the contrary, the stability analysis of the iterative Sherman Morrison formula is discussed in the next section.

\subsection{Stability Analysis}
\label{EnKSM:stability-iterative-method}

The solution of the linear system \eqref{EnKF:system-to-solve} by the iterative Sherman Morrison formula yields the next sequence of matrices during the computation of $\left[\W{\Nens}\right]^{-1}$:
\begin{eqnarray*}
\displaystyle 
\inv{0} &=& \invS{\R} \\
\inv{1} &=& \invS{\R}-\frac{1}{\gamma_{\k}} \cdot \u{1}{0} \cdot \v{1}^{\rm T} \cdot \invS{\R} = \left( {\bf I}-\frac{1}{\gamma_{\k}} \cdot \u{1}{0} \cdot \v{1}^{\rm T}\right ) \inv{0} \\
& \vdots & \\
\inv{\k} &=& \left( {\bf I}-\frac{1}{\gamma_{\k}} \cdot \u{\k}{\k-1} \cdot \v{\k}^{\rm T}\right ) \inv{\k-1} 
\end{eqnarray*}
where 
\begin{eqnarray*}
\displaystyle \gamma_{\k} = 1+\v{\k}^{\rm T} \cdot \u{\k}{\k-1} \in \Re\,,\text{ for $1 \le \k \le \Nens$.}
\end{eqnarray*}

The following situations may affect the proposed method:
\begin{itemize}
\item If any step produces $\gamma_{\k} = 0$, then subsequent steps cannot proceed.
\item Round-off errors can be considerably amplified if $\gamma_{\k} \approx 0$ (numerical instability). 
\item If any matrix $\W{\k}$ in the sequence:
\begin{eqnarray}
\label{EnKSM:sequence-of-matrices}
\displaystyle \left \{  \W{0}, \W{1}, \ldots, \W{\Nens} \right \}\,,
\end{eqnarray}
is singular, the algorithm cannot proceed. 
\end{itemize} 

We now show that the positive definiteness of the covariance matrix $\R$ is a sufficient in order to guarantee the stability of the iterative Sherman Morrison formula.

\begin{theorem}
\label{Theo:sequence-positive}
Assume that $\R$ is positive definite with $\xi^T\, \R \xi \ge \alpha\, \Vert \xi \Vert^2$ for any $\xi \in \mathbb{R}^{\Nobs}$.
Then all matrices $\W{\k}$ are positive definite with $\xi^T\, \W{\k} \xi \ge \alpha\, \Vert \xi \Vert^2$ for any $\xi \in \mathbb{R}^{\Nobs}$.
\end{theorem} 
 
\begin{proof}
First, $\W{0} = \R$ is positive definite. Next, we proceed by finite induction and assume that $\W{\k-1}$ is positive definite
with $\xi^T\, \W{\k-1} \xi \ge \alpha\, \Vert \xi \Vert^2$.
From \eqref{EnKSM:system-to-solve-W} we have that:
\begin{eqnarray*}
\displaystyle \W{\k} &=& \W{\k-1}+\v{\k} \cdot \v{\k}^{\rm T},
\end{eqnarray*}
and therefore $\W{\k}$ is also positive definite:
\begin{eqnarray*}
\xi^T\, \W{\k} \, \xi &=& \underbrace{ \xi^T\,\W{\k-1} \, \xi }_{ \ge \alpha\, \Vert \xi \Vert^2}+ \underbrace{ \left( \xi^T\, \v{\k} \right)^2}_{\ge 0}  \ge \alpha\, \Vert \xi \Vert^2 \quad \forall \, \xi\in \mathbb{R}^\Nobs\,.
\end{eqnarray*}
\end{proof} 
 
\begin{theorem}
\label{Theo:non-zero-gamma-values}
Assume that $\R$ is positive definite. The sequence of values $\gamma_{k}$ generated by the algorithm are strictly greater than one for all $1\le \k \le \Nens$.
\end{theorem}

\begin{proof}
By the iterative Sherman Morrison formula, the common computations ($\u{\k+1}{\k}$) are given by:
\begin{eqnarray}
\displaystyle  \nonumber
\u{\rm 1}{0} &=& \invS{\R} \cdot \v{1} = \left[ \W{0}\right ]^{-1} \cdot \v{1} \\ \nonumber
\u{\rm 2}{1} &=&  \underbrace{\left( {\bf I}-\frac{1}{\gamma_{1}} \cdot  \u{1}{0} \cdot \v{1}^{\rm T} \right ) \cdot \left[ \W{0} \right]^{-1}}_{\left [\W{1} \right ]^{-1}} \cdot \v{2} = \left[ \W{1}\right ]^{-1} \cdot {\v{2}} \\ \nonumber
& \vdots & \\ \label{EnKSM:inverse-iteration-k}
\u{\k+1}{\k} &=&  \underbrace{\left( {\bf I}-\frac{1}{\gamma_{\k}} \cdot \u{\k}{\k-1} \cdot \v{\k}^{\rm T} \right ) \cdot \left[ \W{\k-1}\right ]^{-1}}_{\left [\W{\k} \right ]^{-1}} \cdot \v{\k+1} = \left[ \W{\k}\right ]^{-1} \cdot {\v{\k+1}} 
\end{eqnarray}
Since $\W{\k-1}$ is positive definite we have:
\begin{eqnarray*}
\displaystyle \gamma_{\k} = 1+\v{\k}^{\rm T} \cdot \u{\k}{\k-1} =  1+ \underbrace{\v{\k}^{\rm T} \cdot \inv{k-1} \cdot \v{\k}}_{>0} > 1\, ,
\end{eqnarray*}
consequently $\gamma_{\k}>1$ for all $1 \le \k \le \Nens-1$.
\end{proof}

We have the following direct corollary of Theorem \ref{Theo:sequence-positive}.

\begin{theorem}
\label{Theo:sequence-non-singular}
Assume that $\R$ is positive definite. At iteration $\k$, the linear system:
\begin{eqnarray}
\label{EnKSM:sub-linear-system-to-solve}
\displaystyle \W{\k} \cdot \sZ{\k} = \D\, ,
\end{eqnarray}
has a unique solution, for $1 \le \k \le \Nens$.
\end{theorem}

\subsection{Pivoting}
\label{pivoting}

Theorem \ref{Theo:non-zero-gamma-values} shows that $\gamma_{\k}$ values cannot be near zero. Due to this, we expect that the round-off errors will not increase considerably during an iteration of the iterative Sherman Morrison formula since:
\begin{eqnarray*}
\displaystyle \frac{1}{\gamma_{\k}} \in (0\,,\,1)\,.
\end{eqnarray*}

The following pivoting strategy can be (optionally) applied in order to further decrease round-off error accumulation. It consists of interchanging the columns of matrices $\V$ and $\U{k-1}$ such that the pair $ ( \v{\rm j},\u{\rm j}{\k-1} )$ maximizes $\gamma_{\k}$. Formally, at iteration $\k$, prior the matrix computations \eqref{SM-step-2}, we look for a column index $\ind$ such that:
\begin{eqnarray}
\label{EnKSM:gamma-value}
\displaystyle \ind = \arg \max_{\substack{\rm i}} \left \{ \left | 1+\v{\rm i}^{\rm T} \cdot \u{\rm i}{\k-1} \right |\, , \k \le {\rm i} \le \Nens \right \}\,,
\end{eqnarray}
and then, the columns $\k$ and $\ind$ are interchanged in matrices $\V$ and $\U{\k-1}$. 

The iterative Sherman Morrison formula with pivoting gives the next computational cost:
\begin{eqnarray*} \nonumber
{\rm T}_{\rm SMF}^{\rm PIV} \left ( {\Nens},{\Nobs} \right ) &=& \underbrace{3 \cdot ( \Nens^2 \cdot \Nobs + \Nens \cdot \Nobs)}_{\T{\Nens}{\Nobs}{SMF}} + \underbrace{\sum_{\k=1}^{\Nens} \underbrace{\sum_{\rm i=1}^{\k} \Nobs}_{\eqref{EnKSM:gamma-value}}}_{\text{Pivoting}}\,, \\
&=& \frac{7}{2} \cdot \left( \Nens^2 \cdot \Nobs + \Nens \cdot \Nobs \right )
\end{eqnarray*}
which yields to:
\begin{eqnarray*}
{\rm T}_{\rm SMF}^{\rm PIV} \left ( {\Nens},{\Nobs} \right ) \in \BO \left( \Nens^2 \cdot \Nobs \right )\, ,
\end{eqnarray*}
from which we can conclude that seeking the maximum value of $\gamma_{\k}$ according to \eqref{EnKSM:gamma-value} does not increase the computational cost of the iterative Sherman Morrison formula. Consequently, the overall complexity in the analysis step remains bounded by ~\eqref{EnKSM:overall-complexity-analysis}. 

\subsection{Parallel implementation}
\label{parallelimpl}

In this section we discuss an efficient parallel implementation of the iterative Sherman-Morrison formula. Since the algorithm \eqref{SM-step-1}--\eqref{SM-step-2} can be applied individually to each column of the matrices $\sZ{0} \in \Re^{\Nobs \times \Nens}$ and $\U{0} \in \Re^{\Nobs \times \Nens}$, there are $2\,\Nens$ computations that can be performed in parallel. We define the matrix ${\G{0}} \in \Re^{\Nobs \times 2\cdot \Nens}$ holding the columns of $\V \in \Re^{\Nobs \times \Nens}$ and $\D \in \Re^{\Nobs \times \Nens}$ as follows:
\begin{eqnarray}
\label{ParEnKSM:G0}
\displaystyle \G{0} &=& \left[ \V, \D\right] = \left [ \v{1},\v{2}, \ldots,  \v{\Nens} , \d{1}, \d{2},  \ldots, \d{\Nens} \right ], \\ \nonumber
  &=& \left [\gj{1}{0}, \ldots, \gj{\Nens}{0},  \gj{\Nens+1}{0}, \ldots,\gj{2 \cdot \Nens}{0} \right ] \in \Re^{\Nobs \times 2 \cdot \Nens}, 
\end{eqnarray}
Let $\Nproc{0}$ be the number of available processors at the initial time. The number of operations per processor is 
\begin{eqnarray} \nonumber
\label{EnKSM:assigned-process}
\displaystyle \nc{0} = \frac{2 \cdot \Nens}{\Nproc{0}}\,.
\end{eqnarray}
The matrix \eqref{ParEnKSM:G0}  can be written as 
\begin{eqnarray} \nonumber
\label{EnKSM:block-definition-initial-G}
\displaystyle 
{\G{0}} = \left [ \Bi{1}{0}, \Bi{2}{0}, \ldots , \Bi{\Nproc{t}}{0}\right ] \,,
\end{eqnarray}
where the blocks $\Bi{i}{0} \in \Re^{\Nobs \times \nc{0}}$ are 
\begin{eqnarray} \nonumber
\displaystyle 
\Bi{i}{0} = \left[ \gj{(i-1)\nc{0}+1}{0}, \gj{(i-1)\nc{0}+2}{0}, \ldots \gj{i\nc{0}}{0} \right ] \in \Re^{\Nobs \times \nc{0}} \text{ for  $1 \le {\rm i} \le \Nproc{0}$},
\end{eqnarray}
The parallel, first step \eqref{SM-step-1} of the iterative Sherman-Morrison formula is implemented as an update over the blocks:
\begin{eqnarray} \nonumber
\displaystyle 
\Bi{i}{1} = \invS{\R} \cdot \Bi{i}{0} \in \Re^{\Nobs \times \nc{0}}\,,\quad \text{ for all $1 \le {\rm i} \le \Nproc{0}$} \,,
\end{eqnarray}
which yields
\begin{eqnarray} \nonumber
\displaystyle \G{1} &=& \left [\invS{\R} \cdot \Bi{1}{0},\invS{\R} \cdot \Bi{2}{0},\ldots,\invS{\R} \cdot \Bi{\Nproc{0}}{0} \right ]  \\ \nonumber
&=& \left [\Bi{1}{1},\Bi{2}{1},\ldots,\Bi{\Nproc{0}}{1} \right ]  \\ \nonumber
 &=& \left [ \underbrace{ \gj{1}{1}, \ldots, \gj{\Nens}{1} }_{\U{0}}, \underbrace{\gj{\Nens+1}{1}, \ldots, \gj{2 \cdot \Nens}{1}}_{\sZ{0}} \right ]  \\
 &=& \left[ \u{1}{0}, \u{2}{0}, \ldots, \u{\Nens}{0} , \sZ{0} \right]\, 
\end{eqnarray}

The second step \eqref{SM-step-2} of the iterative Sherman-Morrison formula consists of a sequence of updates applied to the matrices $\sZ{0}$ and $\U{0}$. Such matrices are represented by the columns of matrix $\G{1}$. Thus, consider the computation of level one, each column of the matrix $\G{1}$ can be updated as follows:
\begin{eqnarray} \nonumber
\displaystyle 
\gj{i}{2} = \gj{i}{1}-\gj{1}{1} \cdot \left (1+\v{1}^{\rm T} \cdot \gj{1}{1} \right )^{-1} \cdot \left( \v{1}^{\rm T} \cdot \gj{i}{1} \right ) \in \Re^{\Nobs \times 1}\,,
\quad 2 \le \rm i \le 2 \cdot \Nens\,.
\end{eqnarray}
Similarly to the first step, the computations can be grouped in blocks
\begin{eqnarray} \nonumber
\displaystyle 
\Bi{i}{1} = \left[ \gj{(i-1)\nc{1}+2}{1}, \gj{(i-1)\nc{1}+3}{1}, \ldots \gj{i\nc{1}+1}{1} \right ] \in \Re^{\Nobs \times \nc{1}}\,,\quad \text{ for  $1 \le {\rm i} \le \Nproc{1}$},
\end{eqnarray}
and distributed over the processors:
\begin{eqnarray} \nonumber
\displaystyle 
\Bi{i}{2} = \Bi{i}{1}-\gj{1}{1} \cdot \left (1+\v{1}^{\rm T} \cdot \gj{1}{1} \right )^{-1} \cdot \left( \v{1}^{\rm T} \cdot \Bi{i}{1} \right ) \in \Re^{\Nobs \times \nc{1}}\,,
\quad \mbox{for all }1 \le \rm i \le \Nproc{1}\,.
\end{eqnarray}
Note that $\gj{1}{1}$ ($\u{1}{0}$) is not updated since it is not required in subsequent computations. 
Thus, for the matrix $\G{2}$
\begin{eqnarray} \nonumber
\displaystyle 
\G{2} = \left[ \gj{1}{1}, \gj{2}{2}, \gj{3}{2}, \ldots, \gj{2 \cdot \Nens}{2}\right ]\, ,
\end{eqnarray}
the next common computation is $\gj{2}{2}$ ($\u{2}{1}$), and for the same reasons, this vector is not updated. 

In general, at time step $\rm t$, $1 \le \rm t \le \Nens$, the first $\rm t$ columns of the matrix $\G{t}$ are not included in the update process:
\begin{eqnarray} \nonumber
\displaystyle 
\G{t} = \left[ \gj{1}{1}, \gj{2}{2}, \ldots, \gj{t-1}{t-1}, \gj{t}{t}, \gj{t+1}{t}, \ldots, \gj{2 \cdot \Nens}{t}\right ]\, ,
\end{eqnarray}
The parallel computation of \eqref{SM-step-2} at time step $t$ is performed as follows:
\begin{itemize}
\item Compute the number of computation units (columns of matrix $\G{t} \in \Re^{\Nobs \times \Nens}$) per processor:
\begin{eqnarray} \nonumber
\displaystyle \nc{t} = \frac{2 \cdot \Nens - {\rm t}}{\Nproc{t}}\,.
\end{eqnarray}
\item Perform the update in parallel over the blocks:
\begin{eqnarray} \nonumber
\displaystyle 
\Bi{i}{t} = \Bi{i}{t}-\gj{t}{t} \cdot \left (1+\v{t}^{\rm T} \cdot \gj{t}{t} \right )^{-1} \cdot \left( \v{t}^{\rm T} \cdot \Bi{i}{t} \right ) \in \Re^{\Nobs \times \nc{t}}\,,\quad
\mbox{for all } 1 \le {\rm i} \le \Nproc{t}\,,
\end{eqnarray}
where 
\begin{eqnarray} \nonumber
\displaystyle 
\Bi{i}{t} = \left[ \gj{(i-1)\nc{1}+1+t}{1}, \gj{(i-1)\nc{1}+2+t}{1}, \ldots \gj{i\nc{1}+t}{1} \right ] \in \Re^{\Nobs \times \nc{t}}\,.
\end{eqnarray}
\end{itemize}

This parallel implementation of the iterative Sherman Morrison formula leads to the complexity:
\begin{eqnarray} \nonumber
\displaystyle \rm T_{\rm SMF}^{\rm PAR}(\Nobs,\Nens) &=& \underbrace{\BO \left( \nc{0} \cdot \Nobs \right )}_{\rm step 1} + \sum_{\rm t = 1}^{\Nens}\underbrace{\BO \left( \nc{t} \cdot \Nobs \right )}_{\rm step 2}\, .
\end{eqnarray}

Notice, when the number of processors at time $0 \le \rm t \le \Nens$ is $\Nproc{t} = 2 \, \Nens - \rm t$ then $\nc{t} = 1$. Hence, the corresponding computational cost of the analysis step is bounded by:
\begin{eqnarray} \nonumber
\displaystyle  \BO(\Nobs \cdot \Nens + \Nstate \cdot \Nens) \,,
\end{eqnarray}
therefore, when the number of observations is large enough relative to the number of ensemble members, this parallel approach of the iterative Sherman-Morrison formula exhibits a linear behavior, making this implementation attractive. 

\section{Experimental Results}
\label{sec:results}

In this section several computation tests are conducted in order to assess the accuracy and running time of the EnKF based on iterative Sherman Morrison formula. 
\subsection{Experimental setting}
\label{Results:experimental-settings}

The Sherman-Morrison EnKF implementation as well as the EnKF implementations based on Cholesky and SVD are coded in Fortran 90. The Cholesky and SVD decompositions use functions from the LAPACK library \cite{Anderson90} as follows: 
\begin{itemize}
\item The matrix $\bf W \in \Re^{\Nobs \times \Nobs}$ is built using DSYRK as follows:
\begin{eqnarray} \nonumber
\label{Results:W-computation}
\displaystyle \W{\Nens} = \alpha \cdot \V \cdot \V^{\rm T} + \beta \cdot \R \in \Re^{\Nobs \times \Nobs}\,,\quad
\mbox{with}~~ \alpha = \frac{1}{{\Nens-1}}\,,~~\beta = 1.0\,.
\end{eqnarray}
\item The functions DPOTRF and DPOTRI are used to compute the Cholesky decomposition of matrix $\bf W \in \Re^{\Nobs \times \Nobs}$.
\item The SVD decomposition is performed tby the DGESVD function.
\end{itemize}

In order to measure the quality of the solutions we employ the following performance metrics. The Elapsed Time (ET) measures the overall simulation time
for a method $*$. This metric is defined as follows:
\begin{eqnarray}
\label{Results:Elapsed-time-metric}
\displaystyle \rm ET(*) = Forecast_{*}+Analysis_{*}
\end{eqnarray}
Where $\rm Forecast_{*}$ and $\rm Analysis_{*}$ are the running time for the overall forecast and analysis steps respectively. 

The Root Mean Square Error (RMSE) is defined as follows:
\begin{eqnarray} \nonumber
\displaystyle \varepsilon \left ( * \right ) = {\rm RMSE} = \frac{1}{{\rm N_{steps}}} \cdot \left ( \sum_{{\rm t = 1}}^{{\rm N_{steps}}}{{\rm RSE}_{{\rm t}}} \right )
\end{eqnarray}
where ${\rm N_{steps}}$ is the number of time steps and ${\rm RSE_t}$ is the Root Square Error at time $\rm t$ defined as follows:
\begin{eqnarray} \nonumber
\displaystyle {\rm RSE_t} = \sqrt{\frac{1}{\Nstate} \cdot \left( \xt_{\rm t}-{\bf \overline{x}^C_{\rm t}}\right)^{\bf T} \cdot \left(\xt_{\rm t}-{\bf \overline{x}^C_{\rm t}} \right)}
\end{eqnarray}
where $\xt_{\rm t}$ is the true vector state at time $\rm t$, and ${\bf x^C_t}$ can be either the ensemble mean in the forecast $\xmean$ or analysis $\xmeana$ at time $\rm t$. As can be seen the  RMSE measures in average the distance between a reference solution ($\xt_{\rm t}$) and the given solution (${\bf x^C_{\rm t}}$).

The EnKF implementations are tested on two systems: the Lorenz 96 model \cite{Lorenz98} representing the atmosphere, and a quasi-geostrophic model \cite{Carton94} representing the ocean. They define the model operators ($\mathcal M$) in the EnKF experiments. To compare the performance of different EnKF implementations we measure the elapsed times and the accuracy of analyses for different values of $\Nobs$ and $\Nens$. 

\subsection{Lorenz-96 model ($\bf \N_{ens} \sim \N_{obs}$)}
\label{Results:Lorenz-Model}

The Lorenz 96 model is described by the following system of ordinary differential equations \cite{Lorenz98}:

\begin{eqnarray} \label{eqn:Lorenz96}
\displaystyle \Lorenz,
\end{eqnarray}
which has been heuristically formulated in order to take into in account properties of global atmospheric models such as the advection,  dissipation and forcing. This model exhibits extended chaos with an external forcing value ($\rm F=8$), when the solution is in the form of moving waves. For this reason, the model is adequate to perform basic studies of predictability.

The test assesses how the efficiency of the EnKF implementations depend on the input parameters $\Nobs$ and $\Nens$ when $\Nobs \sim \Nens$ (the number of observations and ensemble members are relatively close). The experimental setting is described below.

\begin{itemize}
\item One time unit of the Lorenz 96 model corresponds to five days of the atmosphere. The observations are made over 500 days (100 time units).
\item The background error is assumed to be $5\% $, i.e., the initial ensemble mean's deviation from the reference solution is drawn from a normal distribution
whose standard deviation is $5\%$ of the reference value.
\item The external forcing is set to $\rm F = 8.0$. 
\item The dimensions of the model state  are $\Nstate \in \{ 500, 1000, 3000, 5000\}$. While the typical dimension for the Lorenz-96 model is $\Nstate=40$,
we scale the system up to assess the performance of different implementations.
\item The number of observations equals the number of states, $\Nobs = \Nstate$. Due to this, the analysis step involves large linear systems of size $\bf W \in \Re^{\Nstate \times \Nstate}$. 
\item The number of ensemble members $\Nens$ depends on the size of the state vector as shown in Table  \ref{Tab:Experimental-Settings-Nens}.

\item At each time $\rm t$, the synthetic observations are constructed as follows:
\begin{eqnarray}
\label{eq:syntethic-observations}
\displaystyle \y{t} = \xt_{\rm t}+\e{t} \in \Re^{\Nstate}
\end{eqnarray}
since the number of observations and variables of the vector state are the same. $\e{t}$ belongs to a normal distribution with zero mean and covariance matrix
\begin{eqnarray*}
\displaystyle \R =  {\bf diag}\left \{ 0.01^2 \right \} \in \Re^{\Nstate \times \Nstate} \, ,
\end{eqnarray*} 
as is usual in practice. The errors are replicated for each compared, EnKF implementation. Due to this, the same data errors are hold for all tests. 

\item The assimilation window is five model days.

\item The localization \eqref{localized-correction-component} is applied using the influence factors
\begin{eqnarray} \nonumber
\label{Results:Lorenz-Multiplication}
\delta_{\rm i, \rm j} = \rm exp \left( -\frac{min\{i,j\}}{\Nstate}\right ) \in \Re^{\Nstate \times \Nens}
\end{eqnarray}
where $\rm min\{i,j\}$ is the minimum distance between the indexes $\rm i$ and $\rm j$ of the vector state; this distance accounts for the
periodic boundary conditions in \eqref{eqn:Lorenz96}.

\end{itemize}

The RMSE results are shown in Table  \ref{Tab:Lorenz-Results-RMS}.
All methods provide virtually identical analyses. As expected, the analysis improves when the size of the ensemble is increased. 

\begin{table}[H]
\centering
{
\begin{tabular}{|c|c|} \hline
${\bf {N_{state}}}$ &  ${\bf {N_{ens}}}$ \\ \hline
$500$ & \{200,250,300,350,400\} \\ \hline
$1000$ & \{400,450,500,550,600\} \\ \hline
$3000$ & \{900,950,1000,1050,1100\} \\ \hline
$5000$ & \{1500,1550,1600\} \\ \hline
\end{tabular}
}
\caption{Number of ensemble members $\Nens$ with respect to the dimension of the vector state $\Nstate$.}
\label{Tab:Experimental-Settings-Nens}
\end{table}

\begin{table}[H]
\centering
{\scriptsize
\begin{tabular}{|c|c|c|c|c|c|} \hline
${\bf N_{state}/N_{obs}}$ & {\bf Step} & $\bf N_{ens}$ & ${\bf EnKF_{Sher}}$ & ${\bf  EnKF_{Chol}}$ & ${\bf EnKF_{SVD}}$  \\ \hline
\multirow{10}{*}{500/500} & \multirow{5}{*}{Forecast}  & 200 & $0.006491846885685$ & $0.006491846885685$ & $0.006491846885685$  \\ 
&  & 250 & $0.003089844737585$ & $0.003089844737585$ & $0.003089844737585$ \\ 
&  & 300 & $ 0.001923620680204$ & $ 0.001923620680204$ & $ 0.001923620680204$  \\ 
&  & 350 & $0.001501324727984$ & $0.001501324727984$ & $0.001501324727984$  \\ 
&  & 400 & $0.001238879857327$ & $0.001238879857327$ & $0.001238879857327$ \\ \cline{2-6}
& \multirow{5}{*}{Analysis} &  200 & $0.005406046859821$ & $0.005406046859821$ & $0.005406046859821$  \\ 
& & 250 & $ 0.002577923867043$ & $ 0.002577923867043$ & $ 0.002577923867043$ \\ 
& & 300 & $0.001595214779293$ & $0.001595214779293$ & $0.001595214779293$  \\ 
& & 350 & $0.001239247824936$ & $0.001239247824936$ & $0.001239247824936$  \\ 
& & 400 & $0.001019599347515$ & $0.001019599347515$ & $0.001019599347515$ \\ \hline
\multirow{10}{*}{1000/1000} & \multirow{5}{*}{Forecast}  & 400 & $0.004416601672762$ & $0.004416601672762$ & $0.004416601672762$  \\ 
&  & 450 & $0.002808556499338$ & $0.002808556499338$ & $0.002808556499338$ \\ 
&  & 500 & $0.002256025116209$ & $0.002256025116209$ & $0.002256025116209$  \\ 
&  & 550 & $0.001827747641793$ & $0.001827747641793$ & $0.001827747641793$  \\ 
&  & 600 & $0.001561032221877$ & $0.001561032221877$ & $0.001561032221877$ \\ \cline{2-6}
& \multirow{5}{*}{Analysis} &  400 & $0.003643849246395$ & $0.003643849246395$ & $0.003643849246395$  \\ 
& & 450 & $0.002320171220995$ & $0.002320171220995$ & $0.002320171220995$ \\ 
& & 500 & $0.001863322297334$ & $0.001863322297334$ & $0.001863322297334$  \\ 
& & 550 & $0.001503808242092$ & $0.001503808242092$ & $0.001503808242092$  \\ 
& & 600 & $0.001281928557481$ & $0.001281928557481$ & $0.001281928557481$ \\ \hline
\multirow{10}{*}{3000/3000} & \multirow{5}{*}{Forecast}  & 900 & $0.009439626047886$ & $0.009439626047886$ & $0.009439626047886$  \\ 
&  & 950 & $0.007199551317193$ & $0.007199551317193$ & $0.007199551317193$ \\ 
&  & 1000 & $0.005410752525373$ & $0.005410752525373$ & $0.005410752525373$  \\ 
&  & 1050 & $0.004299142614958$ & $0.004299142614958$ & $0.004299142614958$  \\ 
&  & 1100 & $0.003476460994219$ & $0.003476460994219$ & $0.003476460994219$ \\ \cline{2-6}
& \multirow{5}{*}{Analysis} &  900 & $0.007957954481050$ & $0.007957954481050$ & $0.007957954481050$  \\ 
& & 950 & $0.006074278547494$ & $0.006074278547494$ & $0.006074278547494$ \\ 
& & 1000 & $0.004552394986586$ & $0.004552394986586$ & $0.004552394986586$  \\ 
& & 1050 & $0.003597428757064$ & $0.003597428757064$ & $0.003597428757064$  \\ 
& & 1100 & $0.002902565379165$ & $0.002902565379165$ & $0.002902565379165$ \\ \hline
\multirow{6}{*}{5000/5000} & \multirow{3}{*}{Forecast}  & 1500 & $0.007678111830765$ & $0.007678111830765$ & $0.007678111830765$  \\ 
&  & 1550 & $0.006378636076218$ & $0.006378636076218$ & $0.006378636076218$ \\ 
&  & 1600 & $0.005608079561230$ & $0.005608079561230$ & $0.005608079561230$ \\ \cline{2-6}
& \multirow{3}{*}{Analysis} &  1500 & $0.006490749509210$ & $0.006490749509210$ & $0.006490749509210$  \\ 
& & 1550 & $0.005389012789042$ & $0.005389012789042$ & $0.005389012789042$ \\ 
& & 1600 & $ 0.004728980886456$ & $ 0.004728980886456$ & $0.004728980886456$ \\ \hline
\end{tabular}
}
\caption{RMSE for the Lorenz-96 model with different number of states. When the number of ensemble members is increased the estimation of the true vector state is improved. All EnKF implementations provide virtually identical results.}
\label{Tab:Lorenz-Results-RMS}
\end{table}


\begin{figure}[H]
\centering
\begin{sideways} $\rm RMSE \left( \XA \right )$ \end{sideways} 
\subfigure[Sherman]{
\includegraphics[width=0.25\textwidth]{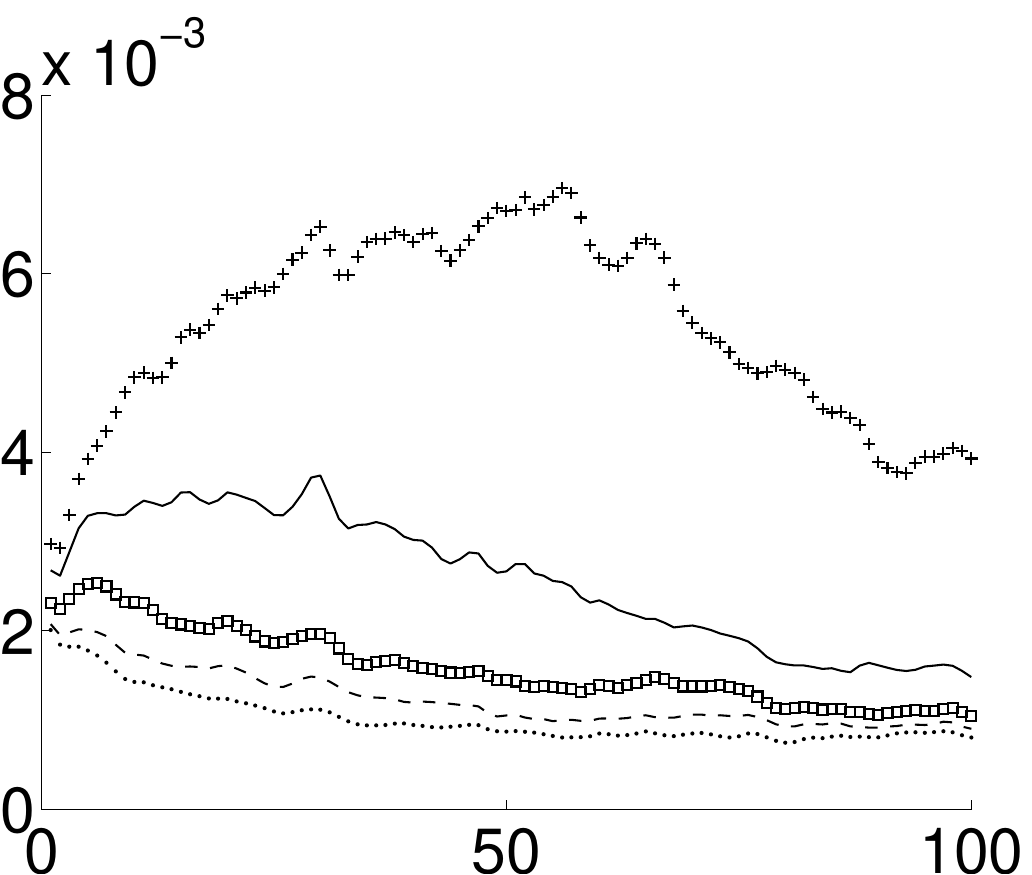} 
}
\subfigure[Cholesky]{
\includegraphics[width=0.25\textwidth]{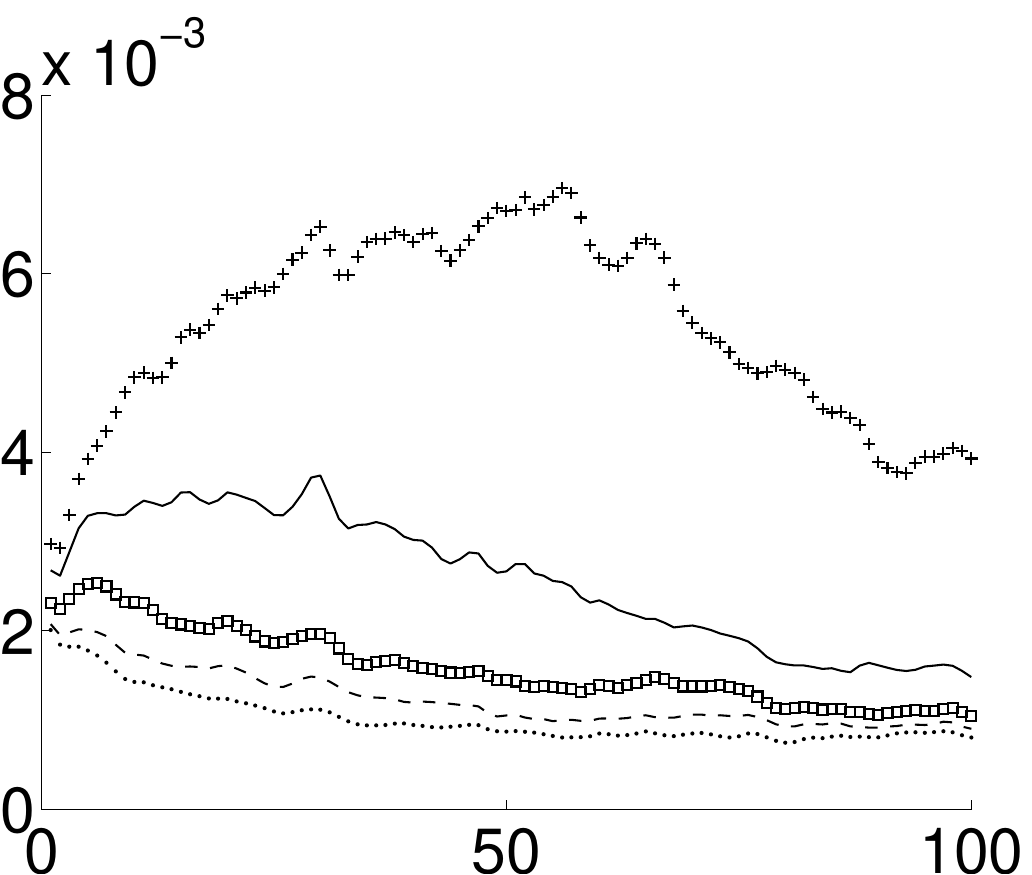}
}
\subfigure[SVD]{
\includegraphics[width=0.25\textwidth]{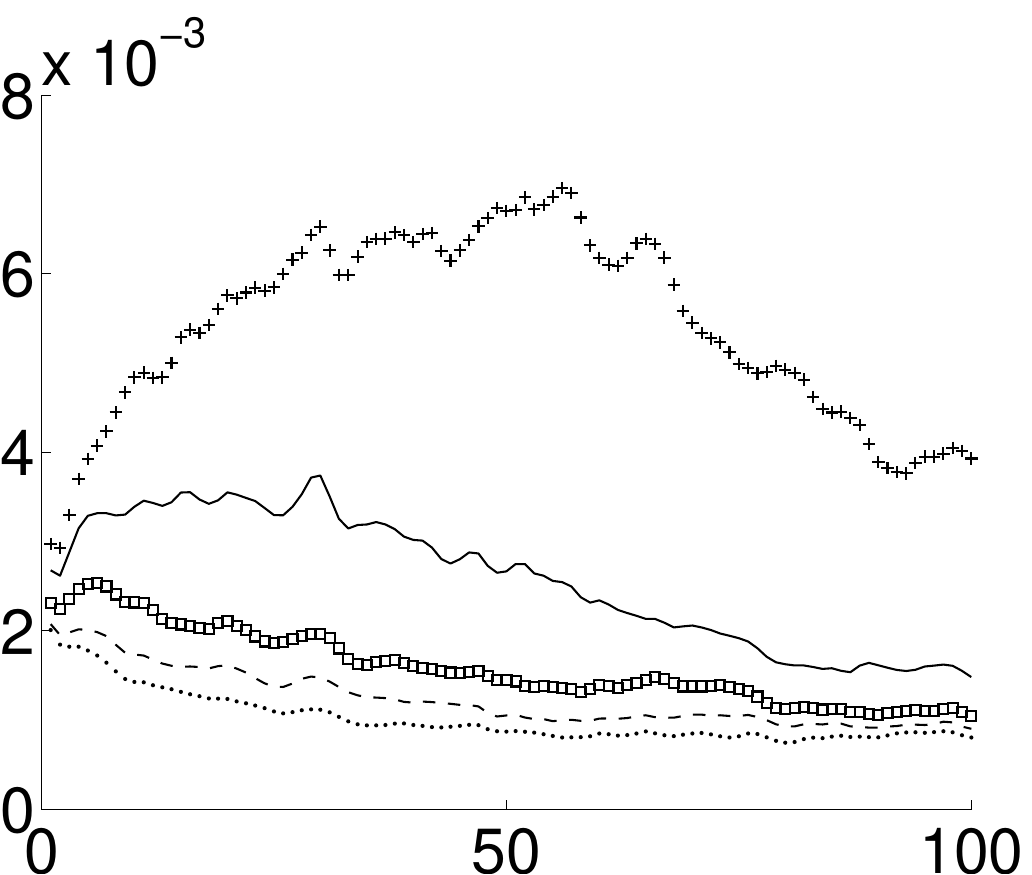} 
}
\caption{Analysis RMSE for the Lorenz model with 500 of variables. Different curves correspond to  different numbers of ensemble members: $200(+), 250(-), 300(\square), 350(--)$ and $400(.)$. When the number of ensemble members is increased the analysis is improved. }
\label{Fig:Lorenz500-Results-Simulation}
\end{figure}

Figures  \ref{Fig:Lorenz500-Results-Simulation} and  \ref{Fig:Lorenz1000-Results-Simulation} show the RMSE decrease over the assimilation window for $\Nstate =$ 500 and 1000, respectively. When the number of ensemble members is increased the analysis errors are smaller, as expected. There is no significant difference in results between different implementations of the EnKF.

\begin{figure}[H]
\centering
\begin{sideways} $\rm RMSE \left( \XA \right )$ \end{sideways} 
\subfigure[Sherman]{
\includegraphics[width=0.25\textwidth]{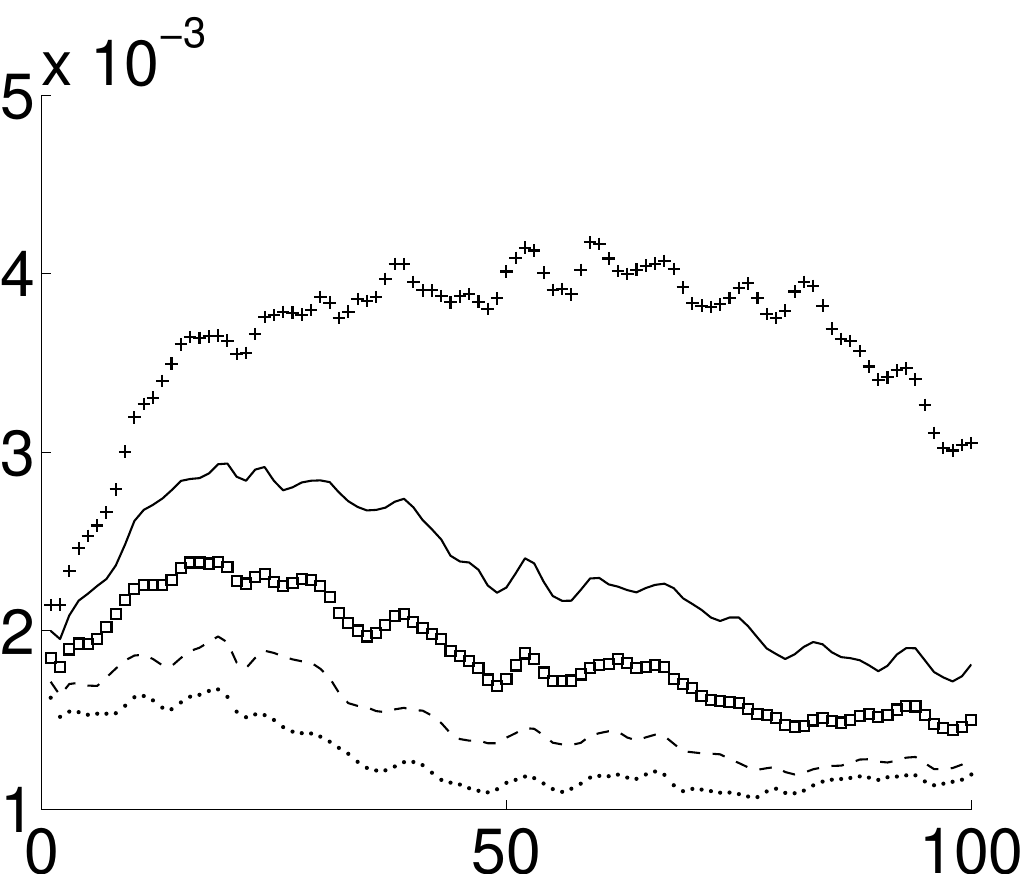} 
}
\subfigure[Cholesky]{
\includegraphics[width=0.25\textwidth]{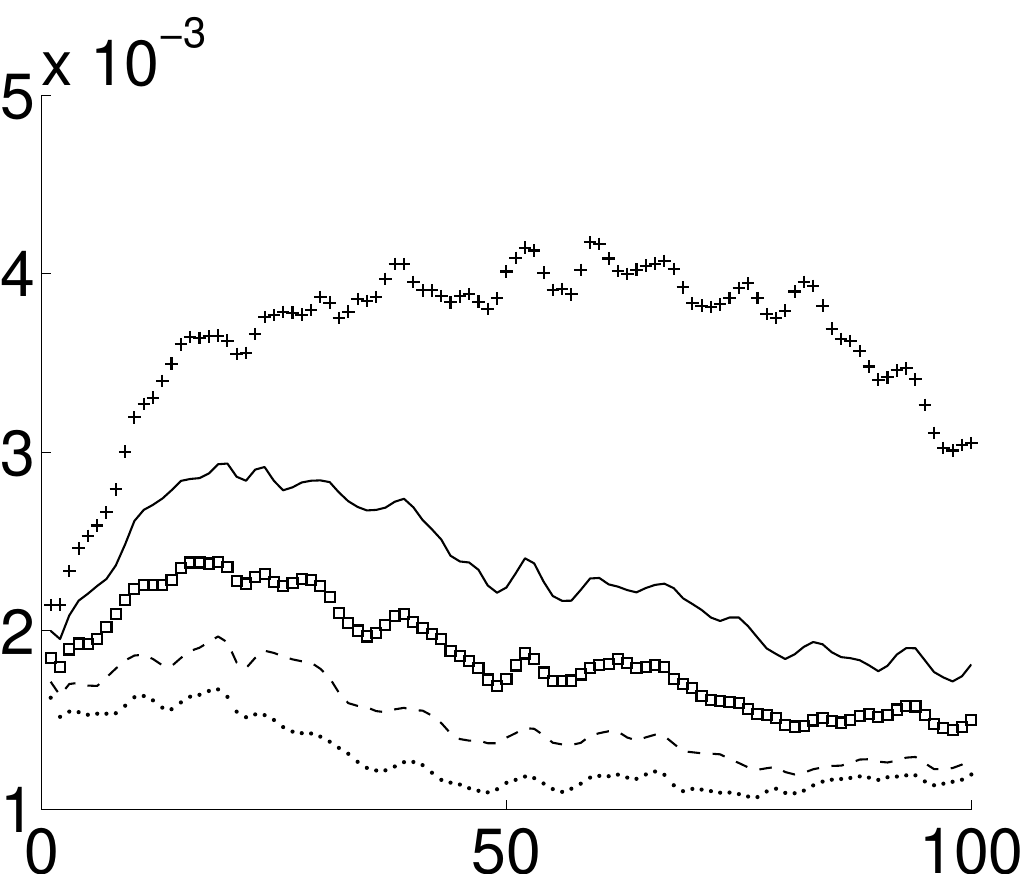}
}
\subfigure[SVD]{
\includegraphics[width=0.25\textwidth]{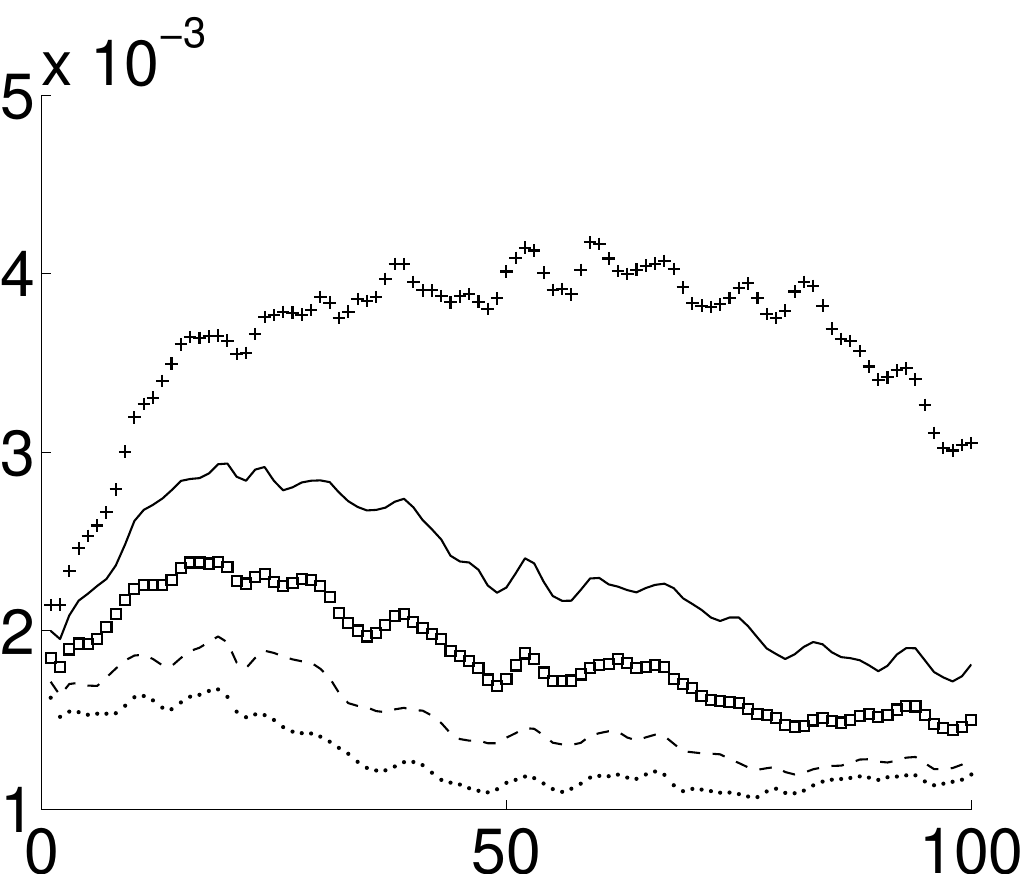} 
}
\caption{Analysis RMSE for the Lorenz model with 1000 of variables. Different curves correspond to  different numbers of ensemble members: $200(+), 250(-), 300(\square), 350(--)$ and $400(.)$. When the number of ensemble members is increased the analysis is improved. }
\label{Fig:Lorenz1000-Results-Simulation}
\end{figure}

The ET results are shown in Table  \ref{Tab:Lorenz-Results-TIME}.
The Cholesky decomposition is the most efficient for a small number of observations and states. 
When the number of observations is increased, the relative performance of Cholesky deteriorates, as expected from the complexity results presented in Section  \ref{intro}. The Cholesky decomposition solution of the linear system \eqref{EnKF:system-to-solve} is not suitable when the number of observations is large. The SVD implementation exhibits a good performance for a small number of ensemble members and observations. However, the ET of the SVD implementation grows faster than that of the Cholesky implementation when the number of ensemble and/or observations are increased, due to the term $\Nens^3 \sim \Nobs^3$ in its complexity formula. Tthe EnKF implementation based on SVD is not suitable  for a large  number of observations or a large number of ensemble members. Finally, the Sherman-Morrison implementation has the best performance for a large number of observations and states. This implementation is suitable for a large number of observations. Since the term $\Nens^3$ does not appear in the cost upper-bound of the iterative-Sherman implementation, when $\Nens \sim \Nobs$, the proposed implementation will exhibit a better performance than those implementations presented in ~\cite{Anderson01,Anderson07,Tippett2003} since they are upper-bounded by  \eqref{EnKF:SVD-overall-complexity} (see table  \ref{Tab:Analysis-Complexity-Contrast}). 

\begin{table}[H]
\centering
{
\begin{tabular}{|c|r|r|r|r|} \hline
${\bf N_{obs}/ N_{state}}$ & $\bf N_{ens}$ & ${\bf EnKF_{Sher}}$ & ${\bf  EnKF_{Chol}}$ & ${\bf EnKF_{SVD} }$  \\ \hline
\multirow{5}{*}{500/500} & 200 & $46.5$ s & $36.9$ s & $66.3$ s \\ 
&   250 & $72.1$ s & $45.1$ s & $91.4$ s \\ 
&   300 & $103.5$ s & $54.2$ s & $118.9$ s \\ 
&   350 & $140.8$ s & $64.4$ s & $138.8$ s \\ 
&   400 & $183.6$ s & $75.4$ s & $174.5$ s \\ \hline
\multirow{5}{*}{1000/1000} & 400 & $377.6$ s & $327.4$ s & $ 687.1$ s \\ 
&   450 & $500.1$ s & $364.6$ s & $715.9$ s \\ 
&   500 & $601.7$ s & $403.5$ s & $ 877.3$  s \\ 
&   550 & $731.8$ s & $462.4$ s & $1038.3$  s \\ 
&   600 & $873.2$ s & $492.4$ s & $1199.3$ s \\ \hline
\multirow{5}{*}{3000/3000} & 900 & $1.7$ h & $2.3$ h & $3.8$ h  \\ 
&   950 & $2.1$ h & $2.4$ h & $3.9$ h \\ 
&   1000 & $2.2$ h & $2.5$ h & $4.1$ h  \\ 
&   1050 & $2.4$ h & $2.7$ h & $4.4$ h  \\ 
&   1100 & $2.6$ h & $2.8$ h & $4.8$ h \\ \hline
\multirow{3}{*}{5000/5000} & 1500 & $8.1$ h & $11.0$ h & $17.3$ h  \\ 
&   1550 & $8.8$ h & $11.2$ h & $17.8$ h\\ 
&   1600 & $9.2$ h & $11.5$ h & $19.9$ h  \\ \hline
\end{tabular}
}
\caption{Computational times for the Lorenz model assimilation performed with different EnKF implementations. 
The Cholesky decomposition is the most efficient for a small number of observations and states. The
Sherman-Morrison implementation is the bet for a large number of observations and states. 
}
\label{Tab:Lorenz-Results-TIME}
\end{table}

\subsection{Quasi-geostrophic model ($\bf N_{obs} \gg \N_{ens}$)}

The Earth's ocean has a complex flow system influenced by the rotation of the Earth, the density stratification due to temperature and salinity, as well as other factors. The quasi-geostrophic (QG) model is a simple model which mimics the real behavior of the ocean. 
It is defined by the following partial differential equations \cite{Carton94}:

\label{Results:QG-Model}
\begin{eqnarray}
\label{QGModel:model} \nonumber
\displaystyle \frac{\partial q}{\partial t}+r \cdot {\mathcal{J}} \left ( {\bf \psi} , q\right ) + \beta \cdot \frac{\partial \psi}{\partial x} &=&  -rkb \cdot \zeta + rkh \cdot \nabla^2 \zeta - rkh2 \cdot \nabla^4 \zeta \\ 
&+& \underbrace{\sin \left( 2 \cdot \pi \cdot y\right )}_{ \text{ External Force}} \, ,
\end{eqnarray}
where 
\begin{eqnarray} 
\displaystyle {\mathcal{J}} \left ( {\bf \psi} , q\right ) = \frac{\partial q}{\partial x} \cdot \frac{\partial \psi}{\partial y}-\frac{\partial q}{\partial y} \cdot \frac{\partial \psi}{\partial x} \, ,
\end{eqnarray}
$q = \zeta - {\rm F} \cdot \psi$ is the potential vorticity, $\psi$ is the stream function, ${\rm F}$ is the Froud number, $\zeta = \nabla^2 \psi$ is the relative vorticity, $r$ is a sort of the Rossby number, $rkb$ is the bottom friction, $rkh$ is the horizontal friction and $rhk2$ is the biharmonic horizontal friction and $x$ and $y$ represent the horizontal and vertical components of the space.

Moreover, $q$ and ${\psi}$ are related to one another through an elliptic operator \cite{Pedl96}:
\begin{eqnarray}
\label{QGModel:relation-vorticity-stream}
\displaystyle q = \tilde{\nabla^2} \psi,
\end{eqnarray}
which yields 
\begin{eqnarray}
\label{QGModel:inversion-vorticity}
\displaystyle  \tilde{\nabla^{-2}} q = \psi\, ,
\end{eqnarray}
where 
\begin{eqnarray*} 
\tilde{\nabla^{2}} = \frac{\partial^2}{\partial x^2}+\frac{\partial^2}{\partial y^2} \,.
\end{eqnarray*}

This elliptic property reflects the assumption that the flow is geostrophically balanced in the horizontal direction, and hydrostatically balanced in the vertical direction.

The QG experiment studies the behavior of EnKF implementations when $\Nobs \gg \Nens$ (the number of observations is much larger than the number of ensemble members) as is usually the case in practice. Moreover, this scenario is more difficult than the previous one (the Lorenz model): large model-errors are considered in the initial ensemble members. Besides, data is available every 10 time units. 

We consider three different grids, denoted QGNM, where the number of horizontal and vertical grid points are $ \rm N$ and $\rm M$, respectively. Specifically, we employ in experiments QG33 (small instance),  QG65 (medium instance) and QG129 (large instance). 
The horizontal and vertical dimensions of the grid are denoted by $\rm L_x$ and $\rm L_y$ respectively.  These instances and the corresponding parameter values are summarized in Table  \ref{Tab:QG-Instances}.
\begin{table}[H]
\centering
{ \footnotesize
\begin{tabular}{|c|c|c|c|c|c|c|c|c|c|} \hline
Instance & ${\rm L_x}$ & ${\rm L_y}$ & ${\rm N}$ & ${\rm M}$ & ${\rm rkb}$ & ${\rm rkh}$ & ${\rm rkh2}$ & ${\rm \beta}$ & $r$ \\ \hline
QG33 & 0.4 & 0.4 & 33 & 33 & $10^{-6}$ & $10^{-7}$ & $2 \times 10^{-12}$ & 1.0 & $10^{-5}$ \\ \hline
QG65 & 1.0 & 1.0 & 65 & 65 & $10^{-6}$ & $10^{-7}$ & $2 \times 10^{-12}$ & 1.0 & $10^{-5}$  \\ \hline
QG129 & 1.0 & 1.0 & 129 & 129 & $10^{-6}$ & $10^{-7}$ & $2 \times 10^{-12}$ & 1.0 & $10^{-5}$  \\ \hline
\end{tabular}
}
\caption{Parameter values for the QG model instances considered. $\rm L_x$ and $\rm L_y$ represent the horizontal and vertical grid sizes, and $\rm N$ and $\rm M$ are the number of horizontal and vertical grid points, respectively.}
\label{Tab:QG-Instances}
\end{table}

The experimental settings are described below.
\begin{itemize}
\item There are 1200 time steps, each of one representing 1.27 days in the ocean.
\item The vorticity of the ocean at each grid point provides a component of the vector state. 
\item The computation of the stream function is done through the solution of the Helmholtz \cite{Otto99} function according to the elliptic property \eqref{QGModel:inversion-vorticity}.
\item Homogeneous Dirichlet boundary conditions are assumed. Due to this, the boundaries of the grid are not mapped into the state vector, and $\Nstate = \rm (N-2) \cdot (M-2)$.
\item The initial ensemble members are constructed as follows:
\begin{eqnarray*}
\displaystyle \xb{i} = \xt + {\varepsilon_{\rm i}^{\rm B}} \cdot  \underbrace{\left (\frac{1}{\Nstate} \cdot \sum_{\rm k = 1}^{\Nstate} \left |  \xt_{\rm k} \right | \right )}_{\bf C} \in \Re^{\Nstate \times 1} \, , \text{ for $1 \le \rm i \le \Nens$}\, ,
\end{eqnarray*}
where ${\varepsilon^{\rm B}}$ is drawn from a Normal distribution with zero mean and covariance matrix
\begin{eqnarray} \nonumber
\displaystyle \Q = {\bf diag} \left \{ \rm STD_{ens}^2 \right \} \in \Re^{\Nstate \times \Nstate} \, .
\end{eqnarray}

For testing purposes, three values are assumed for the standard deviation of model errors ($\rm STD_{ens}$): 2.5, 5.0 and 7.5. Notice the large dispersion of the initial ensemble members, which can make difficult the convergence of any filter since the huge spread out of the initial ensemble members with respect to the typical value $\bf C$. 

\item The number of observation per simulation, for each size ($\Nstate$) of the model state , is defined as follows:
\begin{eqnarray} \nonumber
\displaystyle \Nobs = {\rm P_{obs}} \cdot \Nstate,
\end{eqnarray}
where ${\rm P_{obs}}$ is the percentage of components observed from the model state. The values given to ${\rm P_{obs}}$ are 50\%, 70\% and 90\%. Those, measurements are taken every 10 time units and they are constructed as shown in equation \eqref{eq:syntethic-observations}. Notice, there are 120 analysis steps out of 1200 time steps (10\% of the total simulation time).

\item For the time evolution of the model, zero boundary conditions are assumed and the boundaries are not included onto the ensemble representation.  Due to this, the dimension of the vector state ${\Nstate} = \rm (N-2) \cdot (M-2)$.

\item For each instance we consider simulations with $\Nens \in \{ 20, 60, 100\}$ ensemble members. The number of ensemble members is one to two orders of magnitude smaller than the total number of observations.

%
\end{itemize}

The RSME values for analysis errors for the QG33, QG65 and QG129 instances are shown in Tables  \ref{Tab:QG33-Results-RMSE}, \ref{Tab:QG65-Results-RMSE} and  \ref{Tab:QG129-Results-RMSE}, respectively. The results depend on the number of ensemble members ($\Nens$), the number of observations ($\Nobs$), and the deviation of the initial ensemble mean ($\rm STD_{ens}$). The RSME is quantifiess errors in the stream function $\psi$, whose values are computed through the relation \eqref{QGModel:inversion-vorticity}.  In terms of accuracy there is no significant difference between different EnKF implementations. As expected, when the error in the initial ensemble is increased, the accuracy in the analysis decreases. The error does not show an exponential growth, even when the number of components in the model state ($\Nstate$) is much larger than the number of ensemble members (e.g., for the QG129 instance). When the number of ensemble members is increased, the analysis error is decreased. This is illustrated by the snapshots of the QG33 simulation over 1200 time steps presented in Figure  \ref{Fig:QG33-Snapshots-Simulation}. There, we can clearly see that the ensemble of size 100 provides a better estimation ($\xmean$) to the true state of the model ($\xt$) than the ensembles of sizes 20 and 60. Additionally, the number of observations plays an important role in the estimation of the true model state when the size of the vector state is much larger than the number of ensemble members. 

The ET values for the QG33, QG65 and QG129 simulations are shown in Tables  \ref{Tab:QG33-Results-RMSE}, \ref{Tab:QG65-Results-RMSE} and  \ref{Tab:QG129-Results-RMSE}, respectively. The time is expressed in seconds (s)
if it is below 30 minutes, and otherwise is expressed
in minutes (min) and hours (h). The Cholesky implementation shows good performance when the number of observations is small. From Table  \ref{Tab:QG33-Results-ElapsedTime} (the blocks where the number of observations are 480, 672 and 864) we see that the Cholesky implementation performance is more sensitive to the number of observations than to the number of ensemble members. 
This EnKF implementation is not suitable for a large number of observations. For instance, the Cholesky elapsed time for the QG129 instance is not presented since each simulation takes more than 4 days in order to be completed. 

The SVD implementation shows a better relative performance than for the Lorenz 96 test, since the number of ensemble members is small with respect to the number of observations. For example, for the QG33 instance, the SVD implementation shows a better performance than Cholesky when the number of observations and ensemble members are small. In addition, when the size of vector state is increased, the SVD implementation shows a better performance than Cholesky. This agrees with the computational complexity upper bounds presented  in Section  \ref{intro}. As is expected, the performance of the SVD based methods is better than the Cholesky implementations when the number of observations is much larger than the number of the ensemble members. 

The Sherman-Morrison implementation shows the best performance among the compared methods.
This is true even when the number of observations is much larger than the number of ensemble members, as seen in Table  \ref{Tab:QG129-Results-ElapsedTime}. The results of both test cases (the quasigeostrophic and Lorenz models) lead to the conclusion that the performance of the iterative-Sherman implementation is not sensitive to the increase in the number of observations, making it attractive for implementation with large-scale observational systems.

\begin{table}[H]
\centering
{\footnotesize
\begin{tabular}{|c|c|c|c|c|c|} \hline
$\bf N_{ens}$ & $\bf N_{obs}$ & $\bf STD_{ens}$ & ${\bf  EnKF_{Sher} }$ & ${\bf EnKF_{Chol}  }$ & ${\bf  EnKF_{SVD} }$ \\ \hline

\multirow{9}{*}{20} &  \multirow{3}{*}{480}  & 2.5 & $1.71348538 \times 10^{-4}$ & $1.71348538 \times 10^{-4} $ & $ 1.71348538 \times 10^{-4} $ \\
& & 5.0 & $ 3.42503478 \times 10^{-4} $ & $ 3.42503478 \times 10^{-4} $ & $ 3.42503478 \times 10^{-4} $  \\ 
& & 7.5 & $ 5.13685273 \times 10^{-4} $ & $ 5.13685273 \times 10^{-4} $ & $ 5.13685273 \times 10^{-4} $  \\ 
\cline{2-6}
&  \multirow{3}{*}{672}  & 2.5 & $ 1.68916683 \times 10^{-4} $ & $ 1.68916683 \times 10^{-4} $ & $ 1.68916683 \times 10^{-4} $ \\
& & 5.0 & $ 3.37557073 \times 10^{-4} $ & $ 3.37557073 \times 10^{-4} $ & $ 3.37557073 \times 10^{-4} $  \\ 
& & 7.5 & $ 5.06204927 \times 10^{-4} $ & $ 5.06204927 \times 10^{-4} $ & $ 5.06204927 \times 10^{-4} $  \\ 
\cline{2-6}
&  \multirow{3}{*}{864}  & 2.5 & $ 1.68758284 \times 10^{-4} $ & $ 1.68758284 \times 10^{-4} $ & $ 1.68758284 \times 10^{-4} $ \\
& & 5.0 & $ 3.37604000 \times 10^{-4} $ & $ 3.37604000 \times 10^{-4} $ & $ 3.37604000 \times 10^{-4} $  \\ 
& & 7.5 & $ 5.06437992 \times 10^{-4} $ & $ 5.06437992 \times 10^{-4} $ & $ 5.06437992 \times 10^{-4} $  \\ 
\cline{1-6}
\multirow{9}{*}{60} &  \multirow{3}{*}{480}  & 2.5 & $1.71410180 \times 10^{-4}$ & $1.71410180 \times 10^{-4} $ & $ 1.71410180 \times 10^{-4} $ \\
& & 5.0 & $ 3.43004244 \times 10^{-4} $ & $ 3.43004244 \times 10^{-4} $ & $ 3.43004244 \times 10^{-4} $  \\ 
& & 7.5 & $ 5.14603943 \times 10^{-4} $ & $ 5.14603943 \times 10^{-4} $ & $ 5.14603943 \times 10^{-4} $  \\ 
\cline{2-6}
&  \multirow{3}{*}{672}  & 2.5 & $ 1.64430970 \times 10^{-4} $ & $ 1.64430970 \times 10^{-4} $ & $ 1.64430970 \times 10^{-4} $ \\
& & 5.0 & $ 3.29692664 \times 10^{-4} $ & $ 3.29692664 \times 10^{-4} $ & $ 3.29692664 \times 10^{-4} $  \\ 
& & 7.5 & $ 4.94948120 \times 10^{-4} $ & $ 4.94948120 \times 10^{-4} $ & $ 4.94948120 \times 10^{-4} $  \\ 
\cline{2-6}
&  \multirow{3}{*}{864}  & 2.5 & $ 1.64737541 \times 10^{-4} $ & $ 1.64737541 \times 10^{-4} $ & $ 1.64737541 \times 10^{-4} $ \\
& & 5.0 & $ 3.29861914 \times 10^{-4} $ & $ 3.29861914 \times 10^{-4} $ & $ 3.29861914 \times 10^{-4} $  \\ 
& & 7.5 & $ 4.94925540 \times 10^{-4} $ & $ 4.94925540 \times 10^{-4} $ & $ 4.94925540 \times 10^{-4} $  \\ 
\cline{1-6}
\multirow{9}{*}{100} &  \multirow{3}{*}{480}  & 2.5 & $1.62358824 \times 10^{-4}$ & $1.62358824 \times 10^{-4} $ & $ 1.62358824 \times 10^{-4} $ \\
& & 5.0 & $ 3.23304631 \times 10^{-4} $ & $ 3.23304631 \times 10^{-4} $ & $ 3.23304631 \times 10^{-4} $  \\ 
& & 7.5 & $ 4.84328422 \times 10^{-4} $ & $ 4.84328422 \times 10^{-4} $ & $ 4.84328422 \times 10^{-4} $  \\ 
\cline{2-6}
&  \multirow{3}{*}{672}  & 2.5 & $ 1.54859388 \times 10^{-4} $ & $ 1.54859388 \times 10^{-4} $ & $ 1.54859388 \times 10^{-4} $ \\
& & 5.0 & $ 3.10925745 \times 10^{-4} $ & $ 3.10925745 \times 10^{-4} $ & $ 3.10925745 \times 10^{-4} $  \\ 
& & 7.5 & $ 4.66772642 \times 10^{-4} $ & $ 4.66772642 \times 10^{-4} $ & $ 4.66772642 \times 10^{-4} $  \\ 
\cline{2-6}
&  \multirow{3}{*}{864}  & 2.5 & $ 1.44950729 \times 10^{-4} $ & $ 1.44950729 \times 10^{-4} $ & $ 1.44950729 \times 10^{-4} $ \\
& & 5.0 & $ 2.90313729 \times 10^{-4} $ & $ 2.90313729 \times 10^{-4} $ & $ 2.90313729 \times 10^{-4} $  \\ 
& & 7.5 & $ 4.35737112 \times 10^{-4} $ & $ 4.35737112 \times 10^{-4} $ & $ 4.35737112 \times 10^{-4} $  \\ 
\hline
\end{tabular}
}
\caption{Analysis RMSE for different EnKF implementations applied to the QG33 instance.
All methods give similar results. When the number of ensemble and/or observations is increased, the analysis accuracy is improved.}
\label{Tab:QG33-Results-RMSE}
\end{table}

\begin{figure}[H]
\centering
\begin{tabular}{VNNNN}
& $\xt_{\rm t} $ & $\xmeana_{\rm t}$, $\Nens = 20$ & $\xmeana_{\rm t}$, $\Nens = 60$ & $\xmeana_{\rm t}$, $\Nens = 100$\\ 
\rotatebox{90}{${\rm t = 0}$} &\includegraphics[width=0.2\textwidth]{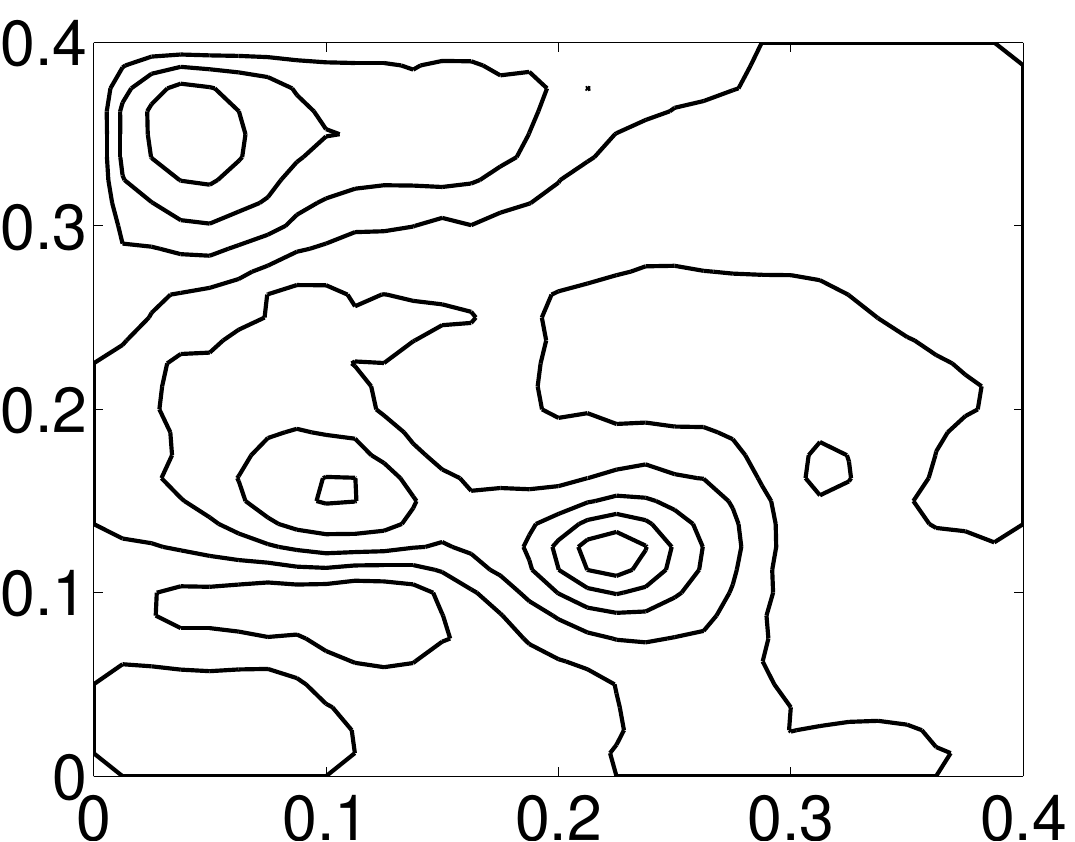} &\includegraphics[width=0.2\textwidth]{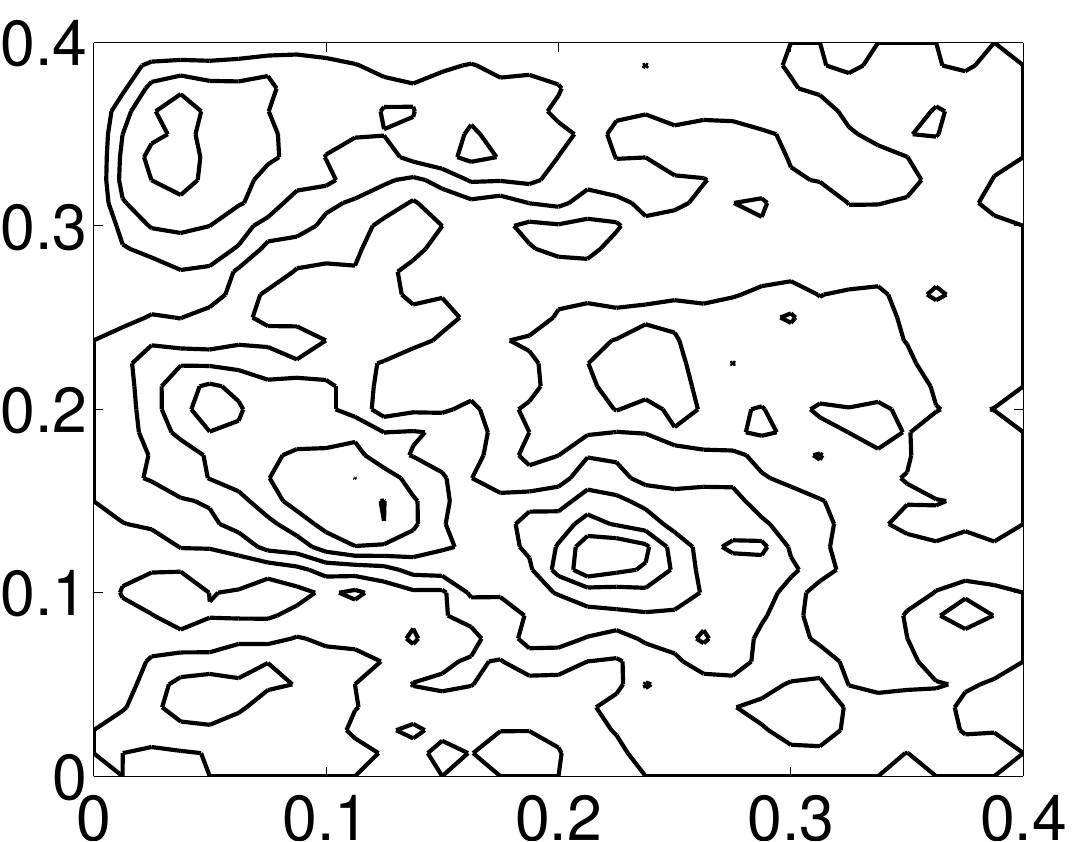} &\includegraphics[width=0.2\textwidth]{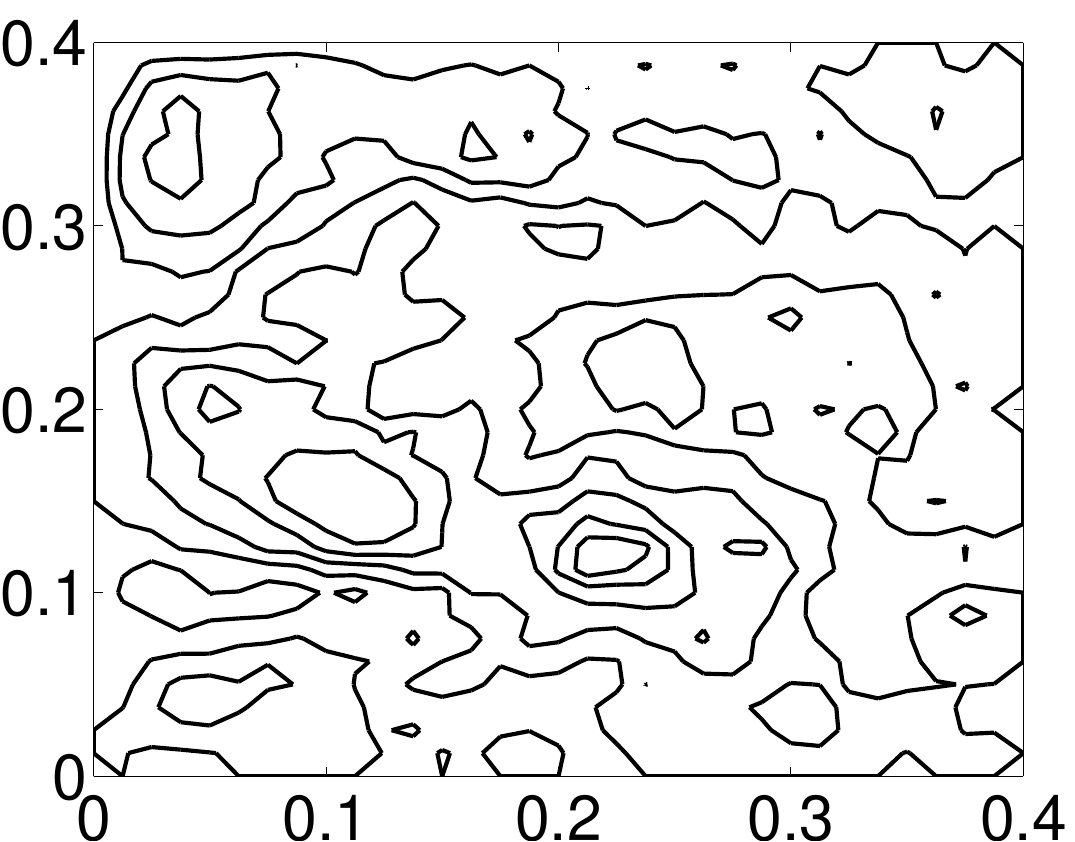} &\includegraphics[width=0.2\textwidth]{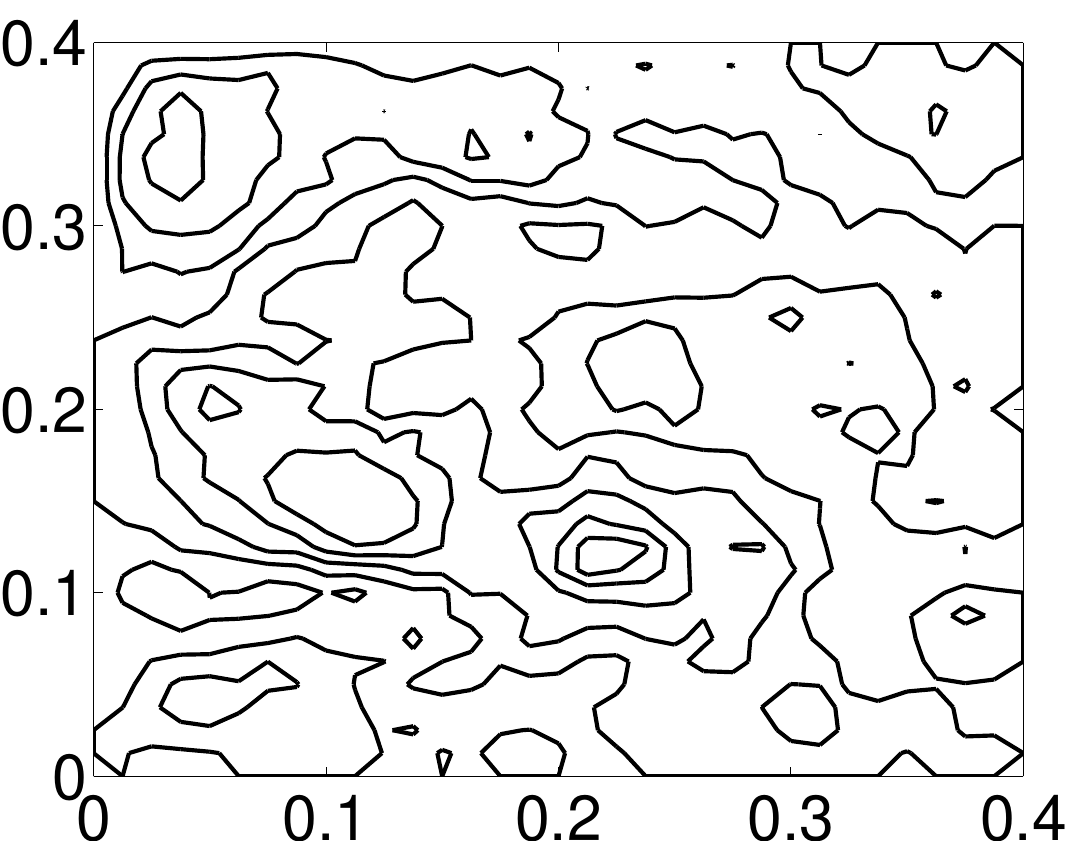} \\
\begin{sideways} ${\rm t = 239}$ \end{sideways} &\includegraphics[width=0.2\textwidth]{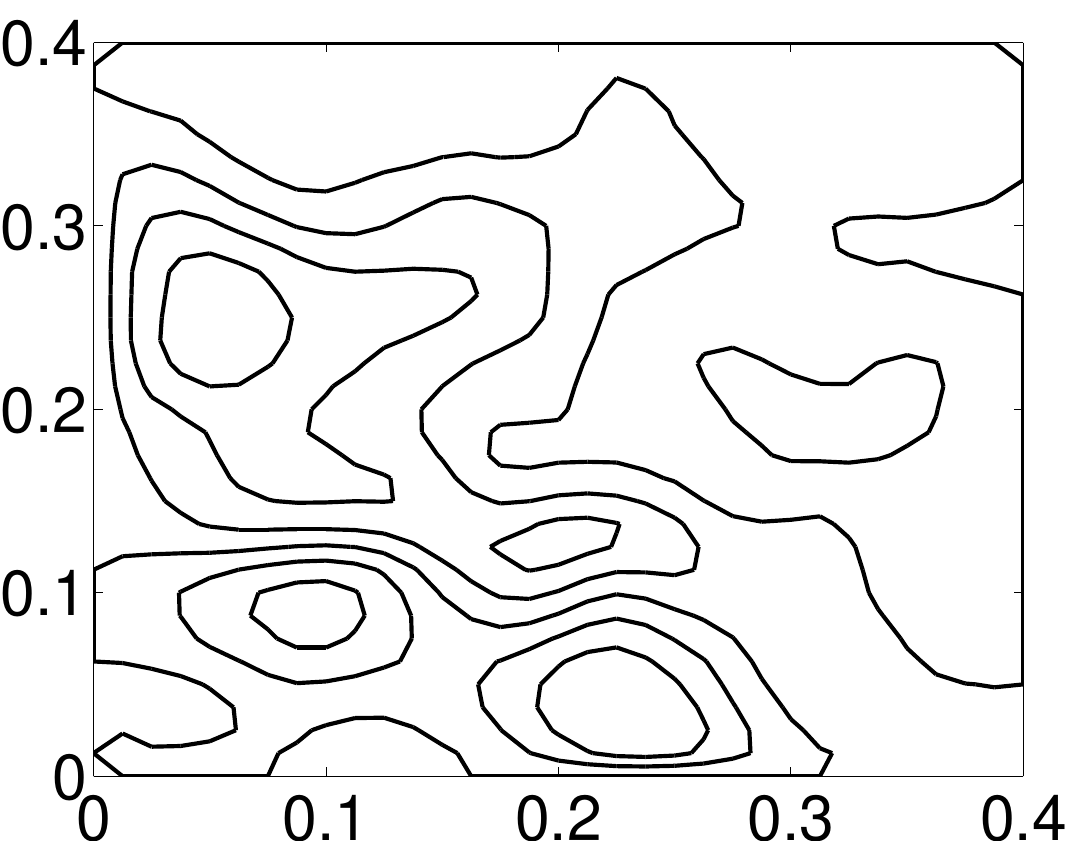} &\includegraphics[width=0.2\textwidth]{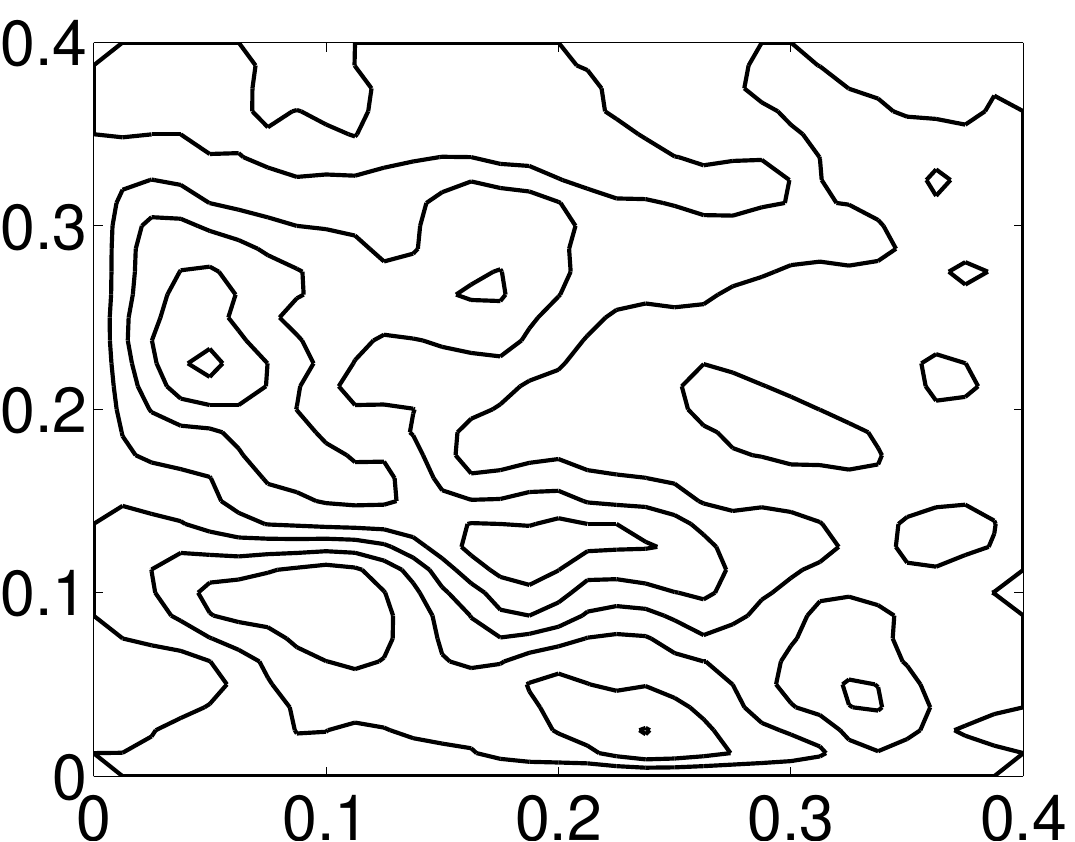} &\includegraphics[width=0.2\textwidth]{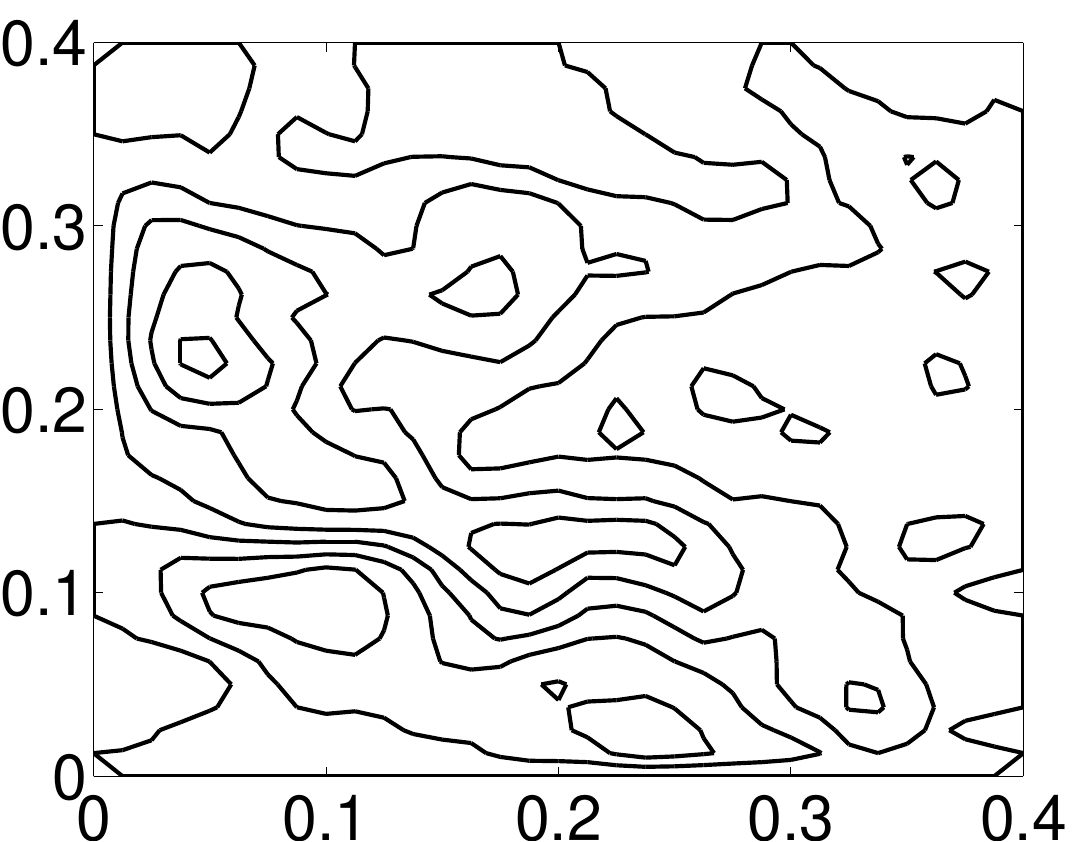} &\includegraphics[width=0.2\textwidth]{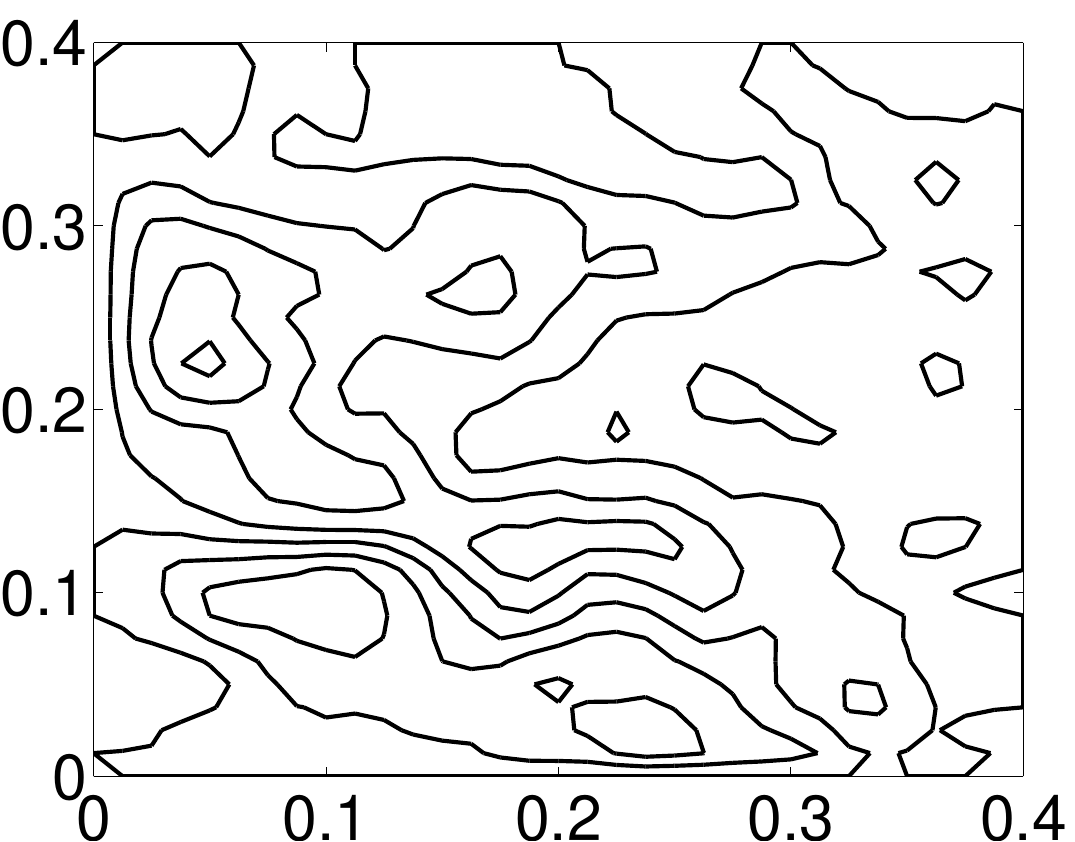} \\
\begin{sideways} ${\rm t = 478}$ \end{sideways} &\includegraphics[width=0.2\textwidth]{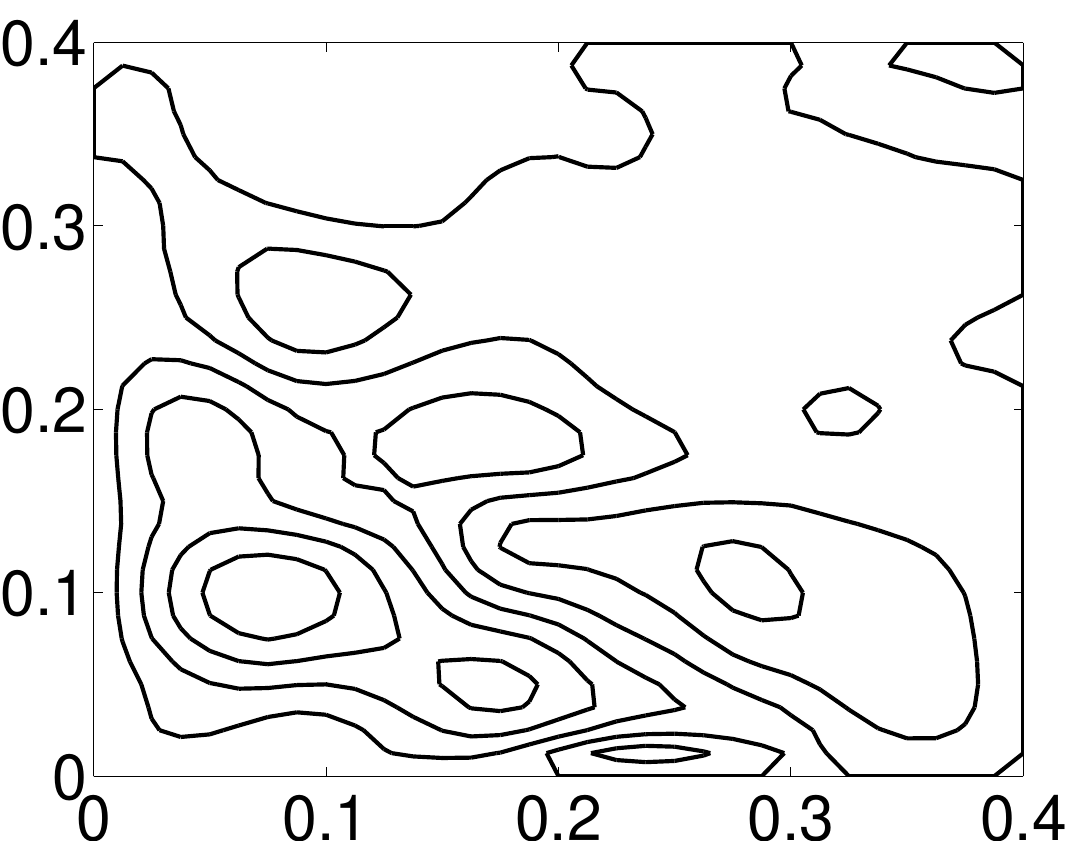} &\includegraphics[width=0.2\textwidth]{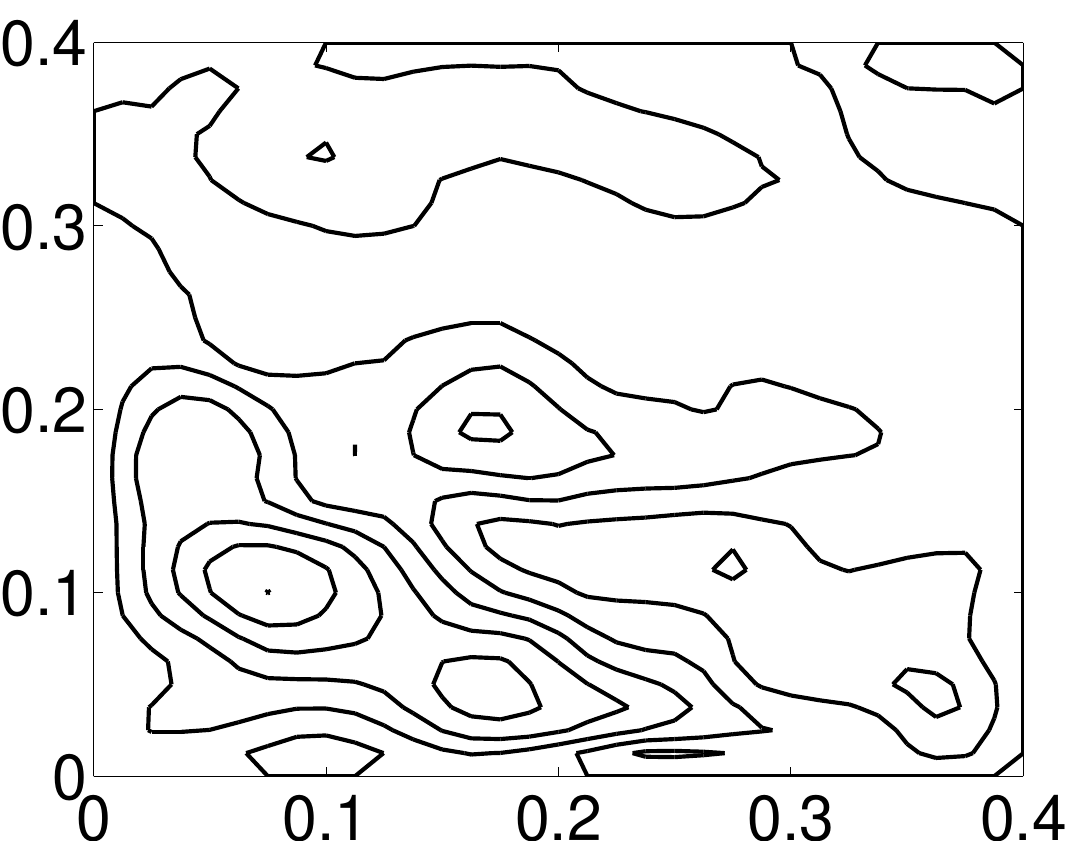} &\includegraphics[width=0.2\textwidth]{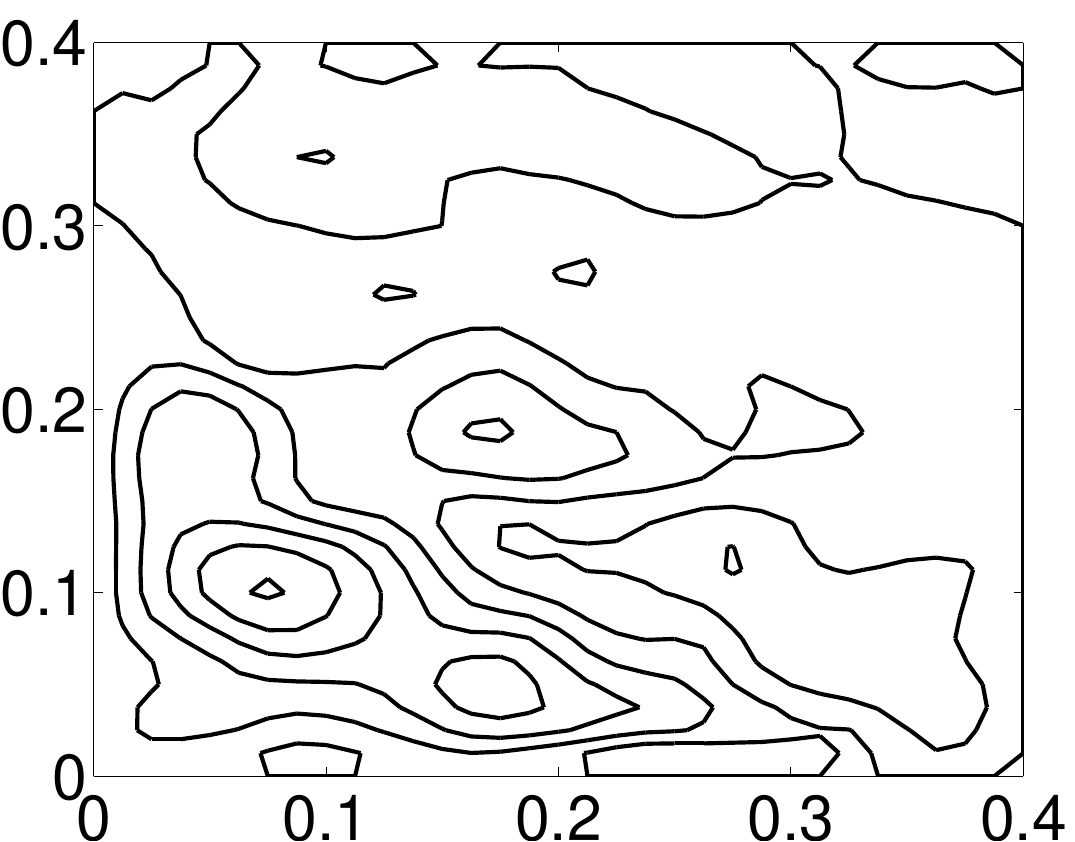} &\includegraphics[width=0.2\textwidth]{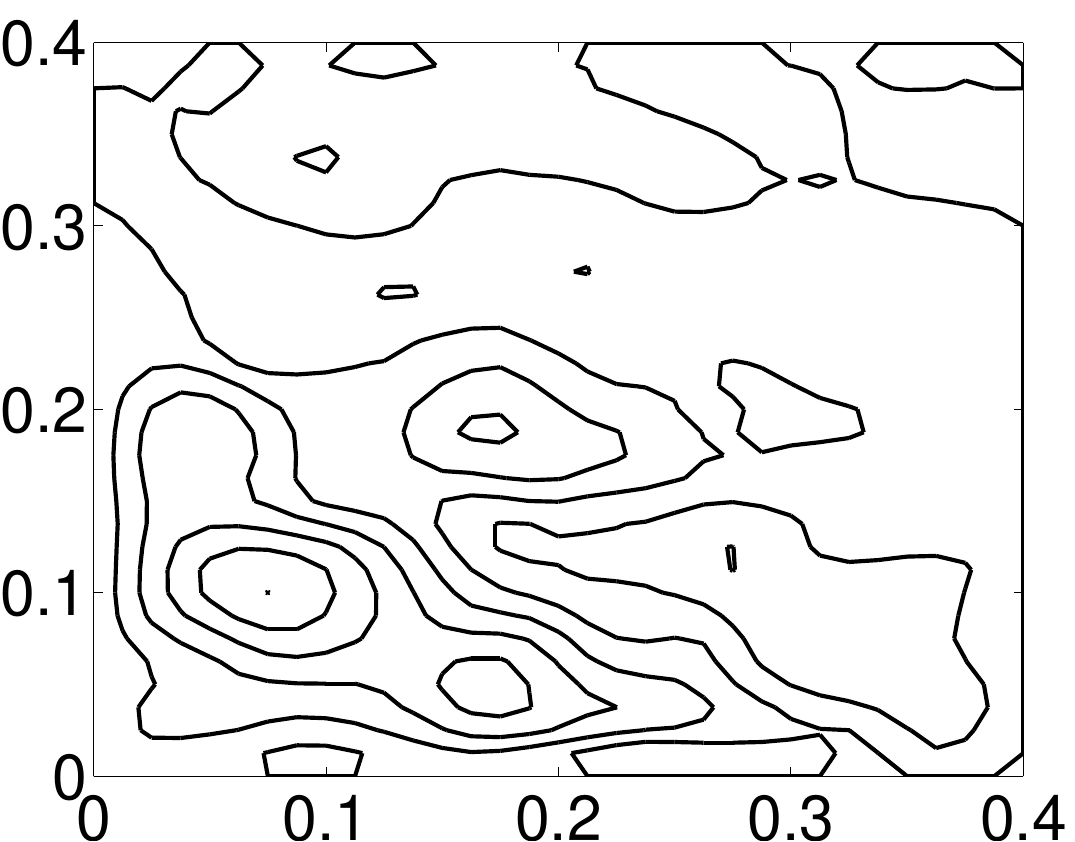} \\
\begin{sideways} ${\rm t = 717}$ \end{sideways} &\includegraphics[width=0.2\textwidth]{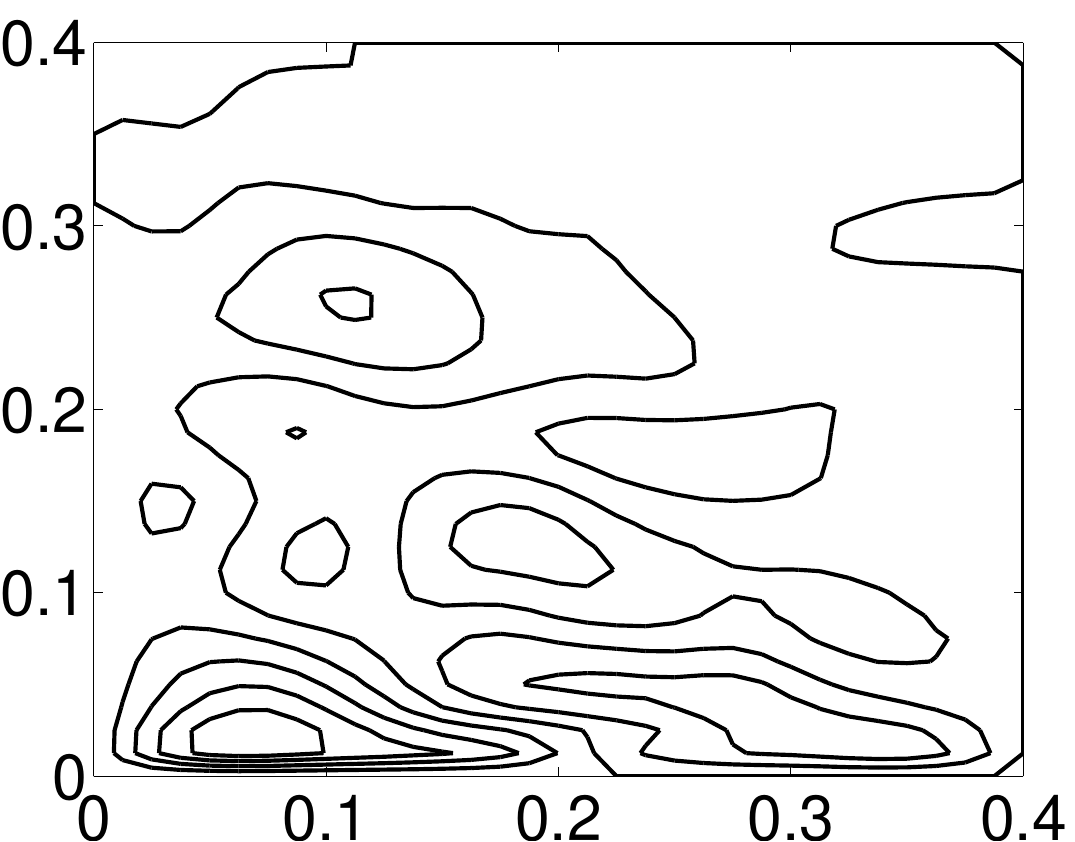} &\includegraphics[width=0.2\textwidth]{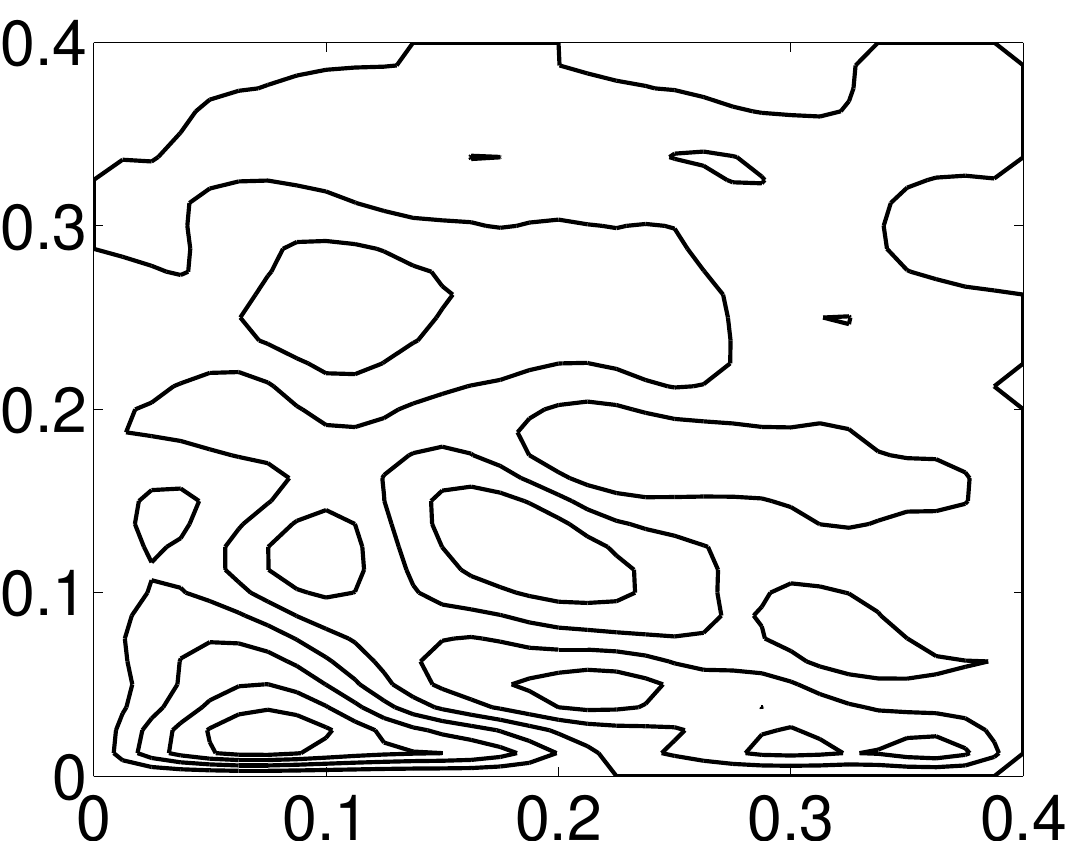} &\includegraphics[width=0.2\textwidth]{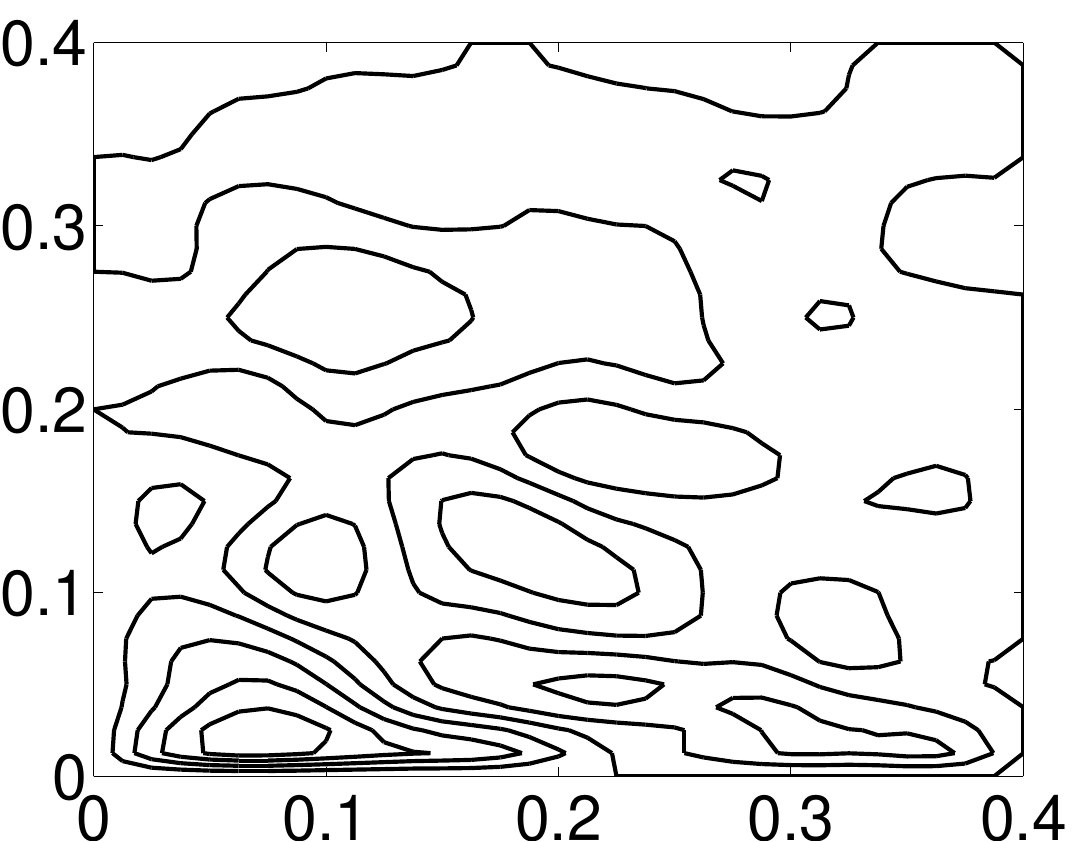} &\includegraphics[width=0.2\textwidth]{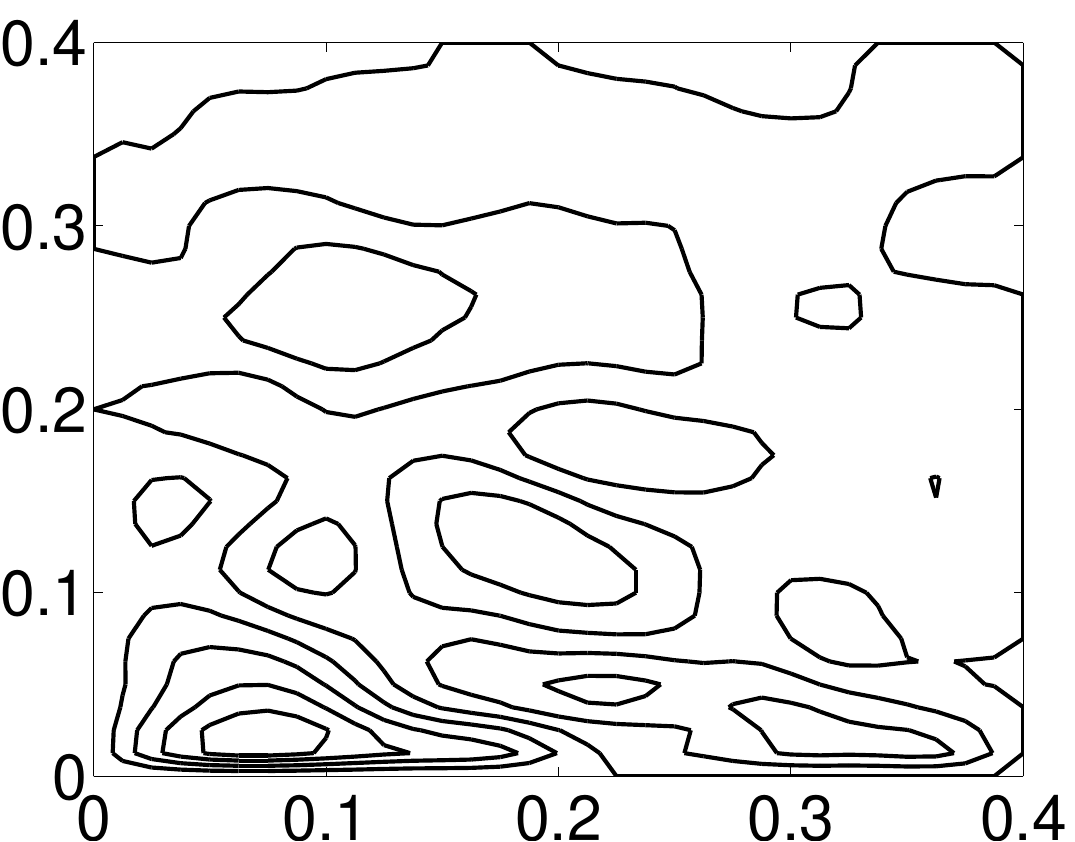} \\
\begin{sideways} ${\rm t = 956}$ \end{sideways} &\includegraphics[width=0.2\textwidth]{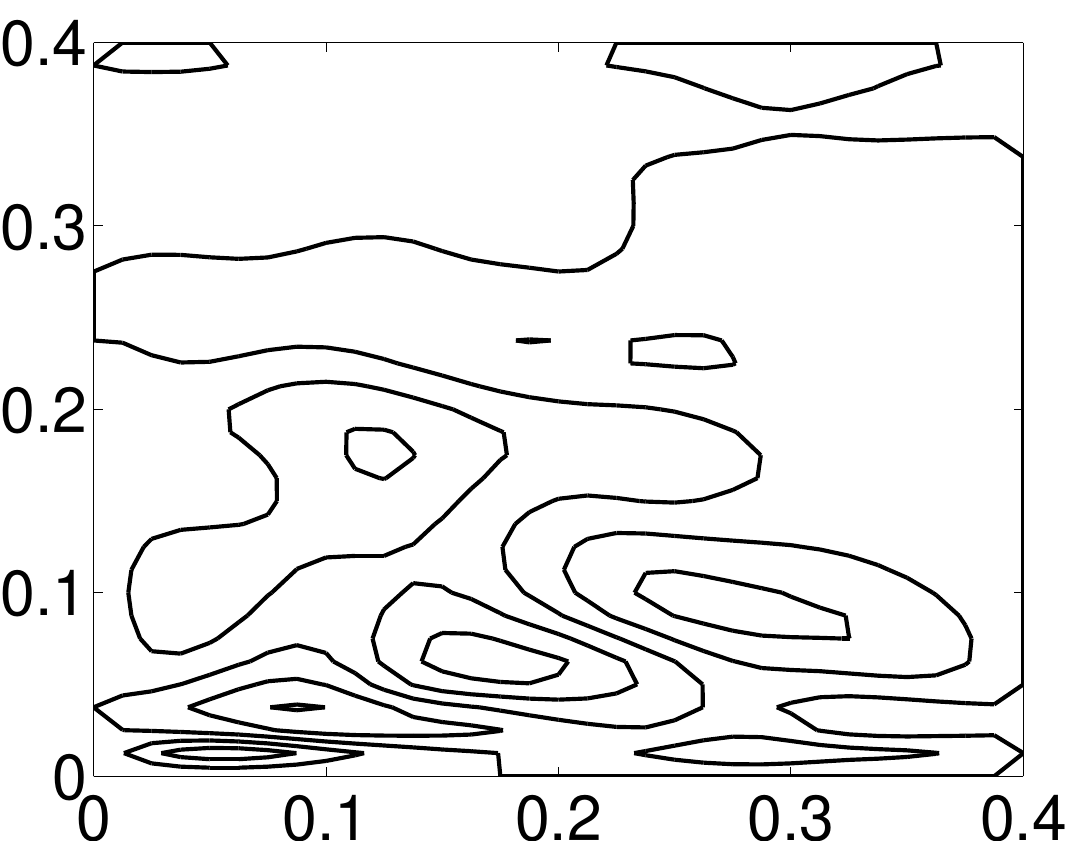} &\includegraphics[width=0.2\textwidth]{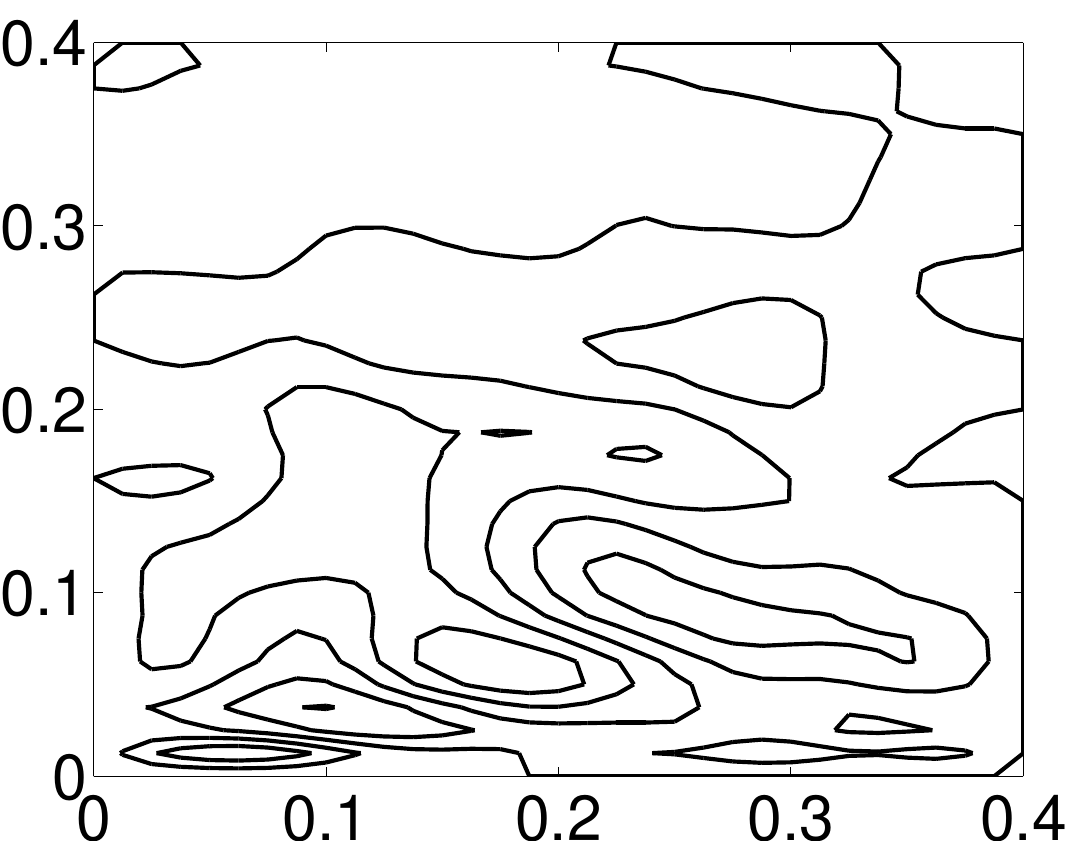} &\includegraphics[width=0.2\textwidth]{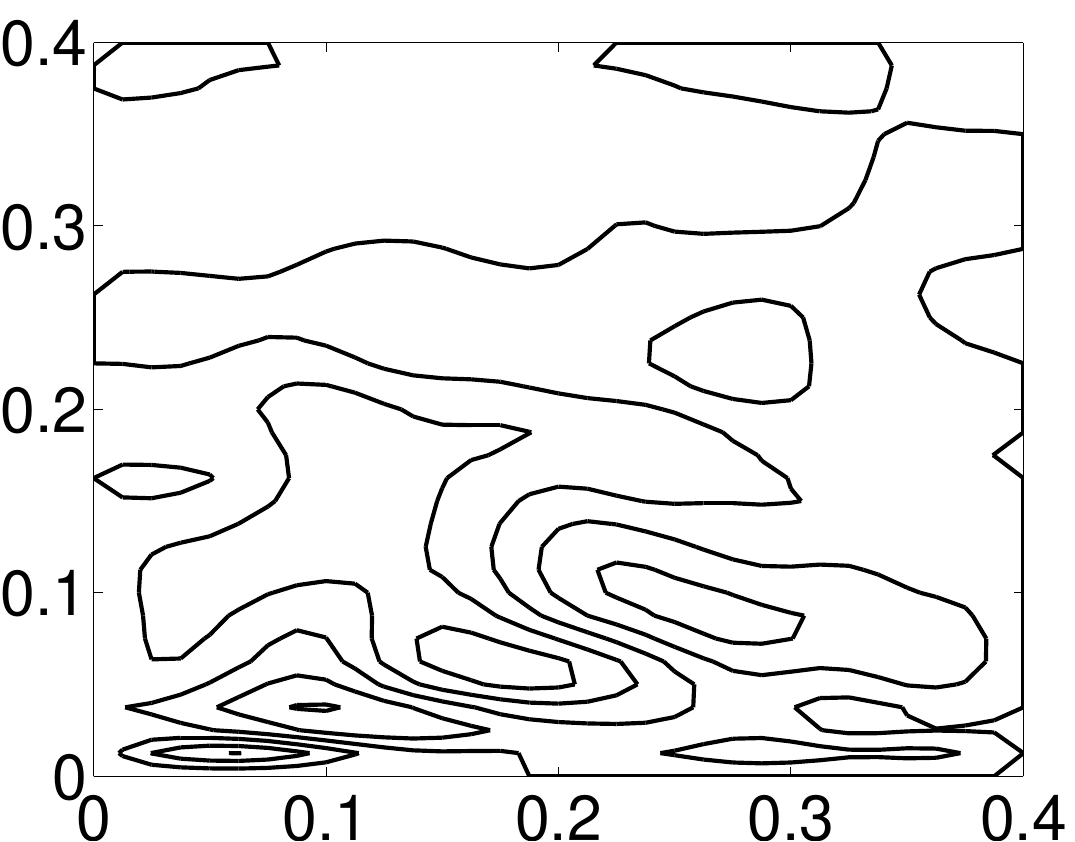} &\includegraphics[width=0.2\textwidth]{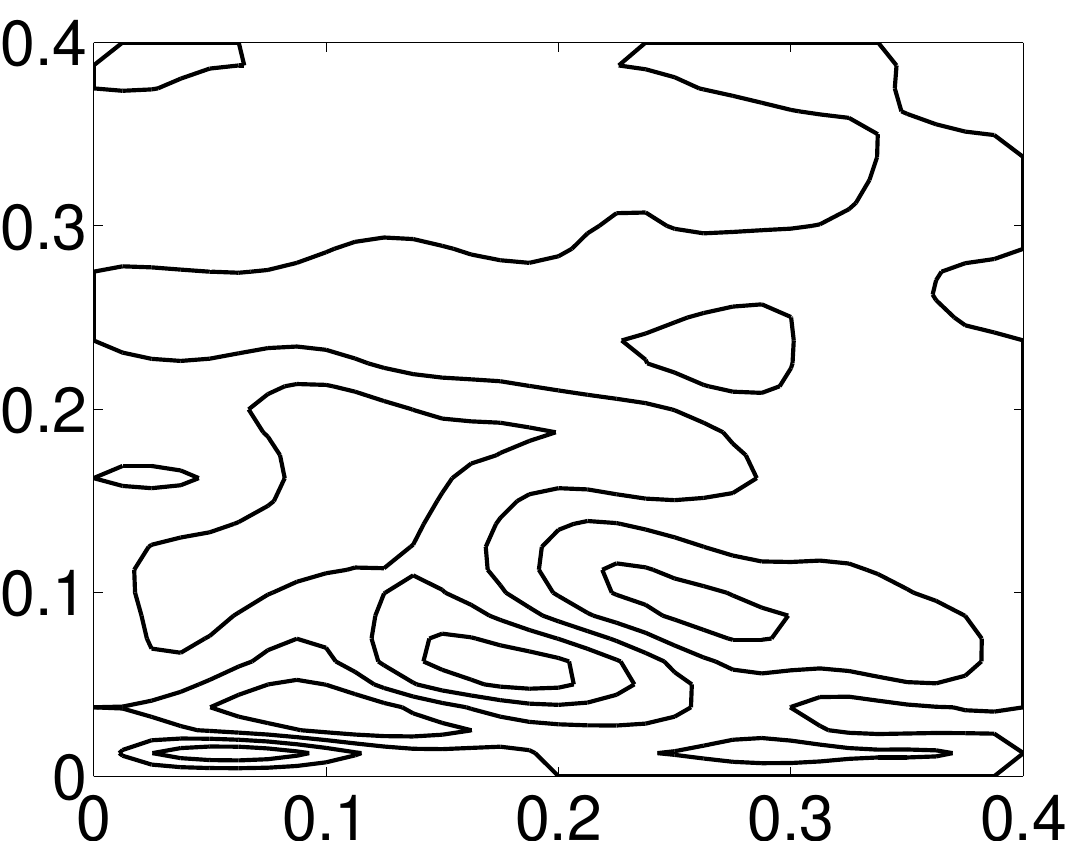} \\
\begin{sideways} ${\rm t = 1195}$ \end{sideways} &\includegraphics[width=0.2\textwidth]{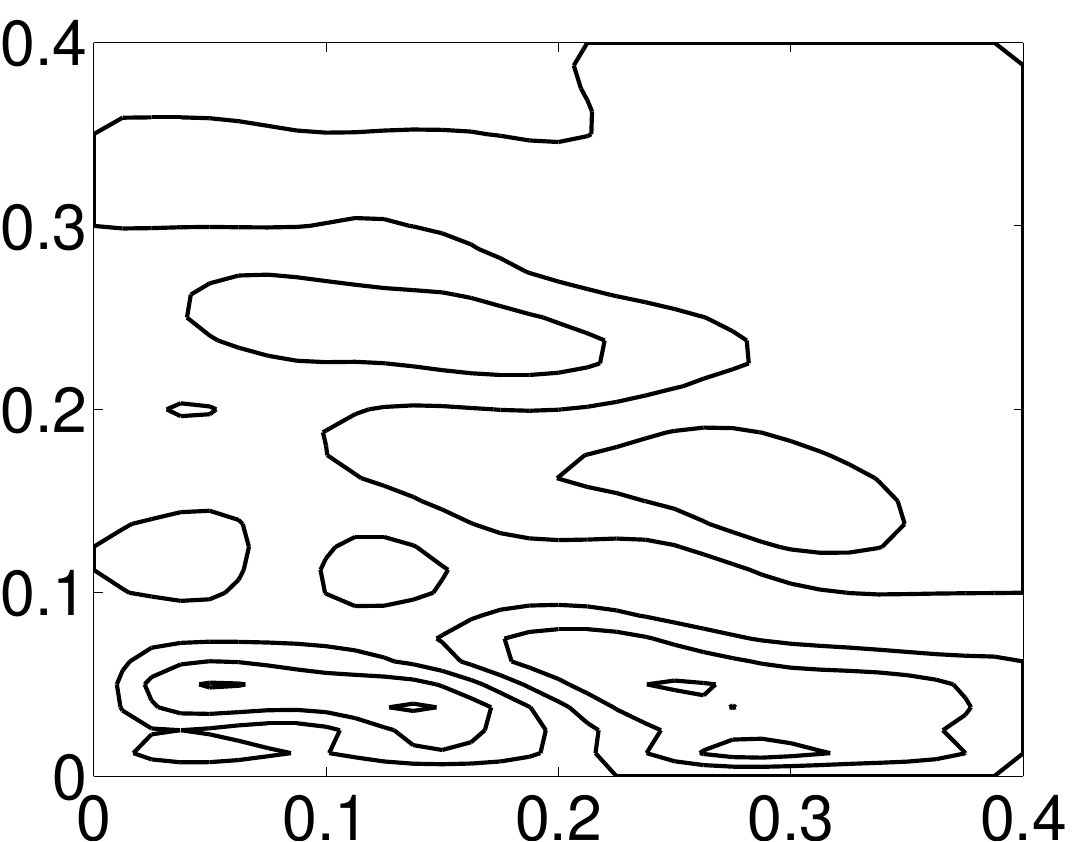} &\includegraphics[width=0.2\textwidth]{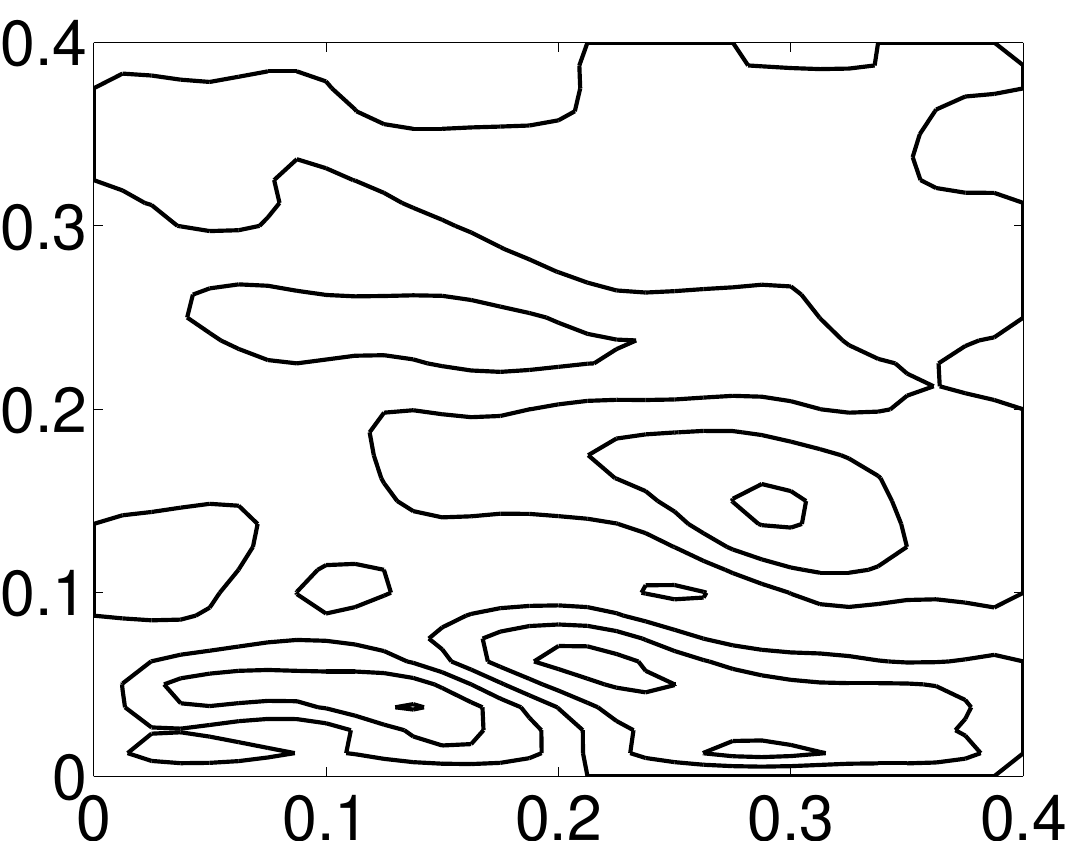} &\includegraphics[width=0.2\textwidth]{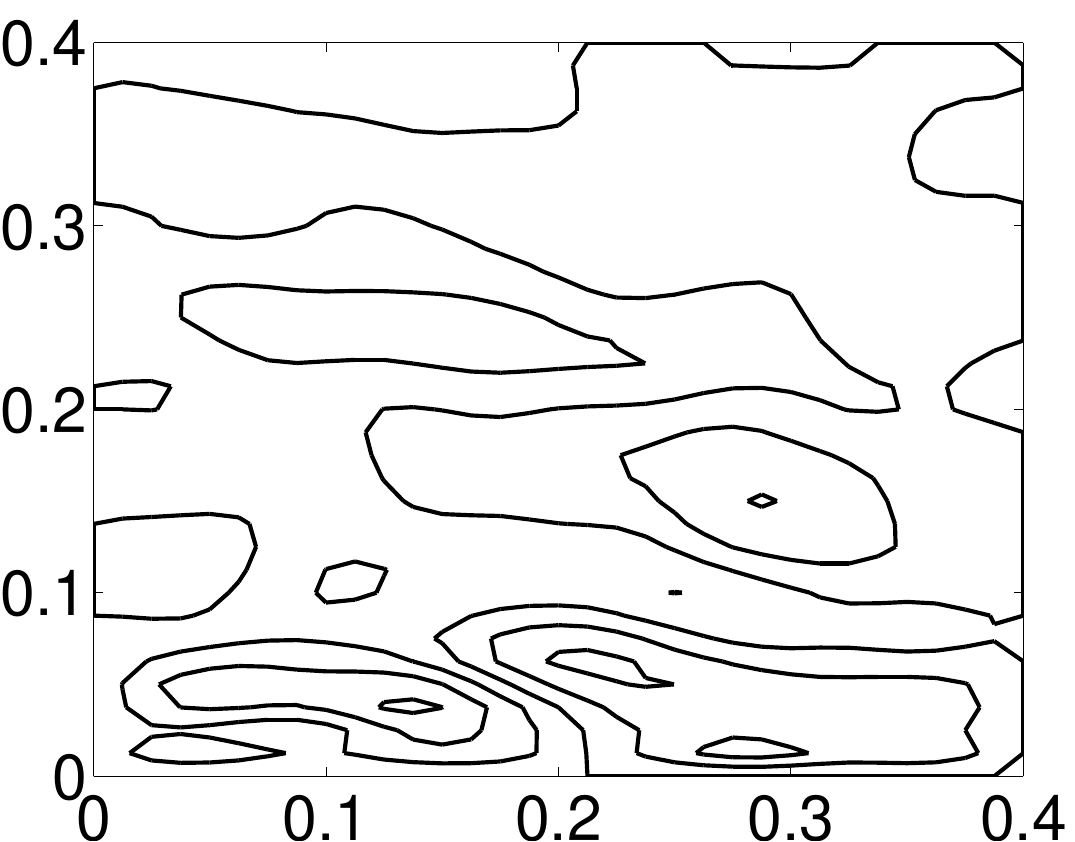} &\includegraphics[width=0.2\textwidth]{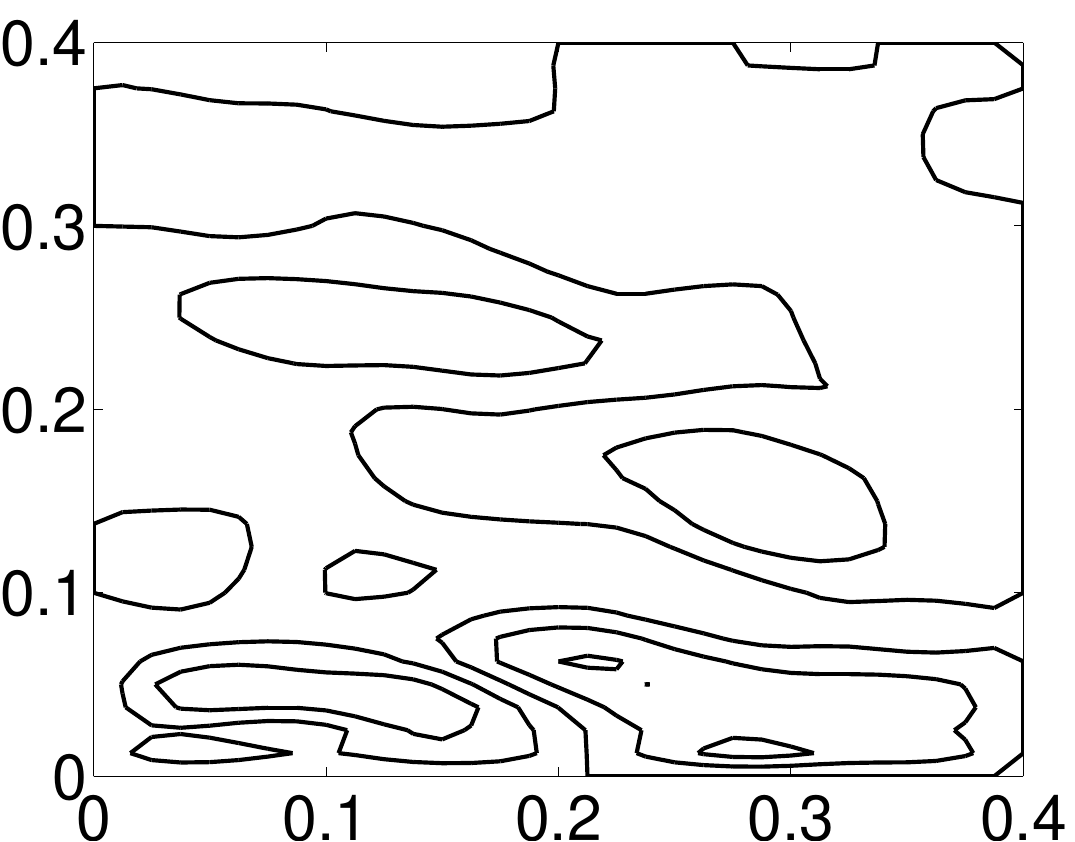} 
\end{tabular}
\caption{Snapshots of the QG33 simulation for $\Nens$ = 20,60 and 100 members, at the time steps $t=0,239,478,717,956$ and $1195$ (out of 1200). As expected, when the number of ensemble members is increased the estimation of the true state ($\xt$) is improved (the RMSE is decreased).
}
\label{Fig:QG33-Snapshots-Simulation}
\end{figure}

\begin{table}[H]
\centering
{\footnotesize
\begin{tabular}{|c|c|c|r|r|r|} \hline
$\bf N_{ens}$ & $\bf N_{obs}$ & $\bf STD_{ens}$ & ${\bf EnKF_{Sher}}$ & ${\bf EnKF_{Chol}}$ & ${\bf EnKF_{SVD}}$ \\ \hline
\multirow{9}{*}{20} &  \multirow{3}{*}{480}  & 2.5 & 17.6 s & 33.4 s & 24.7 s\\
& & 5.0 & 17.2 s & 32.8 s & 25.3 s \\ 
& & 7.5 & 17.3 s & 33.1 s & 28.6 s  \\ 
\cline{2-6}
&  \multirow{3}{*}{672}  & 2.5 & 17.9 s & 61.9 s & 39.4 s \\
& & 5.0 & 17.7 s & 62.9 s & 49.3 s \\ 
& & 7.5 & 17.8 s & 62.1 s & 37.6 s \\ 
\cline{2-6}
&  \multirow{3}{*}{864}  & 2.5 & 17.8 s & 113.7 s & 57.9 s \\
& & 5.0 & 18.2 s & 116.3 s & 79.7 s \\ 
& & 7.5 & 18.1 s & 118.4 s & 61.2 s \\ 
\cline{1-6}
\multirow{9}{*}{60} &  \multirow{3}{*}{480}  & 2.5 & 42.9 s & 57.8 s & 63.7 s \\
& & 5.0 & 42.9 s & 57.5 s & 62.9 s \\ 
& & 7.5 & 42.8 s & 57.8 s & 58.3 s \\ 
\cline{2-6}
&  \multirow{3}{*}{672}  & 2.5 & 44.3 s & 90.3 s & 142.3 s\\
& & 5.0 & 44.5 s & 89.8 s & 92.3 s \\ 
& & 7.5 & 44.5 s & 90.3 s & 91.3 s \\ 
\cline{2-6}
&  \multirow{3}{*}{864}  & 2.5 & 46.8 s & 150.6 s & 187.3 s \\
& & 5.0 & 46.6 s & 156.7 s & 144.8 s \\ 
& & 7.5 & 46.5 s & 154.3 s & 200.9 s \\ 
\cline{1-6}
\multirow{9}{*}{100} &  \multirow{3}{*}{480}  & 2.5 & 72.3 s & 83.2 s & 102.0 s \\
& & 5.0 & 72.4 s & 83.2 s & 96.8 s \\ 
& & 7.5 & 72.6 s & 83.4 s & 94.4 s \\ 
\cline{2-6}
&  \multirow{3}{*}{672}  & 2.5 & 77.8 s & 118.9 s & 140.2 s \\
& & 5.0 & 77.8 s & 119.5 s & 166.8 s \\ 
& & 7.5 & 77.1 s & 120.1 s & 216.4 s \\ 
\cline{2-6}
&  \multirow{3}{*}{864}  & 2.5 & 82.8 s & 209.1 s & 412.9 s \\
& & 5.0 & 83.2 s & 212.7 s & 289.7 s \\ 
& & 7.5 & 82.7 s & 202.6 s & 304.6 s \\ 
\hline
\end{tabular}
}
\caption{Computational times for several EnKF implementations applied to the QG33 instance.
Different numbers of ensemble members and numbers of observations are considered.}
\label{Tab:QG33-Results-ElapsedTime}
\end{table}
\begin{table}[H]
\centering
{\footnotesize
\begin{tabular}{|c|c|c|c|c|c|} \hline
$\bf N_{ens}$ & $\bf N_{obs}$ & $\bf STD_{ens}$ & ${\bf  EnKF_{Sher}  }$ & ${\bf EnKF_{Chol} }$ & ${\bf  EnKF_{SVD} }$ \\ \hline
\multirow{9}{*}{20} &  \multirow{3}{*}{1984}  & 2.5 & $2.38730154 \times 10^{-4}$ & $2.38730154 \times 10^{-4} $ & $ 2.38730154 \times 10^{-4} $ \\
& & 5.0 & $ 4.77366714 \times 10^{-4} $ & $ 4.77366714 \times 10^{-4} $ & $ 4.77366714 \times 10^{-4} $  \\ 
& & 7.5 & $ 7.16006328 \times 10^{-4} $ & $ 7.16006328 \times 10^{-4} $ & $ 7.16006328 \times 10^{-4} $  \\ 
\cline{2-6}
&  \multirow{3}{*}{2778}  & 2.5 & $ 2.38650829 \times 10^{-4} $ & $ 2.38650829 \times 10^{-4} $ & $ 2.38650829 \times 10^{-4} $ \\
& & 5.0 & $ 4.77195277 \times 10^{-4} $ & $ 4.77195277 \times 10^{-4} $ & $ 4.77195277 \times 10^{-4} $  \\ 
& & 7.5 & $ 7.15744082 \times 10^{-4} $ & $ 7.15744082 \times 10^{-4} $ & $ 7.15744082 \times 10^{-4} $  \\ 
\cline{2-6}
&  \multirow{3}{*}{3572}  & 2.5 & $ 2.38731447 \times 10^{-4} $ & $ 2.38731447 \times 10^{-4} $ & $ 2.38731447 \times 10^{-4} $ \\
& & 5.0 & $ 4.77279065 \times 10^{-4} $ & $ 4.77279065 \times 10^{-4} $ & $ 4.77279065 \times 10^{-4} $  \\ 
& & 7.5 & $ 7.15829932 \times 10^{-4} $ & $ 7.15829932 \times 10^{-4} $ & $ 7.15829932 \times 10^{-4} $  \\ 
\cline{1-6}
\multirow{9}{*}{60} &  \multirow{3}{*}{1984}  & 2.5 & $2.34788461 \times 10^{-4}$ & $2.34788461 \times 10^{-4} $ & $ 2.34788461 \times 10^{-4} $ \\
& & 5.0 & $ 4.69640125 \times 10^{-4} $ & $ 4.69640125 \times 10^{-4} $ & $ 4.69640125 \times 10^{-4} $  \\ 
& & 7.5 & $ 7.04505865 \times 10^{-4} $ & $ 7.04505865 \times 10^{-4} $ & $ 7.04505865 \times 10^{-4} $  \\ 
\cline{2-6}
&  \multirow{3}{*}{2778}  & 2.5 & $ 2.34427724 \times 10^{-4} $ & $ 2.34427724 \times 10^{-4} $ & $ 2.34427724 \times 10^{-4} $ \\
& & 5.0 & $ 4.68838293 \times 10^{-4} $ & $ 4.68838293 \times 10^{-4} $ & $ 4.68838293 \times 10^{-4} $  \\ 
& & 7.5 & $ 7.03254635 \times 10^{-4} $ & $ 7.03254635 \times 10^{-4} $ & $ 7.03254635 \times 10^{-4} $  \\ 
\cline{2-6}
&  \multirow{3}{*}{3572}  & 2.5 & $ 2.34303497 \times 10^{-4} $ & $ 2.34303497 \times 10^{-4} $ & $ 2.34303497 \times 10^{-4} $ \\
& & 5.0 & $ 4.68673565 \times 10^{-4} $ & $ 4.68673565 \times 10^{-4} $ & $ 4.68673565 \times 10^{-4} $  \\ 
& & 7.5 & $ 7.03046901 \times 10^{-4} $ & $ 7.03046901 \times 10^{-4} $ & $ 7.03046901 \times 10^{-4} $  \\ 
\cline{1-6}
\multirow{9}{*}{100} &  \multirow{3}{*}{1984}  & 2.5 & $2.37123911 \times 10^{-4}$ & $2.37123911 \times 10^{-4} $ & $ 2.37123911 \times 10^{-4} $ \\
& & 5.0 & $ 4.74051478 \times 10^{-4} $ & $ 4.74051478 \times 10^{-4} $ & $ 4.74051478 \times 10^{-4} $  \\ 
& & 7.5 & $ 7.10951271 \times 10^{-4} $ & $ 7.10951271 \times 10^{-4} $ & $ 7.10951271 \times 10^{-4} $  \\ 
\cline{2-6}
&  \multirow{3}{*}{2778}  & 2.5 & $ 2.34083117 \times 10^{-4} $ & $ 2.34083117 \times 10^{-4} $ & $ 2.34083117 \times 10^{-4} $ \\
& & 5.0 & $ 4.68257831 \times 10^{-4} $ & $ 4.68257831 \times 10^{-4} $ & $ 4.68257831 \times 10^{-4} $  \\ 
& & 7.5 & $ 7.02424460 \times 10^{-4} $ & $ 7.02424460 \times 10^{-4} $ & $ 7.02424460 \times 10^{-4} $  \\ 
\cline{2-6}
&  \multirow{3}{*}{3572}  & 2.5 & $ 2.32378928 \times 10^{-4} $ & $ 2.32378928 \times 10^{-4} $ & $ 2.32378928 \times 10^{-4} $ \\
& & 5.0 & $ 4.64697670 \times 10^{-4} $ & $ 4.64697670 \times 10^{-4} $ & $ 4.64697670 \times 10^{-4} $  \\ 
& & 7.5 & $ 6.97023351 \times 10^{-4} $ & $ 6.97023351 \times 10^{-4} $ & $ 6.97023351 \times 10^{-4} $  \\ 
\hline
\end{tabular}
}
\caption{Analysis RMSE for different EnKF implementations applied to the QG65 instance.
All methods give similar results. }
\label{Tab:QG65-Results-RMSE}
\end{table}


\begin{table}[H]
\centering
{\footnotesize
\begin{tabular}{|c|c|c|r|r|r|} \hline
$\bf N_{ens}$ & $\bf N_{obs}$ & $\bf STD_{ens}$ & $\bf {EnKF_{Sher}}$ & ${\bf EnKF_{Chol}}$ & ${\bf EnKF_{SVD}}$ \\ \hline
\multirow{9}{*}{20} &  \multirow{3}{*}{1984}  & 2.5 & 71.5 s &  44.9 min & 454.7 s\\
& & 5.0 & 70.6 s &  33.4 min & 429.3 s \\ 
& & 7.5 & 70.3 s &  45.4 min & 360.6 s\\ 
\cline{2-6}
&  \multirow{3}{*}{2778}  & 2.5 & 71.0 s & 1.8 h & 831.1 s\\
& & 5.0 & 73.2 s &  1.3 h & 735.7 s \\ 
& & 7.5 & 71.9 s &  1.4 h & 795.8 s \\ 
\cline{2-6}
&  \multirow{3}{*}{3572}  & 2.5 & 72.2 s & 3.8 h & 1602.5 s \\
& & 5.0 & 74.3 s &  3.3 h & 1112.3 s \\ 
& & 7.5 & 72.2 s &  3.0 h & 771.5 s \\ 
\cline{1-6}
\multirow{9}{*}{60} &  \multirow{3}{*}{1984}  & 2.5 & 179.0 s &  45.0 min & 1215.7 s \\
& & 5.0 & 179.5 s & 55.3 min & 1235.6 s \\ 
& & 7.5 & 178.4 s & 51.9 min & 1066.7 s \\ 
\cline{2-6}
&  \multirow{3}{*}{2778}  & 2.5 & 190.6 s &  2.4 h & 1463.3 s \\
& & 5.0 & 188.3 s & 1.9 h &   39.4 min\\ 
& & 7.5 & 189.7 s & 2.5 h &  32.4 min \\ 
\cline{2-6}
&  \multirow{3}{*}{3572}  & 2.5 & 202.6 s &  4.1 h & 1.2 h \\
& & 5.0 & 198.9 s & 2.9 h & 1.1 h  \\ 
& & 7.5 & 201.9 s &  4.3 h & 1.0 h  \\ 
\cline{1-6}
\multirow{9}{*}{100} &  \multirow{3}{*}{1984}  & 2.5 & 313.8 s & 52.4 min  & 37.2 min \\
& & 5.0 & 314.3 s & 40.1 min & 33.2 min \\ 
& & 7.5 & 309.6 s &  58.4 min  & 37.9 min  \\ 
\cline{2-6}
&  \multirow{3}{*}{2778}  & 2.5 & 346.9 s & 1.7 h &  1.1 h \\
& & 5.0 & 342.9 s &  1.7 h & 1.0 h  \\ 
& & 7.5 & 340.7 s &  2.9 h &  1.1 h  \\ 
\cline{2-6}
&  \multirow{3}{*}{3572}  & 2.5 & 373.9 s & 4.8 h & 1.6 h \\
& & 5.0 & 378.3 s & 4.0 h & 1.8 h \\ 
& & 7.5 & 383.7 s & 5.1 h  & 1.3 h \\ 
\hline
\end{tabular}
}
\caption{Computational times for several EnKF implementations applied to the QG65 instance.
Different numbers of ensemble members and numbers of observations are considered.}
\label{Tab:QG65-Results-ElapsedTime}
\end{table}

\begin{table}[H]
\centering
{\footnotesize
\begin{tabular}{|c|c|c|c|c|} \hline
$\bf N_{ens}$ & $\bf N_{obs}$ & $\bf STD_{ens}$ & ${\bf EnKF_{Sher}}$ & ${\bf  EnKF_{SVD}  }$ \\ \hline
\multirow{9}{*}{20} &  \multirow{3}{*}{8064}  & 2.5 & $9.92105438 \times 10^{-5}$ & $9.92105438 \times 10^{-5} $ \\
& & 5.0 & $ 1.98414661 \times 10^{-4} $ & $ 1.98414661 \times 10^{-4} $   \\ 
& & 7.5 & $ 2.97618537 \times 10^{-4} $ & $ 2.97618537 \times 10^{-4} $   \\ 
\cline{2-5}
&  \multirow{3}{*}{11290}  & 2.5 & $ 9.90045121 \times 10^{-5} $ & $ 9.90045121 \times 10^{-5} $  \\
& & 5.0 & $ 1.97996616 \times 10^{-4} $ & $ 1.97996616 \times 10^{-4} $   \\ 
& & 7.5 & $ 2.96988644 \times 10^{-4} $ & $ 2.96988644 \times 10^{-4} $   \\ 
\cline{2-5}
&  \multirow{3}{*}{14516}  & 2.5 & $ 9.87990386 \times 10^{-5} $ & $ 9.87990386 \times 10^{-5} $  \\
& & 5.0 & $ 1.97591049 \times 10^{-4} $ & $ 1.97591049 \times 10^{-4} $   \\ 
& & 7.5 & $ 2.96383002 \times 10^{-4} $ & $ 2.96383002 \times 10^{-4} $   \\ 
\cline{1-5}
\multirow{9}{*}{60} &  \multirow{3}{*}{8064}  & 2.5 & $9.74245572 \times 10^{-5}$ & $9.74245572 \times 10^{-5} $ \\
& & 5.0 & $ 1.94820830 \times 10^{-4} $ & $ 1.94820830 \times 10^{-4} $   \\ 
& & 7.5 & $ 2.92217500 \times 10^{-4} $ & $ 2.92217500 \times 10^{-4} $   \\ 
\cline{2-5}
&  \multirow{3}{*}{11290}  & 2.5 & $ 9.63593685 \times 10^{-5} $ & $ 9.63593685 \times 10^{-5} $  \\
& & 5.0 & $ 1.92682256 \times 10^{-4} $ & $ 1.92682256 \times 10^{-4} $   \\ 
& & 7.5 & $ 2.89005780 \times 10^{-4} $ & $ 2.89005780 \times 10^{-4} $   \\ 
\cline{2-5}
&  \multirow{3}{*}{14516}  & 2.5 & $ 9.67669396 \times 10^{-5} $ & $ 9.67669396 \times 10^{-5} $  \\
& & 5.0 & $ 1.93545013 \times 10^{-4} $ & $ 1.93545013 \times 10^{-4} $   \\ 
& & 7.5 & $ 2.90322987 \times 10^{-4} $ & $ 2.90322987 \times 10^{-4} $   \\ 
\cline{1-5}
\multirow{9}{*}{100} &  \multirow{3}{*}{8064}  & 2.5 & $9.56333807 \times 10^{-5}$ & $9.56333807 \times 10^{-5} $ \\
& & 5.0 & $ 1.91331921 \times 10^{-4} $ & $ 1.91331921 \times 10^{-4} $   \\ 
& & 7.5 & $ 2.87027307 \times 10^{-4} $ & $ 2.87027307 \times 10^{-4} $   \\ 
\cline{2-5}
&  \multirow{3}{*}{11290}  & 2.5 & $ 9.49419202 \times 10^{-5} $ & $ 9.49419202 \times 10^{-5} $  \\
& & 5.0 & $ 1.89929524 \times 10^{-4} $ & $ 1.89929524 \times 10^{-4} $   \\ 
& & 7.5 & $ 2.84918006 \times 10^{-4} $ & $ 2.84918006 \times 10^{-4} $   \\ 
\cline{2-5}
&  \multirow{3}{*}{14516}  & 2.5 & $ 9.47165868 \times 10^{-5} $ & $ 9.47165868 \times 10^{-5} $  \\
& & 5.0 & $ 1.89427095 \times 10^{-4} $ & $ 1.89427095 \times 10^{-4} $   \\ 
& & 7.5 & $ 2.84137686 \times 10^{-4} $ & $ 2.84137686 \times 10^{-4} $   \\ 
\hline
\end{tabular}
}
\caption{Analysis RMSE for different EnKF implementations applied to the QG129 instance.
All methods give similar results.}
\label{Tab:QG129-Results-RMSE}
\end{table}

\begin{table}[H]
\centering
{\footnotesize
\begin{tabular}{|c|c|c|c|c|} \hline
$\bf N_{ens}$ & $\bf N_{obs}$ & $\bf STD_{ens}$ & ${\bf  EnKF_{Sher}  }$ & ${\bf EnKF_{SVD}  }$ \\ \hline

\multirow{9}{*}{20} &  \multirow{3}{*}{8064}  & 5 & 289.9 s & 1.8 h \\
& & 5.0 & 293.5 s & 1.4 h \\ 
& & 7.5 & 286.8 s &  1.9 h  \\ 
\cline{2-5}
&  \multirow{3}{*}{11290}  & 5 & 303.2 s & 2.2 h  \\
& & 5.0 & 302.2 s & 2.0 h   \\ 
& & 7.5 & 303.5 s &  2.9  h  \\ 
\cline{2-5}
&  \multirow{3}{*}{14516}  & 5 & 315.7 s & 4.5 h  \\
& & 5.0 & 308.1 s & 3.8 h  \\ 
& & 7.5 & 309.1 s &   3.7 h  \\ 
\cline{1-5}
\multirow{9}{*}{60} &  \multirow{3}{*}{8064}  & 5 & 764.8 s & 4.1 h \\
& & 5.0 & 795.9 s &  3.2 h   \\ 
& & 7.5 & 764.7 s & 3.9 h   \\ 
\cline{2-5}
&  \multirow{3}{*}{11290}  & 5 & 838.0 s &  6.3 h \\
& & 5.0 & 832.6 s & 5.4 h  \\ 
& & 7.5 & 817.7 s & 6.7 h   \\ 
\cline{2-5}
&  \multirow{3}{*}{14516}  & 5 & 910.7 s &  9.9 h  \\
& & 5.0 & 899.9 s &  9.3 h   \\ 
& & 7.5 & 864.8 s &  11.2 h   \\ 
\cline{1-5}
\multirow{9}{*}{100} &  \multirow{3}{*}{8064}  & 5 & 1381.5 s & 5.3 h \\
& & 5.0 & 1360.3 s & 5.2 h  \\ 
& & 7.5 & 1397.9 s & 5.9 h   \\ 
\cline{2-5}
&  \multirow{3}{*}{11290}  & 5 & 1492.4 s & 10.3 h  \\
& & 5.0 & 1494.5 s &  8.9 h   \\ 
& & 7.5 & 1506.6 s & 6.7 h   \\ 
\cline{2-5}
&  \multirow{3}{*}{14516}  & 5 & 1624.8 s & 13.7 h  \\
& & 5.0 & 1634.9 s & 12.9 h   \\ 
& & 7.5 & 1664.2 s & 14.9 h   \\ 
\hline
\end{tabular}
}
\caption{Computational times for several EnKF implementations applied to the QG129 instance.
Different numbers of ensemble members and numbers of observations are considered.}
\label{Tab:QG129-Results-ElapsedTime}
\end{table}
\section{Conclusions and Future Work}\label{sec:conclusions}

We propose a novel implementation of the EnKF based on an iterative application of the Sherman-Morrison formula. The algorithm exploits the special structure of the background error covariance matrix projected onto the observation space. The computational complexity of the new approach is equivalent to that of the best EnKF implementations  available in the literature. Nevertheless, the performance (elapsed time) of most of the existing methods is strongly dependent from the condition $\Nobs \gg \Nens$ (the number of observations is large compared to the ensemble size); the performance of the new approach is not affected by this condition. In addition, the term $\Nens^3$ is not presented in the upper-bound of the effort of the proposed method, which leads to better performance when the number of observations and of ensemble members are of similar magnitude ($\Nobs \sim \Nens$). A sufficient condition for the stability of the proposed method is the non-singularity of the data error covariance matrix, which, in practice, is always the case. In addition, a pivoting strategy is developed in order to reduce round-off error propagation without increasing the computational effort of the proposed method. The computational cost of this algorithm provides a better theoretical performance than other generic formulations of matrix inversion based on the Sherman Morrison formula available in the literature. To assess the accuracy and performance of the proposed implementation two standard test problems have been employed, namely, the Lorenz 96 and the quasi-geostrophic models. All EnKF implementations tested (Cholesky, SVD, Sherman-Morrison) provide virtually identical analyses. However, the proposed Sherman-Morrison approach is much faster than the others even when the number of observations is large with respect to the number of ensemble members ($\Nobs \gg \Nens$). The parallel version of the new algorithm has a theoretical complexity that grows only linearly with the number of observations, and is therefore well suited for implementation in large scale data assimilation systems.

\section*{Acknowledgment}
This work has been supported in part by NSF through awards NSF OCI-8670904397, NSF CCF-0916493, NSF DMS-0915047, NSF CMMI-1130667, NSF CCF-1218454, AFOSR FA9550-12-1-0293-DEF, AFOSR 12-2640-06, and by the Computational Science Laboratory at Virginia Tech.

\bibliographystyle{abbrv}

\end{document}
This is never printed